%% file: main.tex
\documentclass[final,3p,nopreprintline]{elsarticle}

\def\vdate{March 16, 2025} %

\usepackage{lineno}
\usepackage[colorlinks=true,breaklinks=true,pdftex]{hyperref}
\usepackage{booktabs}
\usepackage{makecell} 

\modulolinenumbers[1]

\usepackage{amsfonts,amsthm,bbm,amssymb,epsf,amsmath,%
latexsym,amscd,amsfonts,enumerate,amsthm,supertabular,color,url}

\usepackage{ulem} 
\usepackage{subfigure} 
\usepackage[nottoc]{tocbibind}

\newdefinition{remark}{Remark}
\newdefinition{algorithm}{Algorithm}

\definecolor{orange}{RGB}{255,127,0}

\biboptions{authoryear}
\def\bxi{\mbox{\boldmath$\xi$}}
\def\RR{{\mathbb{R}}}
\def\ZZ{{\mathbb{Z}}}
\def\NN{{\mathbb{N}}}

\def \bfv {{\bf v}}
\def \z {{\bf z}}
\def \x {{\bf x}}
\def \y {{\bf y}}
\def \w {{\bf w}}
\def \v {{\bf v}}

\def \F {{\bf F}}

\def\boeta{{\boldsymbol\eta}}
\def\bxi{{\boldsymbol\xi}}
\DeclareMathOperator{\spann}{span}
\DeclareMathOperator{\median}{median}

\def\blue#1{{#1}}

\begin{document}
\newcommand{\prf}{\noindent{{\bf Proof} :\ }}
\newcommand{\QED}{\vrule height 1.4ex width 1.0ex depth -.1ex\ \medskip}
\renewcommand{\thetable}{\arabic{table}}

\begin{frontmatter}

\title{
A positive meshless finite difference scheme for scalar conservation laws 
with adaptive artificial viscosity driven by fault detection
}
\author[first]{Cesare~Bracco}%
\ead{cesare.bracco@unifi.it}
\author[last]{Oleg Davydov}
\ead{oleg.davydov@math.uni-giessen.de}
\author[first]{Carlotta Giannelli}
\ead{carlotta.giannelli@unifi.it}
\author[first]{Alessandra Sestini}
\ead{alessandra.sestini@unifi.it}
\address[first]{
Dipartimento di Matematica e Informatica ``U.~Dini'', 
Universit\`a degli Studi di Firenze,\\
Viale Morgagni 67a, 50134 Florence, Italy
}
\address[last]{
Department of Mathematics, Justus Liebig University Giessen,
Arndtstrasse 2, 35392 Giessen, Germany\\[12pt] \vdate\\[-24pt]}

\begin{abstract}
We present a meshless finite difference method for multivariate scalar conservation laws that generates positive schemes 
satisfying a local maximum principle on irregular nodes and relies on artificial viscosity for shock capturing. Coupling 
 two different numerical differentiation formulas and the adaptive selection of the sets of influence allows to meet a local CFL
condition without any {\it a priori}\  time step restriction. The artificial viscosity term is chosen in an adaptive way
by applying it only in the vicinity of the sharp features of the solution identified by an algorithm for fault detection 
on scattered data. Numerical tests demonstrate a robust performance of the method on irregular nodes and advantages of
adaptive artificial viscosity. The accuracy of the obtained solutions is comparable to that for standard
monotone methods available only on Cartesian grids.

\end{abstract}


\end{frontmatter}

\section{Introduction}
\input{section01.tex}

\section{Positive meshless finite difference schemes for scalar conservation laws}
\label{sec:ed}
\input{section02.tex}

\section{Minimal numerical differentiation formulas}
\label{sec:ndf}

\input{section03.tex}

\section{Adaptive fault-based viscosity}
\label{sec:df}

\input{section04.tex}
\section{Summary of algorithms}
\label{sec:algorithm}
\input{sectionALG.tex}

\section{Numerical examples} 
\label{sec:ne}

\input{section05.tex}

\section{Conclusion}
\label{sec:c}

We propose a positive scheme for multivariate scalar conservation laws based on a purely meshless finite difference
discretization in space, explicit Euler time stepping and locally varying artificial viscosity.  To obtain the scheme, we
modify the  $\ell_2$-minimal numerical differentiation formulas of  \cite{davydov2018} by inequality constraints, and apply
them  to  the directional derivatives arising from the quasilinear form of the conservation law as well as to the Laplace
operator that defines artificial viscosity.  Coupling  of the two numerical differentiation formulas and adaptive selection of
the sets of influence allows to meet a local CFL condition without any {\it a priori} time step restriction. Moreover, we
suggest adaptive selection of the amount of local artificial viscosity by applying the fault detection method of
\cite{bracco2019}. This is possible in particular because the positivity and hence stability of  the scheme is achieved
without relying on the artificial viscosity term, whose sole purpose is therefore the shock capturing. 

Numerical experiments for bivariate benchmark test problems confirm the robust performance of the method on irregular nodes
and advantages of   adaptive artificial viscosity. The accuracy of the obtained solutions is comparable to if not better than that for
standard monotone methods available only on Cartesian grids.

There are many directions of research of interest for future work. Purely polynomial numerical differentiation may be
enhanced by kernel-based techniques that are used successfully in the meshless finite difference methods of RBF-FD type, see 
\cite{fornberg2015primer,DavydovSchaback16,bayona2017role,davydov2019optimal}.
The flexibility of the space discretization makes it natural to explore in this setting the local time-stepping combined with
local refinement in space. Positive schemes may be used locally where appropriate, enhanced by high order schemes as part of 
limiter- or ENO/WENO-based algorithms. In particular, an extension to the conservation laws 
of the meshless FCT method for linear convection dominated problems \cite{sokolov2019flux} could be considered.
\blue{An important goal for future work is  
the extension of the introduced positive schemes to general vector-valued conservation problems.}

\section*{Acknowledgments} We thank anonymous reviewers whose comments helped to improve the quality of the paper.
C. Bracco, C. Giannelli and A. Sestini are members of the INdAM Research Group GNCS. This work has been partially supported by
INdAM through Finanziamenti Premiali SUNRISE and GNCS. 

\bibliography{biblio}

\end{document}

%% file: section01.tex
Meshless methods discretize partial differential equations on irregular nodes in space without creating a grid or mesh
\cite{SJD22}. This allows to distribute the degrees of freedom in the most flexible way and rapidly rearrange them if needed
as time evolves, without meshing or re-meshing, and thus eliminating a major time consuming and often algorithmically difficult
step in numerical simulations. Especially, meshless finite difference methods, based on polynomials or radial basis functions,
see for example \cite{onate1996finite,fornberg2015primer,JTAMB19}, attract increasing attention of engineers and mathematicians as a
general purpose numerical technique that combines the basic simplicity and high accuracy inherited from the classical finite
difference method with full spatial flexibility.

Numerical approximation of nonlinear multivariate hyperbolic conservation laws that model compressible fluid flow is an area
where these methods have a lot to offer because of the complexity of the solutions, with shock discontinuities developing and
moving in time. A major challenge for developing effective  schemes of this type is that the standard multivariate finite
difference schemes for conservation laws are obtained by applying univariate methods in a  dimension-by-dimension fashion, see
e.g.\ \cite{hestbook}. %
Univariate numerical fluxes or Riemann solvers may also be naturally applied in the normal direction to the element interfaces 
in partition-based methods. 
However, exploiting such univariate structures is more difficult in the setting of meshless finite difference methods, and may
lead to undesired restrictions on the node generation or selection of sets of influence, reducing the natural geometric
flexibility of these methods. 
Nevertheless, several approaches exploiting an edge-based connectivity structure by either borrowing it
from a mesh or generating it in a meshless cloud of nodes, can be found in 
\cite{lohner2002finite,sridar2003upwind,shu2005upwind,praveen2007kinetic,ortega2009finite,
katz2009comparison,katz2010meshless,huh2018new}.

The idea that the stability of a numerical scheme in the vicinity of the shocks may be achieved with
the help of an appropriate  artificial viscosity term  added to a central type scheme for the divergence goes back to
\cite{vonneumann1950method} for classical finite differences, and is a popular approach  
in the meshless finite difference methods,  see in particular
\cite{batina1993gridless,fornberg2011stabilization,flyer2012guide,shankar2018hyperviscosity,tominec2023residual}.
The viscosity term serves to enforce a stable eigenvalue spectrum of the discrete iterative scheme and also helps to correctly
reproduce the shock movement thanks to the principle of vanishing viscosity. 
\cite{ghosh1995least,kuhnert2003upwind,li2017generalized,seifarth2017numerische} pursue the same goals via upwind schemes, 
where the divergence at a node  is discretized on neighboring nodes in the upstream half-space.
In particular, for scalar conservation laws the artificial viscosity term or upwinding may be organized such that
the weights of the scheme are \textit{positive} if the neighborhoods of all nodes are well-balanced. %
Then the numerical approximation of the state variable at each 
node is a convex combination of the state variables of the previous time step at neighboring nodes if the time step is
sufficiently small, which ensures local maximum principle and hence the stability. 
Note that in the meshless setting numerical results in the literature are often presented for nodes obtained  by highly
optimized algorithms based for example on the minimization of the repulsion energy applied to  simple domains such as the 
square. 
On the other hand, the full advantage of meshless methods may only be realized if the methods are sufficiently robust to
guarantee stability on rather irregular nodes obtained by inexpensive algorithms, see \cite{SJD22}.

Numerical schemes for bivariate scalar conservation laws that are guaranteed to be positive on irregular nodes  have been
suggested in \cite{FS01} and \cite{yin2011five}. In this approach the individual weights of a numerical  differentiation
formula for the divergence on irregular nodes are first modified to eliminate the unwanted signs, and then an artificial
viscosity term is added to restore the first order consistency of the formula. Therefore the amount of local viscosity is
driven by the stability goal, and only by chance may help the shock capturing.

In this paper we suggest a positive meshless scheme on irregular nodes for scalar conservation laws with any number of 
independent variables that does not require any additional structures in the node sets, and employs artificial viscosity for
shock capturing rather than enforcing positivity. The conservation law is discretized in the quasilinear formulation via
numerical differentiation formulas for the gradient operator on irregular nodes by a polynomial local  least squares scheme
combining minimization of an error bound for numerical differentiation suggested in \cite{davydov2018} with inequality
constraints that ensure the positivity of the scheme. The solvability of the small local quadratic optimization problems
arising from this approach is achieved by consecutive enlargements of the set of influence  until it meets a local CFL
condition (that is, the positivity of the central weight of the scheme). Therefore, no global CFL condition is needed, and no
\textit{a priori} time step restriction is used. Artificial viscosity term is chosen in a adaptive way  by applying 
it only in the vicinity of the sharp features of the solution, in order to ensure shock capturing, which  reduces the overall
numerical dissipation. We identify the regions of discontinuity by the recent fault detection method of \cite{bracco2019}
designed for  scattered data.

We test numerically three versions of our scheme (no artificial viscosity, constant or adaptive artificial viscosity) on
irregular nodes for three bivariate benchmark problems: inviscid Burgers' equation with either smooth or discontinuous
initial condition, and the rotating wave problem with a non-convex flux. In all cases solutions are non-oscillating as
expected due to the positivity of the schemes. The shocks are reproduced by both methods that employ artificial viscosity,
while the version with adaptive viscosity is less dissipative. The accuracy of our solutions is comparable to those obtained
by classical monotone schemes on Cartesian grids, and mostly somewhat better, despite using irregular nodes.

We expect that our positive first order scheme may be complemented by high order discretizations in conjuction with 
limiters, ENO/WENO and further ideas known in the literature on conservation laws.

The structure of the paper is as follows. Section~\ref{sec:ed} introduces positive meshless finite difference schemes and
discusses their stability, consistency and the role of the artificial viscosity term. Section~\ref{sec:ndf} describes our
algorithm for choosing the weights of the inequality constrained minimal  numerical differentiation formulas for the space
discretization of the conservation law and a viscosity term, which serves as the main tool to ensure a local CFL  condition
and the positivity of the scheme. Section~\ref{sec:df} presents our adaptive algorithm that determines the amount of local 
artificial viscosity for each node. In Section~\ref{sec:algorithm} we formulate three methods to obtain positive schemes that
differ in how they handle the artificial viscosity. Numerical tests are presented in Section~\ref{sec:ne}, and a conclusion is
given in Section~\ref{sec:c}.

%% file: section02.tex
In this paper we deal with the meshless numerical solution of a general scalar conservation law on a bounded domain
$\Omega \subset \RR^d$, $d\ge 2$. (However, as usual, a problem on a hyperrectangle  with
periodic boundary conditions may be interpreted as being posed on $\Omega = \RR^d$.) 
Thus, the considered problem is modeled by the following PDE
\begin{equation} \label{Claw}
u_t + \nabla_\x \cdot \F (u) = 0 \,, \qquad t \in (0\,,\,T)\,, \quad \x = (x_1,\ldots,x_d) \in \Omega\,, 
\end{equation}
where $\nabla_\x \cdot$ denotes the divergence operator in space and the PDE is endowed with suitable boundary conditions and with
the initial condition
\begin{equation} \label{initcond}
u(0,\x) = u_0(\x)\,, \quad \x \in \Omega\,.
\end{equation}
As usual, the unknown function 
$u: (0\,,\,T) \times \Omega \rightarrow \RR,\,\, u =u(t\,,\,\x),$ is called the {\it state variable}, 
and the  given vector function 
$\F : \RR \rightarrow \RR^d$, $\F(u) =(F_1(u),\ldots,F_d(u))$, with $F_i:\RR \rightarrow \RR$, $i=1,\ldots,d$, is the {\it flux}. 
We assume that $\F$ is differentiable, so that the {\it conservative form}  (\ref{Claw}) of the equation  may be replaced
by the equivalent \textit{quasilinear} equation
\begin{equation} 
\label{Claweq} u_t + \F' (u)  \nabla_\x u = 0,
\end{equation}
where  $\nabla_\x$ is the space gradient and $\F'(u) :=(F'_1(u),\ldots,F'_d(u))$.
In particular, in the case when $\F(u) = u \bfv$, with $\bfv$ denoting a nonvanishing constant vector in $\RR^d$,
we obtain the equation of the  {\it constant linear transport}. For any conservation law (\ref{Claw}), 
the state variable $u$  
is conserved with respect to time. In particular, for each $\x_0 \in \Omega, \, u$ remains constantly equal to $u_0(\x_0)$
along  the corresponding outgoing {\it characteristic} line $(t, \x(t))$, $t\ge0$, of the $d+1$-dimensional $t$-$\x$ space,
where
$$
\x(t) = t \F'(\hat u_0) + \x_0\,, \quad \hat u_0 = u_0(\x_0)\,.$$
Clearly, the presence of characteristics implies first of all that  every finite jump in $u_0$ is transported along the
characteristics when time goes on. Furthermore, in the nonlinear case two characteristics may intersect and this implies that
the strong solution of (\ref{Claw}) is only well defined for sufficiently small $T$. For later times, discontinuities (shocks)
may appear, even when $u_0$ is smooth, after which the equation gets many weak solutions, such that special efforts are
required  to identify a physically meaningful continuation of the state variable, so called {\it entropy solution},  
see for example \cite{hestbook}.

Because of all these difficulties, sophisticated adaptive algorithms are needed in order to obtain satisfactory numerical
approximations of challenging conservation law problems arising in applications. Therefore numerical methods capable to
address these challenges must be in position to quickly refine, coarsen and rearrange the discretization of the spatial
computational domain when time evolves. Even if mesh-based numerical techniques such as finite volume method or discontinuous 
Galerkin method possess a great deal of flexibility, they are often slowed down by the pitfalls of the meshing technology. 
Meshless methods, in particular meshless finite differences, that rely solely on nodes distributed irregularly in
the space, without the need to connect them into meshes, provide a promising alternative.  
 
A variety of numerical algorithms for conservation laws has been developed for mesh-based methods, in particular for  the classical finite difference method. They include low order monotone methods, high order limiter-based schemes, as well as high order ENO and WENO type approaches. Typically methods based on finite differences are  developed for the univariate case and then generalized to Cartesian grids by a dimension-by-dimension construction, which makes mesh adaptation difficult.  

In this paper we suggest a genuinely multivariate positive  scheme 
applicable on arbitrary irregular nodes. It relies on the meshless finite difference method and therefore
requires no  grids, triangulations or any other meshes connecting the nodes.

Let $X=\{\x_1,\ldots,\x_N\}\subset\overline{\Omega}$ be an arbitrary finite set of discretization nodes, and $\Delta t>0$ a step 
size in time. Note that the term $\F' (u)  \nabla_\x u$ in \eqref{Claweq} for each fixed $\x$ is the directional derivative 
of $u$ at $\x$,
$$
D_{\boeta}u(t,\x)=\frac{\partial u}{\partial \boeta}(t,\x)=\eta_1\frac{\partial u}{\partial x_1}(t,\x)
+\cdots+\eta_d\frac{\partial u}{\partial x_d}(t,\x)$$ 
in the direction $\boeta=\F' (u(t,\x))$, a vector with components $F_1'(u(t,\x)),\ldots,F_d'(u(t,\x))$. 
(Note that $\boeta$ is not normalized.) 
The approximate solution $U$ is computed recursively by applying the Euler method in time and bivariate numerical
differentiation formulas in space according to the rule 
\begin{equation} \label{schm}
U(t+\Delta t,\x_i)=U(t,\x_i) - \Delta t \sum_{\x_j\in X_i}w_{ij}U(t,\x_j),\qquad U(0,\x_i)=u_0(\x_i),%
\end{equation}
where $X_i$ is a subset of $X$ called  the {\it influence set} of the node $\x_i$, and $w_{ij}\in\RR$ are the 
{\it weights} of a numerical differentiation formula for the directional
derivative at $\x_i$,
\begin{equation} \label{ndfeta}
D_{\boeta_i}f(\x_i)\approx \sum_{\x_j\in X_i}w_{ij} f(\x_j),\quad 
\quad \boeta_i:=\F' (U(t,\x_i)).
\end{equation}
The influence set $X_i$ is a small selection of nodes in the neighborhood of $\x_i$.
We stress that in contrast to the classical finite difference method, the influence sets $X_i$ are irregular and 
usually have different geometric shapes for different nodes $\x_i$. 
Hence, the weights $w_{ij}$ have to be determined individually for each
$\x_i$ rather than being some scaled versions of a fixed stencil. Moreover, the weights also depend on the direction 
$\boeta_i$, therefore the weights  $w_{ij}$, $j\in X_i$, are in general different for different $i$ and different times.

Assuming without loss of generality that $\x_i\in X_i$ 
(otherwise we add $\x_i$ to $X_i$ and set $w_{ii}=0$), we write \eqref{schm} in the form
$$
U(t+\Delta t,\x_i)=(1- \Delta t\,w_{ii})U(t,\x_i) - \Delta t \sum_{\x_j\in X_i\setminus\{\x_i\}}w_{ij}U(t,\x_j).$$
The scheme \eqref{schm} is said to be {\it positive} if the right hand side expression is a positive linear combination
of $U(t,\x_j)$, $\x_j\in X_i$, that is, if
\begin{equation} \label{m1}
\Delta t\,w_{ii} \le 1
\end{equation} 
and 
\begin{equation}
\label{m2}
w_{ij} \le 0,\quad \forall\x_j\in X_i\setminus\{\x_i\}.
\end{equation} 
Note that the weights $w_{ij}$ may depend on $U(t,\x_i)$ since the direction $\boeta_i$ in \eqref{ndfeta} in general depends
on it. Therefore, a positive scheme is not necessarily monotone, for which the right hand side of \eqref{schm} would be
required to be non-decreasing with respect to the values $U(t,\x_i)$ and $U(t,\x_j)$, $\x_j\in X_i$, of the approximate state
variable at the previous time step, see \cite{hestbook}. Under additional assumption
\begin{equation} \label{sumw}
\sum_{\x_j\in X_i}w_{ij} =0
\end{equation}
the approximated value $U(t+\Delta t,\x_i)$ of the state variable at time  $t+\Delta t$ is a convex combination of the values
$U(t,\x_j)$, $\x_j\in X_i\cup \{\x_i\}$. %
Therefore the solution satisfies a {\it local maximum principle} in the form
\begin{equation}\label{locmax}
\min_{\x_j\in X_i\cup \{\x_i\}}U(t,\x_j)\le U(t+\Delta t,\x_i)\le \max_{\x_j\in X_i\cup \{\x_i\}}U(t,\x_j),
\end{equation}
which in particular implies the following inequalities for all time steps and all nodes $\x_i$
$$
\min_{\x\in \overline{\Omega}}  u_0(\x)\le \min_{x_j\in X} U(0,\x_j)\le U(t,\x_i)\le \max_{x_j\in X} U(0,\x_j)
\le \max_{\x\in \overline{\Omega}}  u_0(\x),$$
reproducing the maximum principle 
$$
\min_{\bxi\in \overline{\Omega}}  u_0(\bxi)\le u(t,\x)\le \max_{\bxi\in \overline{\Omega}}  u_0(\bxi)$$
available for the entropy solution of the scalar conservation law \eqref{Claw}, see for example  \cite{ZhangXiaShu12}. 
The local maximum principle \eqref{locmax} ensures the {\it stability} of the scheme and prevents  unphysical behavior, such as
artificial oscillations.

Positive schemes for scalar conservation laws originated in the classical upwind method by \cite{CIR52}. Positive and related 
\textit{local extremum diminishing} schemes of higher order have been developed for mesh-based methods, 
see for example \cite{jameson1995analysis,guermond2014second,kuzmin2020monolithic}, as well as for meshless finite difference
methods, using an underlying mesh or an edge-based connectivity  structure \cite{katz2009comparison,katz2010meshless}. 
Without imposing any connectivity of meshless discretization nodes, or their well-balanced distribution,
positive first order methods have been suggested by
\cite{FS01,yin2011five} in the bivariate case $d=2$. In contrast to the latter, our method described below 
 is not restricted to two space dimensions, and does not use artificial viscosity to achieve positivity,
which allows to devote it to the shock capturing.

We say that the scheme \eqref{schm} is {\it consistent} in space, if  the weights $w_{ij}$ are chosen such that
\eqref{ndfeta} is exact for all constant and linear polynomials, 
\begin{equation} \label{exa1}
D_{\boeta_i}p(\x_i)= \sum_{\x_j\in X_i}w_{ij} p(\x_j)\qquad \forall p\in \Pi^d_2,
\end{equation}
where $\Pi^d_q$ denotes the linear space of all  polynomials of order at most $q$ in $d$ variables, 
that is of degree strictly less than $q$,
$$
\Pi^d_q=\spann\{\x^\alpha:\alpha\in\ZZ^d_+,\;|\alpha|<q\},\quad \x^\alpha:=x_1^{\alpha_1}\cdots x_d^{\alpha_d},
\quad |\alpha|:=\alpha_1+\cdots+\alpha_d. $$
Consistent schemes satisfy \eqref{sumw}, which follows from the  exactness of \eqref{ndfeta} 
for the constant $f=\x^0\equiv 1$.
 
Stability and consistency of a positive scheme do not yet guarantee that the numerical solution  $U(t,\x)$ 
approximates the entropy solution rather than some other (unphysical) solution of \eqref{Claw} if shocks appear.
A common approach to enforce this is based on the {\it principle of vanishing viscosity}, see for example \cite{Evans10}, 
that says that the physical solution
$u$ of \eqref{Claw} is obtained as a limit for $\mu\to0^+$ of the smooth solution $u_\mu$ of the parabolic 
initial value problem 
\begin{equation} 
\label{paraeq} u_t + \F' (u)  \nabla_\x u = \mu \Delta_\x u,\qquad u(0,\x) = u_{0,\mu}(\x),\quad \mu>0,
\end{equation}
where $u_{0,\mu}$ is a smooth approximation of $u_0$ converging to $u_0$ as $\mu\to0+$ in $L_2$-norm. Therefore, we also
choose a numerical differentiation formula for the space Laplacian 
$\Delta_\x=\frac{\partial^2}{\partial x_1^2}+\cdots+\frac{\partial^2}{\partial x_d^2}$ in the form
\begin{equation} \label{ndfL}
\Delta_\x f(\x_i)\approx \sum_{\x_j\in X_i}v_{ij} f(\x_j),
\end{equation}
and require exactness for all quadratic polynomials,
\begin{equation} \label{exa2}
\Delta_\x  p(\x_i)= \sum_{\x_j\in X_i}v_{ij} p(\x_j)\qquad \forall p\in \Pi^d_3.
\end{equation}
The scheme \eqref{schm} is modified by adding an {\it artificial viscosity} term multiplied by a small positive $\mu_i$, 
which will be allowed to be chosen individually for each $\x_i$,
\begin{equation} \label{schmnu}
U(t+\Delta t,\x_i)=U(t,\x_i) - \Delta t \sum_{\x_j\in X_i}w_{ij}U(t,\x_j)+\mu_i \Delta t \sum_{\x_j\in X_i}v_{ij}U(t,\x_j),
\qquad \mu_i>0.
\end{equation}
The modified scheme is positive if and only if
\begin{equation} \label{mnu1}
\Delta t\,(w_{ii}-\mu_i v_{ii}) \le 1
\end{equation} 
and 
\begin{equation}
\label{mnu2}
w_{ij}-\mu_i v_{ij} \le 0,\quad \forall\x_j\in X_i\setminus\{\x_i\}.
\end{equation} 
Exactness of \eqref{ndfL} for quadratic polynomials implies in particular exactness for constants, that is
\begin{equation} \label{sumv}
\sum_{\x_j\in X_i}v_{ij} =0,
\end{equation}
hence the modified scheme retains the convex combination property and also satisfies the local maximum principle  as soon as
it is positive. Moreover, $\Delta p=0$ for all linear polynomials, hence
\begin{equation} \label{zlin}
\sum_{\x_j\in X_i}v_{ij} p(\x_j)=0\qquad \forall p\in \Pi^d_2,
\end{equation}
which implies that the scheme \eqref{schmnu} written as \eqref{schm} with $w_{ij}$ replaced by 
$\tilde w_{ij}:=w_{ij}-\mu_i v_{ij}$ remains a consistent discretization of \eqref{Claweq} since
\begin{equation} \label{exa3}
D_{\boeta_i}p(\x_i)= \sum_{\x_j\in X_i}\tilde w_{ij} p(\x_j)\qquad \forall p\in \Pi^d_2.
\end{equation}
Note that we require exactness of \eqref{ndfL} for quadratic polynomials, rather than only for linear polynomials,
in order to ensure that $\sum_{\x_j\in X_i}v_{ij}U(t,\x_j)$ in the artificial viscosity term provides an approximation
of the physical viscosity $\Delta_\x u$ at $\x_i$.

In what follows we obtain the weights $w_{ij}$ and $v_{ij}$ by solving for each $\x_i$ a small quadratic optimization problem
with  linear equality and inequality  constraints, see Section~\ref{sec:ndf}. 
It is clear that  constraints \eqref{m1} and \eqref{mnu1} become less restrictive if $\Delta t$ is smaller.
We refer to \eqref{m1} and \eqref{mnu1} as {\it CFL conditions} even if we do not provide any bound for $\Delta t$
that guarantees positivity and hence stability of the scheme, which is normally associated with
CFL type conditions going back to \cite{CFL28}. 
Indeed, since we do not assume any regularity of the distribution of the nodes in $X$, it would not be efficient, 
for example, to first choose the weights $w_{ij}$ and then take $\Delta t$ so small that \eqref{m1} is satisfied 
for all nodes $\x_i$.

\blue{Even if the choice of the time step $\Delta t$ and the viscosity coefficients $\mu_i$ are in general problem dependent,
and fine tuning them may significantly improve the results in the applications, it is important to take into account their
natural rescaling depending on the size of the  characteristic velocity. It is easy to see that the solution of (\ref{Claweq})
becomes $u(\gamma t,\x)$ if $\F'(u)$ is replaced by $\gamma \F'(u)$, and the solution of (\ref{paraeq})
becomes $u(\gamma t,\x)$ if $\F'(u)$ is replaced by $\gamma \F'(u)$, while $\mu$ is replaced by $\gamma \mu$,
for any $\gamma>0$. Therefore, in order to obtain the same behavior also in the discrete case, one needs to ensure that
the discrete solution $U(t,\x_i)$ becomes $U(\gamma t,\x_i)$ under the same scaling of $\F'(u)$ and $\mu_i$.
Fortunately, this is easy to achieve for the scheme \eqref{schmnu} by replacing the time step $\Delta t$ with 
$\gamma^{-1}\Delta t$, and $\mu_i$ with $\gamma\mu_i$, because the coefficients $w_{ij}$ will become $\gamma w_{ij}$, and
hence $\gamma$ cancels out in both the divergence and the viscosity terms of (\ref{paraeq}). This justifies the usual 
approach of making $\mu$ directly and $\Delta t$ inversely proportional to the size of the characteristic velocity $\F'(u)$.}

Note that in the case of constant linear transport ($\F'=$ const)  and smooth initial
condition $u_0$ our scheme is monotone since its weights $w_{ij}$ only depend on the local distribution of nodes near $\x_i$, 
see Section~\ref{sec:ndf}.
Hence,  we may only expect convergence of first
order for smooth solutions, see for example Theorem~4.12  in \cite{hestbook}, and thus the Euler method of \eqref{schm} is appropriate for time
discretization.

%% file: section03.tex
We obtain numerical differentiation formulas \eqref{ndfeta} and \eqref{ndfL} by minimizing certain weighted $\ell_2$-seminorms
of the weight vectors $[w_{ij}]_{\x_j\in X_i}$ and  $[v_{ij}]_{\x_j\in X_i}$ under equality constraints that ensure the
exactness for linear, respectively, quadratic polynomials, and inequality constraints that guarantee
the positivity of the scheme.

Minimization of $\ell_2$-seminorms is a standard approach in the meshless finite difference methods of polynomial type,
it is equivalent to obtaining numerical differentiation formulas by applying a differential operator to a weighted 
 least squares polynomial approximation of the data.  
There are many choices for the least square weights employed in the literature, see for example  \cite{JTAMB19}. 
We use the weights of the \textit{$\ell_2$-minimal formulas} suggested in \cite{davydov2018} on the basis of the error
bounds of numerical differentiation. 
In this section we first review general notions and results related to these formulas. 
Subsequently, we present our approach to enforcing required properties of the weights.

Let 
\begin{equation} \label{diffop}
 Df(\x) := \sum_{\alpha \in \ZZ_+^d,\, |\alpha| \le k} c_\alpha(\x)
 \frac{\partial^{|\alpha|} f}{\partial x_1^{\alpha_1}\cdots\partial x_d^{\alpha_d}}(\x)\,
\end{equation}
be a linear differential operator of order $k$  applied to a sufficiently regular function $f$ at a point $\x$. 
\blue{For any $\x\in\overline{\Omega}$ and a set $X_{loc}\subset X$ of its neighbors in $X$, consider the 
numerical approximation 
\begin{equation}\label{diff1}
\hat Df(\x):=\sum_{\x_j\in X_{loc}} w_j f(\x_j)\,
\end{equation}
of $Df(\x) $ defined as a linear combination of the values of $f$ on $X_{loc}$. 
In this paper we use two types of operators: directional derivatives $D=D_{\boeta_i}$ and the Laplace operator 
$D=\Delta$, and perform numerical differentiation at nodes $\x=\x_i$ on the influence sets $X_{loc}=X_i$ with 
weights denoted $w_j=w_{ij}$ for $D_{\boeta_i}$ and $w_j=v_{ij}$ for $\Delta$.}

 If the weight vector $\w=[w_j]_{\x_j\in X_{loc}}$  is chosen such that \eqref{diff1} is exact for polynomials up to a
certain order $q>k$,\footnote{The existence of the solution obtained by imposition of the exactness conditions, 
which may not be achieved for some special geometry of the set $X_{loc}$, is guaranteed for example if $X_{loc}$ is  
unisolvent for the polynomial space of interest.} 
then the following upper bound for the error holds \cite[Theorem 7]{davydov2018},
\begin{equation}  
\label{erf} 
\vert \hat Df(\x) - Df(\x) \vert \le \sigma(\x, X_{loc},\w)\, h_{\x,X_{loc}}^{r+\gamma-k}  \vert f 
\vert_{r,\gamma,\Omega_{loc}}\,, 
\quad f \in C^{r,\gamma}(\Omega_{loc}), \quad k\le r \le q-1, 
\end{equation} 
where  $\Omega_{loc}$ denotes a domain containing the set
$\displaystyle{\bigcup_{\x_j\in X_{loc}} [\x\,,\,\x_j] }\,,$
\[
\sigma(\x, X_{loc}, \w):= h_{\x,X_{loc}}^{k-r-\gamma} \sum_{\x_j\in X_{loc}}\vert w_j \vert \|\x_j-\x\|_2^{r+\gamma}\,,
\qquad
h_{\x,X_{loc}} := \max_{\x_j\in X_{loc}} \|\x_j - \x\|_2 \,,
\]
$$ \vert f \vert_{r,\gamma,\Omega_{loc}} \,:=\,  \frac{1}{(\gamma+1)\cdots (\gamma+r)}\, 
\Big( \sum_{\vert \alpha \vert = r}  {r\choose{\alpha}} 
\Big\vert \frac{\partial^{|\alpha|} f}{\partial x_1^{\alpha_1}\cdots\partial x_d^{\alpha_d}}
\Big\vert^2_{0,\gamma,\Omega_{loc}} \Big)^{1/2},
\quad |f|_{0,\gamma,\Omega_{loc}}:=\sup_{\x,\y\in \Omega_{loc}\atop \x\ne\y}\frac{|f(\x)-f(\y)|}{\|\x-\y\|_2^\gamma}, $$
and $C^{r,\gamma}(\Omega_{loc})$, $r\ge 0$, $\gamma\in(0,1]$, denotes the space consisting of all $r$ times continuously
differentiable functions $f$ on $\Omega_{loc}$ such that $\vert f \vert_{r,\gamma,\Omega_{loc}}<\infty$.
Note that the value of $\sigma(\x, X_{loc},\w)$ does not depend on the size $h_{\x,X_{loc}}$ of the neighborhood of $\x$,
but only on the cardinality and the shape of $X_{loc}$, as well as on the choice of the weights. Therefore the convergence
order of the approximation as $h_{\x,X_{loc}}\to0$ is determined by the power of $h_{\x,X_{loc}}$ in \eqref{erf} 
unless the sets $X_{loc}$ are in specific bad positions  (for example, contained in zero sets of some polynomials of order
$q$), that normally do not realize on irregular nodes.

As long as the exactness for polynomials of order up to $q>k$ does not determine the weights $w_j$ uniquely, in particular
when the number $n_{loc}$ of nodes in $X_{loc}$ is greater than$\dim \Pi^d_q={d+q-1\choose d}$,  different strategies for a suitable
selection of the weights can be considered.
The obvious choice of minimizing $\sigma(\x, X_{loc},\w)$ leads to the {\it $\ell_1$-minimal formulas} that can be computed 
by linear programming. A computationally more attractive alternative is to minimize the weighted $\ell_2$-seminorm
\begin{equation} 
\label{ell2w}  
 \vert {\bf w} \vert_{2,r+\gamma} :=\Big(\sum_{\x_j\in X_{loc}}  w_j^2 \Vert { \x_j}-{\x}\Vert_2^{2(r+\gamma)}\Big)^{1/2},
\end{equation}
which leads to the {\it $\ell_2$-minimal formula} \eqref{diff1} satisfying
\begin{equation} 
\label{erf2} 
\vert \hat Df(\x) - Df(\x) \vert \le \sqrt{n_{loc}}\,\rho(\x, X_{loc})\, h_{\x,X_{loc}}^{r+\gamma-k}  
\vert f \vert_{r,\gamma,\Omega_{loc}}\,, 
\quad f \in C^{r,\gamma}(\Omega_{loc}), \quad k\le r \le q-1, 
\end{equation} 
where
$$
\rho(\x, X_{loc}):=\inf\Big\{\sigma(\x, X_{loc},\w):\,\w\in \RR^{n_{loc}},\; Dp(\x)=\sum_{\x_j\in X_{loc}} w_j p(\x_j)\;
\text{ for all }p\in \Pi^d_q\Big\}.$$
Despite the additional factor $\sqrt{n_{loc}}$ in \eqref{erf2} in
comparison to the  $\ell_1$-minimal formulas, both types of formulas  behave
similarly in numerical experiments, see further details in \cite{davydov2018}.

\subsection*{Inequality constrained minimal formulas}

In the case of Laplace operator $\Delta f=\frac{\partial^2}{\partial x_1^2}+\cdots+\frac{\partial^2}{\partial x_d^2}$
it is often attractive to enforce the {\it positivity} of the weight vector $\w$ in \eqref{diff1}, 
by requiring that $\x\in X_{loc}$ and $w_j\ge 0$ for all $\x_j\in X_{loc}\setminus\{\x\}$. 
(Then $w_j<0$ for $\x_j=\x$ because of the exactness for constant.) 
This for example helps to ensure that system matrices of the resulting meshless finite difference method 
for the Poisson equation are M-matrices, with various benefits for the numerical solution, see  
\cite{Seibold08} and references therein. In \cite[Section 4.3]{davydov2018} an error bound is obtained for positive formulas
for elliptic differential operators of second order.

We generate inequality constrained $\ell_2$-minimal formulas in the form \eqref{ndfeta} 
and \eqref{ndfL} for the operators $D_{\boeta_i}(\x_i)$ and $\Delta_\x(\x_i)$ 
\blue{that minimize \eqref{ell2w} with parameters $r,\gamma$ chosen to take into account the
error bound \eqref{erf}} as follows. The algorithm depends 
on the prescribed bounds $n_{\min},n_{\max}$ for the size of the set of influence and the artificial viscosity parameter 
$\mu_i\ge0$.
We start by choosing an initial set of influence $X_i=X_i^{init}$ consisting of $|X_i^{init}|=n_{\min}$ nearest neighbors of
$\x_i$ in $X$, including $\x_i$ itself. (We use $n_{\min}=10$ in the numerical experiments in the bivariate case $d=2$.) 

If $\mu_i>0$, then we first compute  the weight vector $\v_i=[v_{ij}]_{\x_j\in X_i}$, $X_i=X_i^{visc}$, as the solution of the
quadratic minimization problem
\begin{equation}
\min_{\v}|\v|_{2,3}^2\quad\text{subject to \eqref{exa2} and the inequalities }v_{ij}\ge 0\;\,\forall \x_j\in X_i\setminus\{\x_i\},
\label{weiv}
\end{equation}
where 
$$
|\v|_{2,3}=\Big(\sum_{\x_j\in X_i}  v_{ij}^2 \| { \x_j}-{\x_i}\|_2^{6}\Big)^{1/2},$$ 
which corresponds to the case $q=3,r=2,\gamma=1$ of the $\ell_2$-minimal formulas, and $X_i=X_i^{visc}$ is either
$X_i^{init}$, or a larger set obtained by repeatedly expanding $X_i$ to the set of $\lceil 1.2\, |X_i| \rceil$  nearest
neighbors of $\x_i$ should \eqref{weiv} be infeasible. In order to avoid an infinite loop we terminate the process if the new
size is greater than $n_{\max}$, in which case we refrain from using artificial viscosity for the node $\x_i$ 
and set $\mu_i=0$. In our bivariate experiments with  $n_{\max}=100$ this never happened for interior nodes $\x_i$. Since \eqref{weiv}
is often infeasible for boundary nodes even for a large $X_i$, we set $\mu_i=0$ whenever $\x_i\in\partial\Omega$.

After $\v_i$ is chosen as explained above (or set formally to zero if $\mu_i=0$), we compute  the weight vector $\w_i=[w_{ij}]_{\x_j\in X_i}$ as the solution of the
quadratic minimization problem
\begin{equation}
\min_{\w}|\w|_{2,2}^2\quad\text{subject to \eqref{exa1} and the inequalities $w_{ii}\le B$, }w_{ij}\le 0\;\,\forall \x_j\in X_i\setminus\{\x_i\},
\label{weiw}
\end{equation}
where 
$$
|\w|_{2,2}=\Big(\sum_{\x_j\in X_i}  w_{ij}^2 \| { \x_j}-{\x_i}\|_2^{4}\Big)^{1/2},$$ 
which corresponds to the case $q=2,r=1,\gamma=1$ of the $\ell_2$-minimal formulas, the constant $B$ is given by
\begin{equation}
B=\max\big\{\tfrac{1}{2\Delta t},\tfrac{1}{\Delta t}-\mu_i |v_{ii}|\big\}
\label{w1bound}
\end{equation}
and $X_i$  is initialized as $X_i^{init}$ if $\mu_i=0$ or $X_i^{visc}$ otherwise, and repeatedly expanded  to the set of
$\lceil 1.2\, |X_i| \rceil$  nearest neighbors of $\x_i$ as long as \eqref{weiw} is infeasible.  Note that $v_{ii}$ is
negative as soon as $\mu_i>0$ because of \eqref{sumv} and since $v_{ij}$, $\x_j\in X_i$, cannot all be zero thanks to
\eqref{exa2}. Again, we terminate the process if the new size is greater than $n_{\max}$, meaning the algorithm exits with
failure. In this case we repeat the same procedure  as above after removing the inequality constraints in \eqref{weiw} and
resetting $X_i$ to $X_i^{init}$ if $\mu_i=0$ or $X_i^{visc}$ otherwise. With $n_{\max}=100$ in the bivariate case  this
happened very rarely and only for certain interior nodes situated very close to the inflow boundary nodes in Example 1 below,
see   Section~\ref{ex1}.  
\blue{We set $v_{ij}=0$ for all $j$ such that $\x_j\in X_i\setminus X_i^{visc}$.}
Note that the largest size of $X_i$ observed in the numerical experiments of Section~\ref{sec:ne}
was 27 for both the viscosity term  ($X_i^{visc}$) and  the divergence term (final $X_i$) \blue{when using the nodes of the 
Cartesian grid or Halton points, and reaching at most 40 for the viscosity term and 48 for the divergence term in the case of
random nodes of Example~2 in Section~\ref{ex2}}.

Finally, we correct $\mu_i$ if needed, by setting
\begin{equation}
\mu_i=\min\big\{\mu_i,\tfrac{1}{2\Delta t|v_{ii}|}\big\}.
\label{nucor}
\end{equation}
It is easy to check that the resulting weights satisfy the positivity conditions \eqref{mnu1} and \eqref{mnu2}
with updated $\mu_i$.

Clearly, the algorithm described here is not the only possible scenario to ensure the positivity of the scheme. Minimization
of $|\v|_{2,3}$ and $|\w|_{2,2}$ is motivated by the error bound \eqref{erf2}. The weights $v_{ij}$ are computed before 
$w_{ij}$ because after knowing $v_{ii}$ we often obtain via \eqref{w1bound} a weaker constraint 
$w_{ii}\le B$ in \eqref{weiw}, which may help to somewhat reduce the size of $X_i$ and hence reduce the factor 
$\sigma(\z, X_{\z},\w)$  in \eqref{erf}, improving the accuracy of the numerical differentiation.
Different approaches to computing $v_{ij}$, $w_{ij}$ and organizing their interaction and adjustment of $\mu_i$ 
are possible and should be explored in the future. The version presented here works well in the numerical experiments
presented in Section~\ref{sec:ne}.

\begin{remark}\label{posvisc}
Note that the positivity of the scheme may be enforced via the choice of the positive artificial viscosity weight vector $\v$ 
and the coefficient $\mu_i>0$, which would allow the weight vectors  $\w$ for the directional derivatives to
fail the condition $w_{ij}\le 0$, $j\ne i$. However, this  would mean that we have to apply
stronger inequality constraints on $\v$ in the form $\mu_iv_{ij}\ge w_{ij}$ whenever $w_{ij}>0$ (assuming $\w$ is computed
before $\v$). In particular, we would not be able to choose $\mu_i=0$, which would rule out adaptive viscosity as in
Algorithm~\ref{adaptnu} below, where $\mu_i$ vanishes at some distance from the faults in the solution. Moreover, a sufficiently large $\mu_i>0$ may be needed in this case, increasing the numerical
dissipation unnecessarily. In contrast to this, our approach allows to add just enough artificial viscosity to ensure 
shock capturing, without burdening it by a second
task of stabilizing the scheme.

\end{remark}

%% file: section04.tex
Since the stability of our scheme is guaranteed by its positivity, the artificial viscosity term is only responsible  for the 
correct resolution of shocks. Therefore it is desirable to choose the viscosity parameter $\mu_i$ in \eqref{schmnu} 
adaptively, in particular setting it to zero for nodes $\x_i$ that are far from the shock front, which improves the local 
accuracy of the overall scheme because local numerical differentiation weights $\tilde w_{ij}$ in this case coincide 
with the weights $w_{ij}$ of the unconstrained $\ell_2$-minimal formula, in contrast to the combined weights 
$\tilde w_{ij}=w_{ij}-\mu_iv_{ij}$ when $\mu_i>0$, that lose the $\ell_2$-minimality. 
Since the solution $u$ exhibits  jump discontinuities at the shock fronts, we choose $\mu_i>0$  only for nodes that are close to the discontinuities of the discrete numerical
solution $U(t,\x_j)$, $j=1,\ldots,N$, of the previous time step, identified by a fault detection method. We refer to
\cite{guermond17}  and references therein for algorithms relying on adaptive viscosity in the context of the finite element
and finite volume  methods for conservation laws.

We make use of the fault indicator introduced in \cite{bracco2019} that also relies on the minimal differentiation formulas 
described in Section~\ref{sec:ndf}. Given a set of scattered nodes $X\subset\RR^d$, with associated function values
$f(\x)$, $\x \in X$, it allows us to determine which nodes are close to a discontinuity of $f$. In many applications
discontinuities form curves, which are referred to as {\it faults} in the literature.

\blue{For each $\x_i\in X$, let $X_{i,\cal F}\subseteq X$ 
be the set of $n_{\cal F}$ points $\x_j\in X$ closest to $\x_i$.} We define the
{\it fault indicator}
\begin{equation} \label{fault_ind1}
I_{\Delta}(\x_i,n_{\cal F}):=
\frac{| \hat\Delta f (\x_i)|}{\displaystyle \sum_{\x_j\in X_{i,\cal F}}  | w_j |\, \| \x_j - \x_i \|^2_2 }\,, 
\end{equation}
where 
\begin{equation} \label{fault_ind1b} 
\hat\Delta f (\z) := \sum_{\x_j\in X_{i,\cal F} }  w_j f(\x_j)
\end{equation}
is the  $\ell_2$-minimal differentiation formula approximating the Laplacian $\Delta f(\x_i)$ with polynomial exactness order
$q=3$ and minimizing the seminorm $|\w|_{2,3}$. \blue{Note that the second fault indicator introduced
in \cite{bracco2019}, employing the gradient instead of the Laplace operator, can also be used and gives comparable results in
numerical experiments. We opted for the Laplacian-based one because it is cheaper. Indeed, 
the approximation of the scalar-valued Laplace operator  by a numerical differentiation formula requires significantly
less computational effort than approximating the vector-valued gradient operator. }

In \cite{bracco2019} we have shown (see Theorem 2 there) the following property of the indicator \eqref{fault_ind1}, 
\begin{equation}\label{fiineq}
I_{\Delta}(\x_i,n_{\cal F})\le |f|_{1,1,\Omega_{loc}},
\end{equation}
with $\Omega_{loc}$ as in \eqref{erf}. Although \eqref{fiineq} is stated in  \cite{bracco2019} only for the bivariate case, it is easy to
extend it to any number of variables $d$. Indeed, this follows from the fact that Theorem 3 in \cite{bracco2019} holds for any
$d$ with the same proof.
It is remarkable that the estimate \eqref{fiineq} does not contain any
unknown constants and therefore allows for comparison between different $\x_i$ independently of how irregular each local  node
set $X_{i,\cal F}$ is.  Therefore, a large value of $I_{\Delta}(\x_i,n_{\cal F})$ indicates that $f$  is likely to  be
discontinuous in $\Omega_{loc}$, or at least its gradient tends to change abruptly there, that is the behavior expected from 
the solution of a conservation law in the vicinity of a  shock front. Moreover, it is shown in the bivariate case in Theorem~6
of \cite{bracco2019}  that $I_{\Delta}(\x_i,n_{\cal F})$ behaves like $O(h^{-2})$ as $\x_i$ approaches a fault and
$\hbox{diam}(X_{i,\cal F})\rightarrow 0$, while it is bounded if $\x_i$ stays sufficiently far from all faults.
In other words, the indicator tends to take large values at points close to faults and bounded values
at points far from them, which justifies classifying as close to a fault the points of the set
\begin{equation*}
{\cal F}(\alpha,X):=\{\x_i\in X:\, I_{\Delta}(\x_i,n_{\cal F}) > \alpha\},
\end{equation*}
with $\alpha>0$ being a suitable threshold. 

Thanks to these properties of the indicator, a reasonable choice of
$\alpha$ is the median value  of $I_{\Delta}$ computed on the nodes in $X$. This works well, except when $f$ is constant or
linear in large subareas of the domain, since  $\Delta f(\z)=0$ and so $I_{\Delta}(\x_i,n_{\cal F})\approx 0$ for any 
$\x_i$ belonging to these areas, and $\alpha$, chosen as the median value, would be very close to $0$. As a result, all the
points where the value of the indicator is slightly above 0 would be marked as close to a fault even if they are not. For this
reason, we use the following two-step fault detection algorithm, that has shown excellent performance in extensive numerical
tests  on scattered data presented in \cite{bracco2019}.

\medskip

Computation of the set of fault nodes ${\cal F}(X)$.
 \begin{itemize}
\item Set $\alpha_1= C_1\median\{I_{\Delta}(\x_i,n_{\cal F}):\, \x_i\in X\}$,  and  compute ${\cal F}(\alpha_1,X)$.
\item Set $\alpha_2=C_2 \median\{I_{\Delta}(\x_i,n_{\cal F}):\, \x_i\in {\cal F}(\alpha_1,X)\}$, and 
 compute the final set of fault nodes as ${\cal F}(X):= {\cal F}(\alpha_2,{\cal F}(\alpha_1,X))$.
\end{itemize} 
Note that  in all numerical experiments of Section~\ref{sec:ne} we use $C_1=1$ and $C_2=2$, 
\blue{the same values that proved to be effective in the experiments 
of \cite{bracco2019}, as our goal is to provide a robust choice, avoiding a fine-tuning of the parameters example by example.}

\subsection*{Adaptive artificial viscosity}
We now describe our approach to choosing the factor $\mu_i$ in \eqref{schmnu} when adaptive fault-based viscosity is used.
As input we  need the set of nodes $X$, the values $U(t,\x_i)$ for all $\x_i\in X$,
and the {\it spacing} $h$ for the node set $X$. Spacing $h$ is a parameter that we use in the node generation as
explained in  Section~\ref{sec:ne}, with $h\approx(|\Omega|/N)^{1/d}$. In addition we need three user specified 
parameters: $\mu$, the maximum factor for the artificial viscosity terms;  the number $n_{\cal F}$ of nearest nodes 
for the computation of  the fault indicator, which we choose to be the same number $n_{\cal F}=n_{\min}$ as in 
Section~\ref{sec:ndf}; and the constant $C_3$ that governs the transition zone between $\mu_i>0$ and $\mu_i=0$ in 
order to avoid abrupt changes in $\mu_i$, \blue{that otherwise have an adverse effect on the performance
of the method. Note that the advantages of  a continuous transition of the amount of the artificial viscosity between zero in the
regions where the solution is smooth and its maximum value where the solution is discontinuous 
have been observed for other numerical methods, see for example Section 12.2.3 
 in \cite{hestbook}.} We use $C_3=5$ in the experiments of Section~\ref{sec:ne}.

\medskip

Local viscosity parameter $\mu_i$.
\begin{itemize}
\item 
For each $\x_i\in X$, choose the set $X_{i,\cal F}=X_i^{init}$ consisting of $n_{\cal F}$ nearest neighbors of $\x_i$ in $X$,
and compute ${\cal F}(X)\subset X$ according to the above fault detection algorithm, with $f$ on $X$ given by
$f(\x_j)=U(t,\x_j)$, $j=1,\ldots,N$. 
\item 
For each $\x_i\in X$, find the distance $\rho_i$ to the set ${\cal F}(X)$, 
$\rho_i:=\min\big\{\|\x_i-\x_j\|_2:\,\x_j\in {\cal F}(X)\big\}$, and set 
\begin{equation}\label{nui}
\mu_i=\max\Big\{0,1-\tfrac{\rho_i}{C_3h}\Big\}\,\mu.
\end{equation}
In particular, $\mu_i=0$ as long as $\rho_i\ge C_3h$.
\end{itemize}
\blue{Fault points and the values of $\mu_i$ are illustrated in Figure~\ref{fig:RW_fault} in the numerical
examples section of the paper.}

%% file: sectionALG.tex
We now describe three versions of  the positive schemes for \eqref{Claweq}, differing in their approaches to the artificial
viscosity. All three algorithms rely on the following input data, while artificial viscosity requires further
setup parameters detailed for Algorithms~\ref{constnu} and \ref{adaptnu} that make use of it.
\medskip

{\bf Input:}
\begin{itemize}
\item Data of the problem: derivative of the flux $\F'$, domain $\Omega\subseteq\RR^d$, the function $u_0$ that
defines the initial condition \eqref{initcond}, final time $T$, data for the boundary conditions whenever appropriate.
\item Parameters of the numerical scheme: discretization nodes $X$, step-size in time $\Delta t>0$, 
minimum and maximum number of closest neighbors  $n_{\min}$ and $n_{\max}$ for the sets of influence.
\end{itemize} 

The {\bf output} of each scheme is the set of values
$$
U(k\Delta t,\x_i),\quad i=1,\ldots,N, \;\text{ and }\;k\in \NN\;\text{ with }\;k\Delta t\le T,$$
as approximations of $u(k\Delta t,\x_i)$, for the  solution $u$ of the conservation law \eqref{Claw}.
\medskip 

Our first algorithm does not add any artificial viscosity, which is appropriate when the  shocks are absent.

\begin{algorithm}\label{nonu}
{\bf Positive scheme without artificial viscosity.}
\begin{itemize}
\item Initialization: Set $t=0$ and $U(0,\x_i)=u(0,\x_i)$, $i=1,\ldots,N$. 
\item While $t+\Delta t\le T$: 
\begin{itemize}
\item	For each $\x_i$, $i=1,\ldots,N$:
		\begin{enumerate}
		\item Compute the set of influence $X_i$ and the weights $w_{ij}$, $j\in X_i$, of an
		inequality constrained minimal formula, as described in Section~\ref{sec:ndf}, with $\mu_i=0$.
		\item Compute $U(t+\Delta t,\x_i)$ by \eqref{schm}.
		\end{enumerate}
\item Set $t:=t+\Delta t$.

\end{itemize}
\end{itemize}
\end{algorithm}

As mentioned above, adding small artificial viscosity is motivated by the principle of vanishing viscosity, and its goal is
to allow accurate computation of the solution in the presence of shocks. As a side effect, solutions computed with 
artificial viscosity are smoother, which is normally not desired, but sometimes also produces better shapes by reducing
non-physical behavior of the numerical solution near discontinuities. This latter positive effect is not significant for our
method because consistent positive schemes satisfy the local maximum principle \eqref{locmax} and therefore do not suffer from
artificial oscillations. 

\begin{algorithm}\label{constnu}
{\bf Positive scheme with constant artificial viscosity.}
\begin{itemize}
\item Additional input parameter: viscosity factor $\mu>0$.
\item Initialization: Set $t=0$ and $U(0,\x_i)=u(0,\x_i)$, $i=1,\ldots,N$. 
\item While $t+\Delta t\le T$: 
\begin{itemize}
\item
	For each $\x_i$, $i=1,\ldots,N$:
		\begin{enumerate}
		\item Compute the set of influence $X_i$ and the weights $v_{ij}$ and $w_{ij}$, $j\in X_i$, of the
		inequality constrained minimal formulas, as described in Section~\ref{sec:ndf}, with $\mu_i=\mu$.
		\item Compute $U(t+\Delta t,\x_i)$ by \eqref{schmnu} with $\mu_i=\mu$.
		\end{enumerate}
\item Set $t:=t+\Delta t$.
\end{itemize}
\end{itemize}
\end{algorithm}

Finally, we consider a scheme where viscosity is only added for points that are close to discontinuities (faults) detected
by the method described in Section~\ref{sec:df}. The goal is to improve approximation quality in the smooth regions because
the resulting discretization of $D_{\boeta_i}$ in \eqref{ndfeta} will allow a smaller set of influence $X_i$ and smaller
seminorm $|\w_i|_{2,2}$ of the weight vector $\w_i=[w_{ij}]_{\x_j\in X_i}$ minimized in  \eqref{weiw}, in comparison to the 
weights $\tilde w_{ij}=w_{ij}-\mu_i v_{ij}$ effectively used in \eqref{schmnu} when $\mu_i>0$. In addition, this scheme avoids
the computation of the positive weights for the approximation of the Laplacian in most of the domain, leading to a significant
computational saving. In fact, the computation of the approximation of the Laplacian involved in the expression of the fault
indicator \eqref{fault_ind1} is less expensive since a minimal differentiation formula without sign constraints on the weights
is employed, without an adaptive expansion of the influence set needed for the
inequality constrained weights.

\begin{algorithm}\label{adaptnu}
{\bf Positive scheme with adaptive fault-based artificial viscosity.}
\begin{itemize}
\item Additional input parameters: maximum viscosity factor $\mu>0$, spacing $h$ of the node set $X$, fault detection
parameters $n_{\cal F},C_1,C_2$ and transition  zone parameter $C_3$.
\item Initialization: Set $t=0$ and $U(0,\x_i)=u(0,\x_i)$, $i=1,\ldots,N$. 
\item While $t+\Delta t\le T$:
\begin{itemize}
\item
Apply the fault detection method presented in Section~\ref{sec:df} to the point set $X$ and the corresponding solution
values $\{U(t,\x),\, \x\in X\}$ with parameters $n_{\cal F},C_1,C_2$, obtaining the set of fault nodes
${\cal F}(X)\subset X$.
\item
For each $\x_i$, $i=1,\ldots,N$:
		\begin{enumerate}
		\item Compute the individual viscosity factor $\mu_i$ by \eqref{nui}, using ${\cal F}(X)$, $\mu$, $h$ and $C_3$
		as described in Section~\ref{sec:df}. 
		\item Compute the set of influence $X_i$ and the weights $v_{ij}$ and $w_{ij}$, $j\in X_i$, of the
		inequality constrained minimal formulas, as described in Section~\ref{sec:ndf}.
		\item Compute $U(t+\Delta t,\x_i)$ by \eqref{schmnu}.
		\end{enumerate}
\item Set $t:=t+\Delta t$.
\end{itemize}
\end{itemize}
\end{algorithm}

%% file: section05.tex
\renewcommand\thefigure{\arabic{section}.\arabic{figure}}  
\renewcommand\thetable{\arabic{section}.\arabic{table}}  

In this section we report numerical results obtained by applying our positive schemes to three bivariate benchmark problems selected to
highlight  different aspects of the method.

 In order to test our method on irregular nodes we initialize $X$ as a set of $M$ quasi-random Halton points in a square
$\Omega\subset\RR^2$, where $M$ is determined by the {\it spacing} parameter $h>0$ and side length  $\ell$ of the square by rounding
$(\ell/h)^2$ to the nearest integer. Then all nodes with distance to the boundary of $\Omega$ at most $0.25 h$ are removed
from $X$. 
\blue{Although Halton points already provide a good representative of irregularities expected in meshless node generation, 
see \cite{SJD22}, in Example 2 we also test the performance of our method on the more severe random points obtained in a similar manner 
by MATLAB's pseudo-random number generator.}
In the case of periodic boundary conditions (Examples 2 and 3)  we extend $X$ periodically to a larger domain, which
simplifies the computation of the sets of influence.
Example 1 requires additional nodes on the boundary, as explained in Section~\ref{ex1}. Whenever we use a Cartesian grid, $h$
is the usual step size. Therefore the same spacing $h$ corresponds to comparable numbers $N$ of nodes for the grids and Halton
points.

The parameters of our schemes are chosen in a unified way in all examples. 
\blue{In particular, we make the  viscosity  factor $\mu$ directly and the time step $\Delta t$ inversely proportional to the
maximum size of the characteristic velocity $\F'(u_0)$ by setting $\mu = 0.5 h \, v_0$ and 
$\Delta t = 0.2h\, / v_0$,  where $v_0 := \max_{\x \in \bar \Omega} \Vert \F'(u_0) \Vert_\infty$.}
 We always take $n_{\min}=10$, $n_{\max}=100$, $\Delta t=0.2 h$ and choose the final time $T$ such that $K=T/\Delta t$ is an
integer.
 The viscosity or maximum viscosity factor is always $\mu=0.5h$, and additional parameters of the fault-based 
viscosity are chosen as $n_{\cal F}=10,C_1=1,C_2=2,C_3=5$. 

In addition to the qualitative assessment of the results by comparing 3D plots of the solutions, we also compute 
the error measures
$$
E_1=\tfrac{1}{N}\sum_{i=1}^N|U(T,\x_i)-U_{\rm ref}(T,\x_i)|,
\quad\text{and}\quad E_2=\Big(\tfrac{1}{N}\sum_{i=1}^N(U(T,\x_i)-U_{\rm ref}(T,\x_i))^2\Big)^{1/2}$$
with respect to a reference solution $U_{\rm ref}$ given either by an exact solution (available in Example~1), 
or by a high resolution solution, obtained on a dense Cartesian grid with an appropriate higher order method 
provided in the code supplement of the book \cite{hestbook}. In the case of an irregular set $X$, the values $U_{\rm ref}(T,\x_i)$ are obtained by piecewise linear
interpolation of  $U_{\rm ref}(T,\cdot)$ on the Delaunay triangulation of the Cartesian grid. Note that $E_1$ and $E_2$ may be considered as discrete 
replacements for the $L^1$ and $L^2$ errors on $\Omega$, respectively.

\subsection{Example~1: Inviscid Burgers' equation with discontinuous initial condition}\label{ex1}
\begin{figure}[!t]
\subfigure[Exact solution.]{\includegraphics[scale=.37]{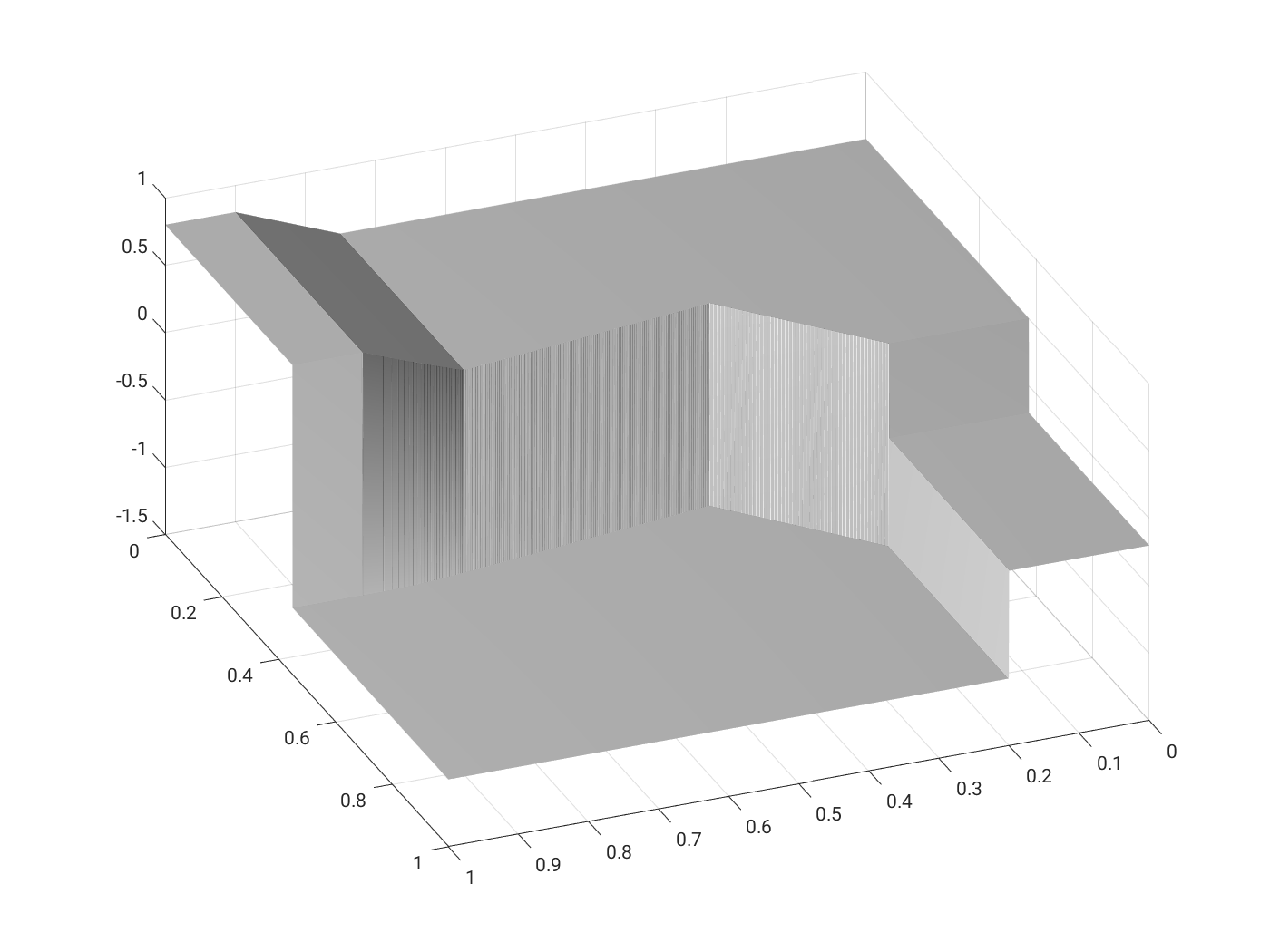}}
\subfigure[Halton nodes and influence sets.]{\includegraphics[scale=.57]{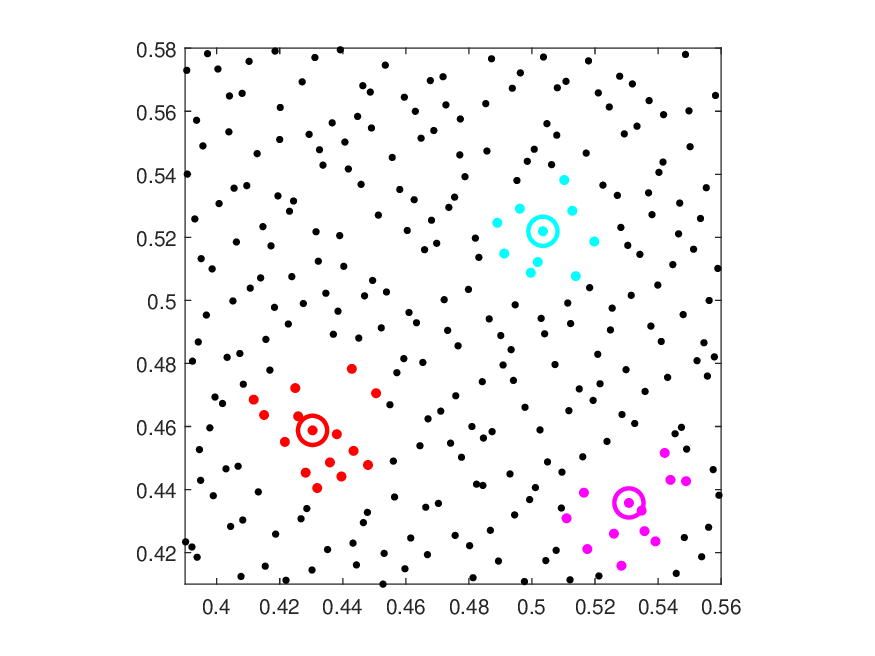}} 
\caption{\blue{Example 1: Exact solution and a zoom of the node set, with three highlighted influence sets $X_i$ 
(the encircled point is $\x_i$) of size $10$ (in cyan), $12$ (in magenta), and $15$ (in red).}}
\label{fig:ex2nodes}     
\end{figure}

In the first example we consider the inviscid Burgers' equation with  flux $\F(u) = \frac{1}{2}u^2\bfv$, $\bfv = (1\,,\,1)$ 
on the spatial domain $\Omega = [0,\,1]^2$, with discontinuous initial condition
\begin{equation*}
u_0(\x) = u_0(x,y) =  \begin{cases}{} 
-0.2 & \hbox{if}\, x<0.5\, \hbox{and}\, y>0.5,\\
-1 & \hbox{if}\, x>0.5\, \hbox{and}\, y>0.5,\\
0.5 & \hbox{if}\, x<0.5\, \hbox{and}\, y<0.5,\\
0.8 & \hbox{if}\, x>0.5\, \hbox{and}\, y<0.5,
\end{cases} 
\end{equation*}
and exact solution given for $t>0$ by
\begin{equation*}
u(\x) = u(x,y) = \begin{cases}{} 
-0.2 & \hbox{if}\, x<\frac{1}{2}-\frac{3t}{5}\, \hbox{and}\, y>\frac{1}{2}+\frac{3t}{20},\\
0.5 & \hbox{if}\, x<\frac{1}{2}-\frac{3t}{5}\, \hbox{and}\, y<\frac{1}{2}+\frac{3t}{20},\\
-1 & \hbox{if}\, \frac{1}{2}-\frac{3t}{5}<x<\frac{1}{2}-\frac{t}{4}\, \hbox{and}\, y>\frac{-8x}{7}+\frac{15}{4}-\frac{15t}{28},\\
0.5 & \hbox{if}\, \frac{1}{2}-\frac{3t}{5}<x<\frac{1}{2}-\frac{t}{4}\, \hbox{and}\, y<\frac{-8x}{7}+\frac{15}{4}-\frac{15t}{28},\\
-1 & \hbox{if}\, \frac{1}{2}-\frac{t}{4}<x<\frac{1}{2}+\frac{t}{2}\, \hbox{and}\, y>\frac{x}{6}+\frac{5}{12}-\frac{5t}{24},\\
0.5 & \hbox{if}\, \frac{1}{2}-\frac{t}{4}<x<\frac{1}{2}+\frac{t}{2}\, \hbox{and}\, y<\frac{x}{6}+\frac{5}{12}-\frac{5t}{24},\\
-1 & \hbox{if}\, \frac{1}{2}+\frac{t}{2}<x<\frac{1}{2}+\frac{4t}{5}\, \hbox{and}\, y>x-\frac{5}{18t}(x+t-\frac{1}{2})^2,\\
\frac{2x-1}{2t} & \hbox{if}\, \frac{1}{2}+\frac{t}{2}<x<\frac{1}{2}+\frac{4t}{5}\, \hbox{and}\, y<x-\frac{5}{18t}(x+t-\frac{1}{2})^2,\\
-1 & \hbox{if}\, x>\frac{1}{2}+\frac{4t}{5}\, \hbox{and}\, y>\frac{1}{2}-\frac{t}{10},\\
0.8 & \hbox{if}\, x>\frac{1}{2}+\frac{4t}{5}\, \hbox{and}\, y<\frac{1}{2}-\frac{t}{10},
\end{cases}
\end{equation*}
see \cite{guermond11}. Inflow boundary conditions are derived from the above formula for the exact solution. \blue{The solution at
time $T=0.5$ is presented in Figure~\ref{fig:ex2nodes}(a).}

\begin{figure}[!p]
\subfigure[No viscosity (Algorithm~\ref{nonu}).]%
{\includegraphics[scale=.339]{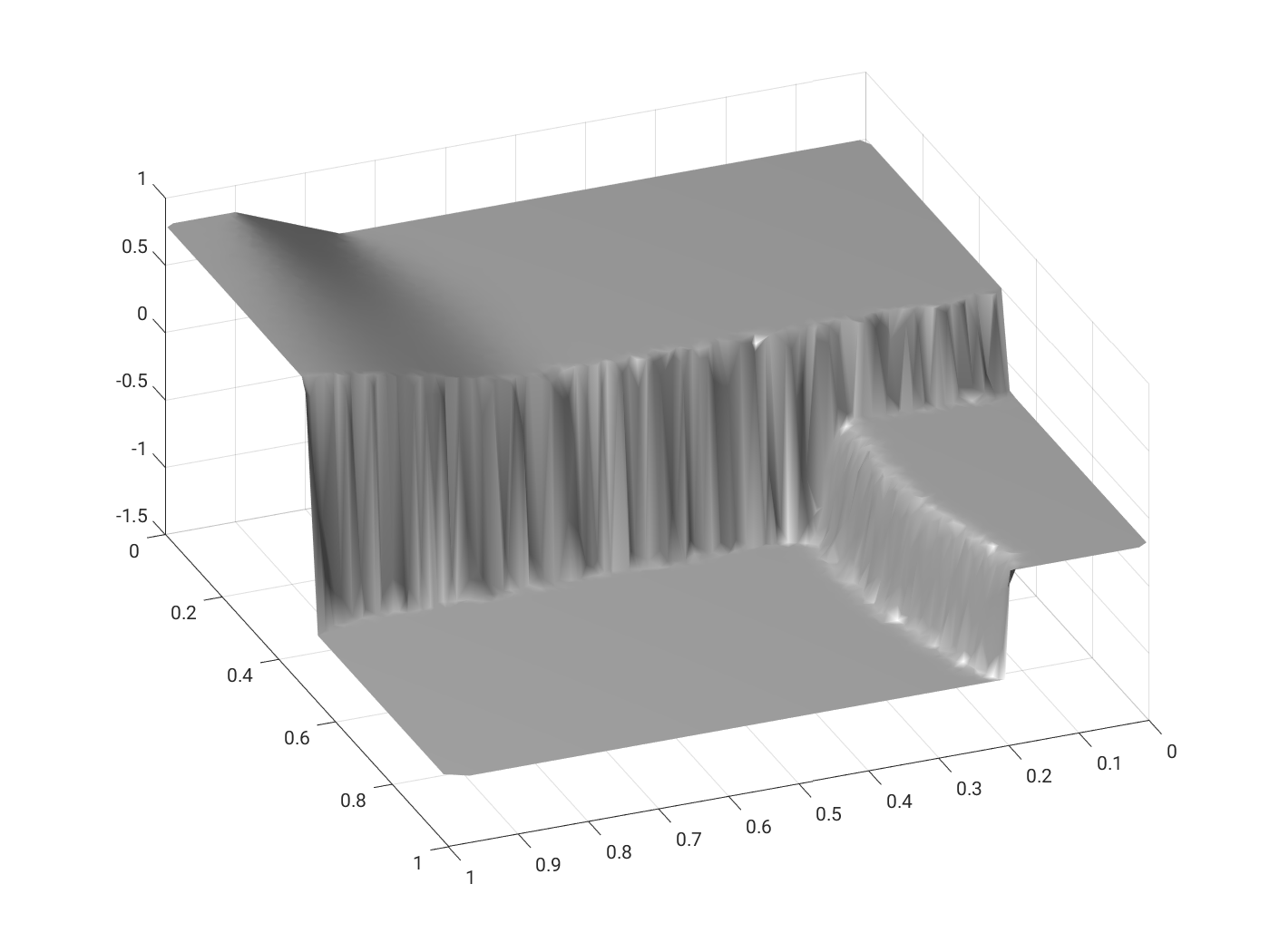}
\includegraphics[scale=.339]{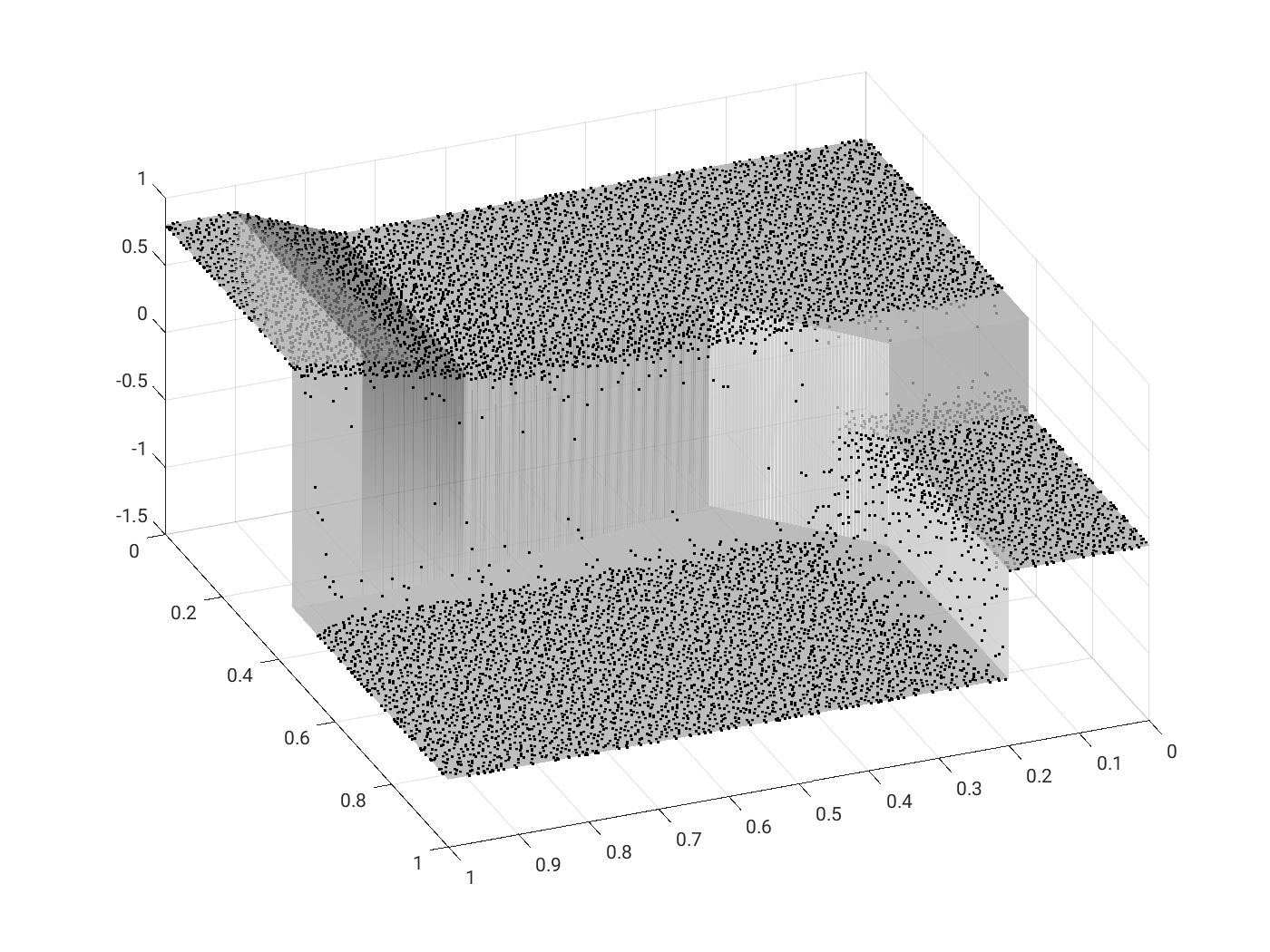}
}
\subfigure[Constant viscosity (Algorithm~\ref{constnu}).]%
{\includegraphics[scale=.339]{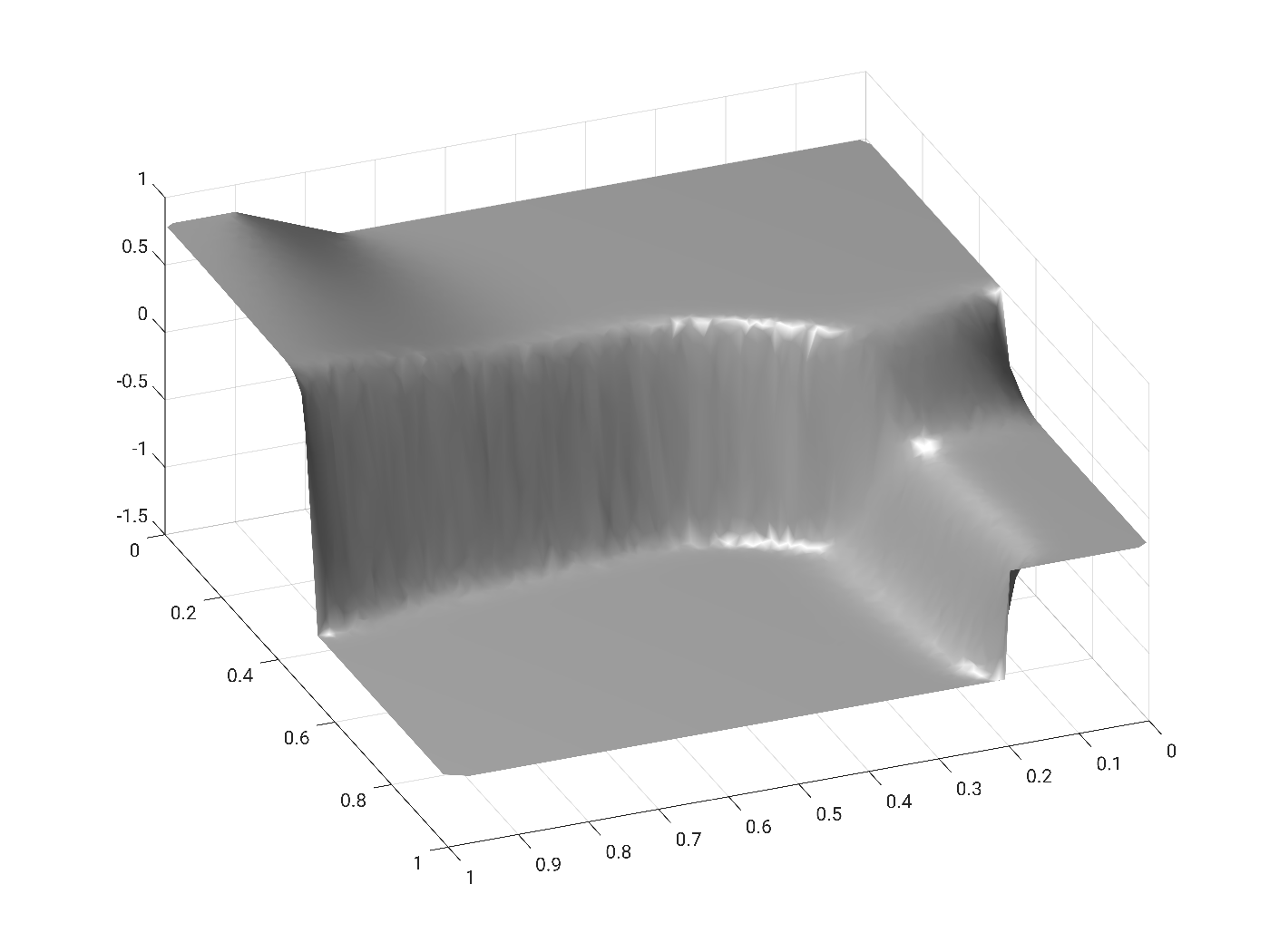}
\includegraphics[scale=.339]{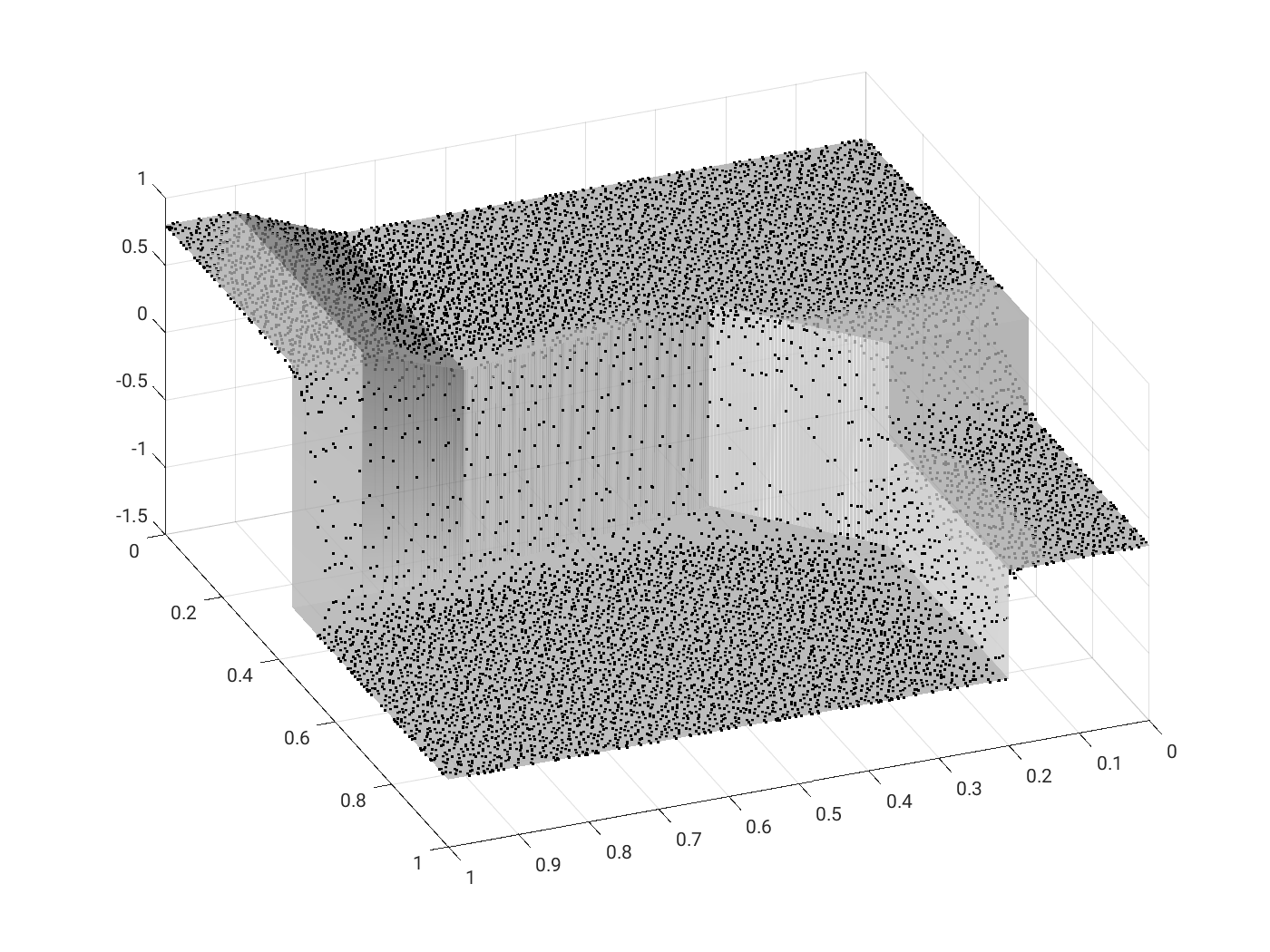} 
}
\subfigure[Adaptive viscosity (Algorithm~\ref{adaptnu}).]%
{
\includegraphics[scale=.339]{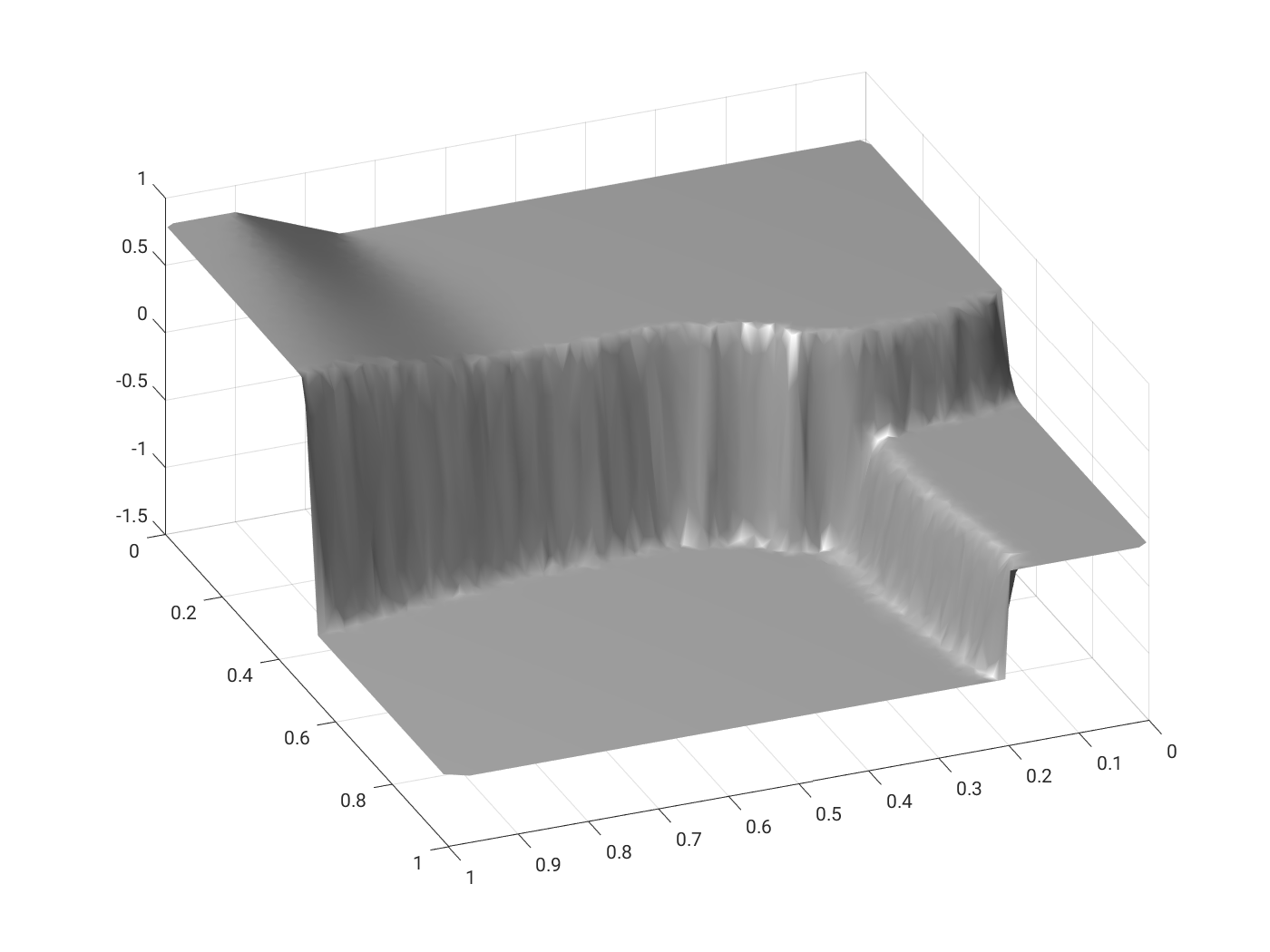}
\includegraphics[scale=.339]{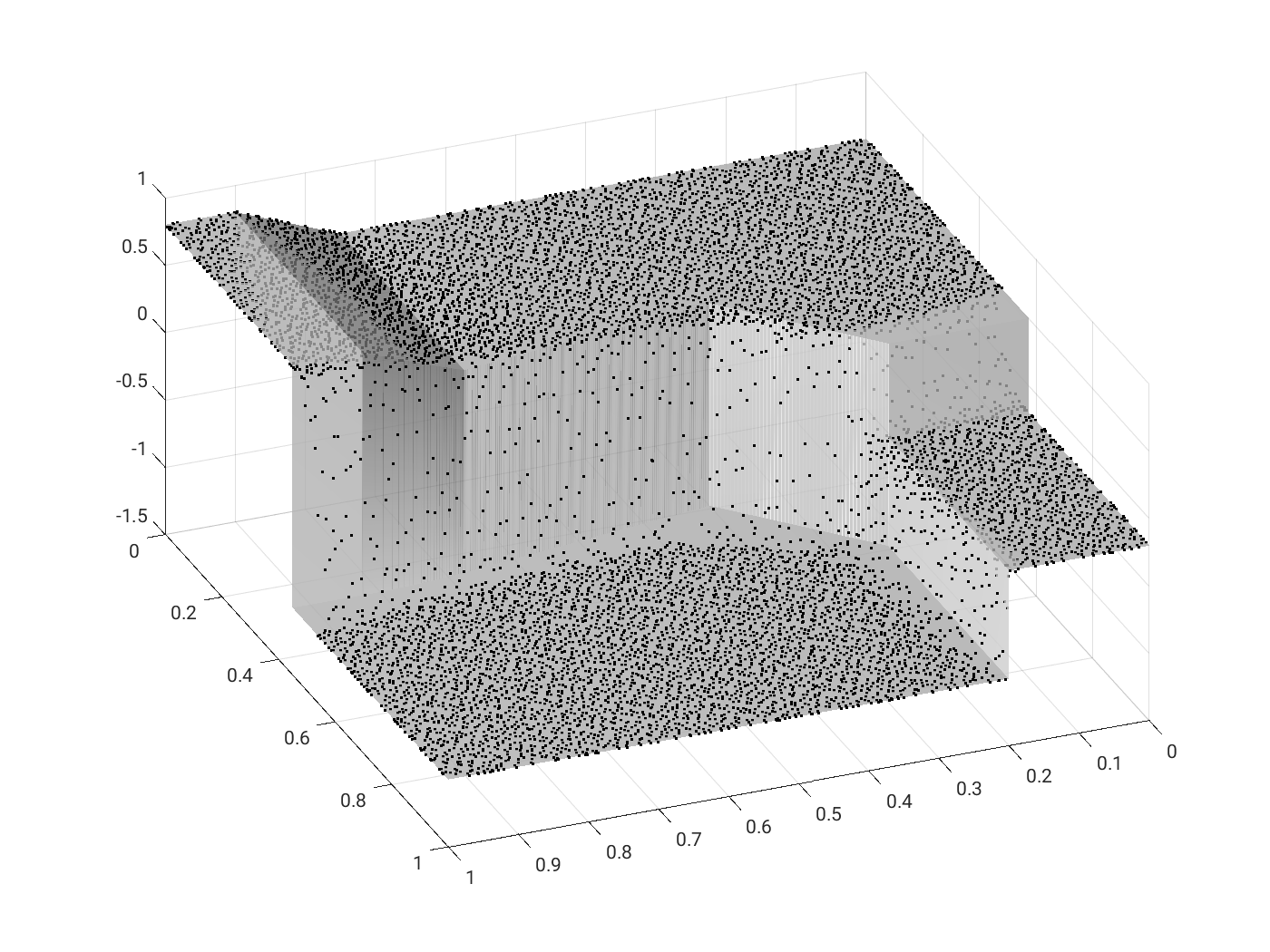}
}
\caption{Example 1: Visual comparison of numerical solutions for $T=0.5$.
\blue{For each algorithm, the 3D visualization of the discrete solution (left) is shown 
together with the point cloud-based scatter plot over a transparent exact solution
(right).} }
\label{fig:ex2b}   
\end{figure}

We choose spacing $h=0.01$ and generate $10000=1/h^2$ Halton points in the unit square, \blue{from which we remove
all points with distance less than $0.25h$ from the boundary, while we add the points obtained by projecting to the 
boundary  the Halton  points at distance less than $h$ from it. This makes a rather irregular distribution of the boundary
nodes, as seen in the right plots of Figure~\ref{fig:ex2b}.}

At each time step we identify inflow boundary nodes as those $\x_i\in X\cap\partial\Omega$, for which the approximated flow
direction $\boeta_i=\F' (U(t,\x_i))$ points inside $\Omega$. For these nodes, we set  $U(t+\Delta t,\x_i):=u(t+\Delta t,\x_i)$
according to the above formula for the exact solution. For all other boundary and interior nodes the value $U(t+\Delta
t,\x_i)$ is computed according to the respective algorithm of  Section~\ref{sec:algorithm} in use. Since the minimization
problem \eqref{weiv} turns out to be infeasible even with  $n_{\max}=100$ neighbors for many boundary nodes, we always set
$\mu_i=0$ if $\x_i\in\partial\Omega$ in order to switch off the artificial viscosity for these nodes.  For  certain nodes (fewer than
20) situated in the immediate neighborhood of the inflow boundary, the minimization problem \eqref{weiw} is infeasible, and  we
use minimal numerical differentiation without inequality constraints to compute the weights $w_{ij}$. This  does not generate
instability in our experiments. 
\blue{In Figure \ref{fig:ex2nodes}(b) we show three examples of employed influence sets.}

\begin{figure}[!t]
\begin{center}
\subfigure[Exact solution.]%
{\includegraphics[scale=.30]{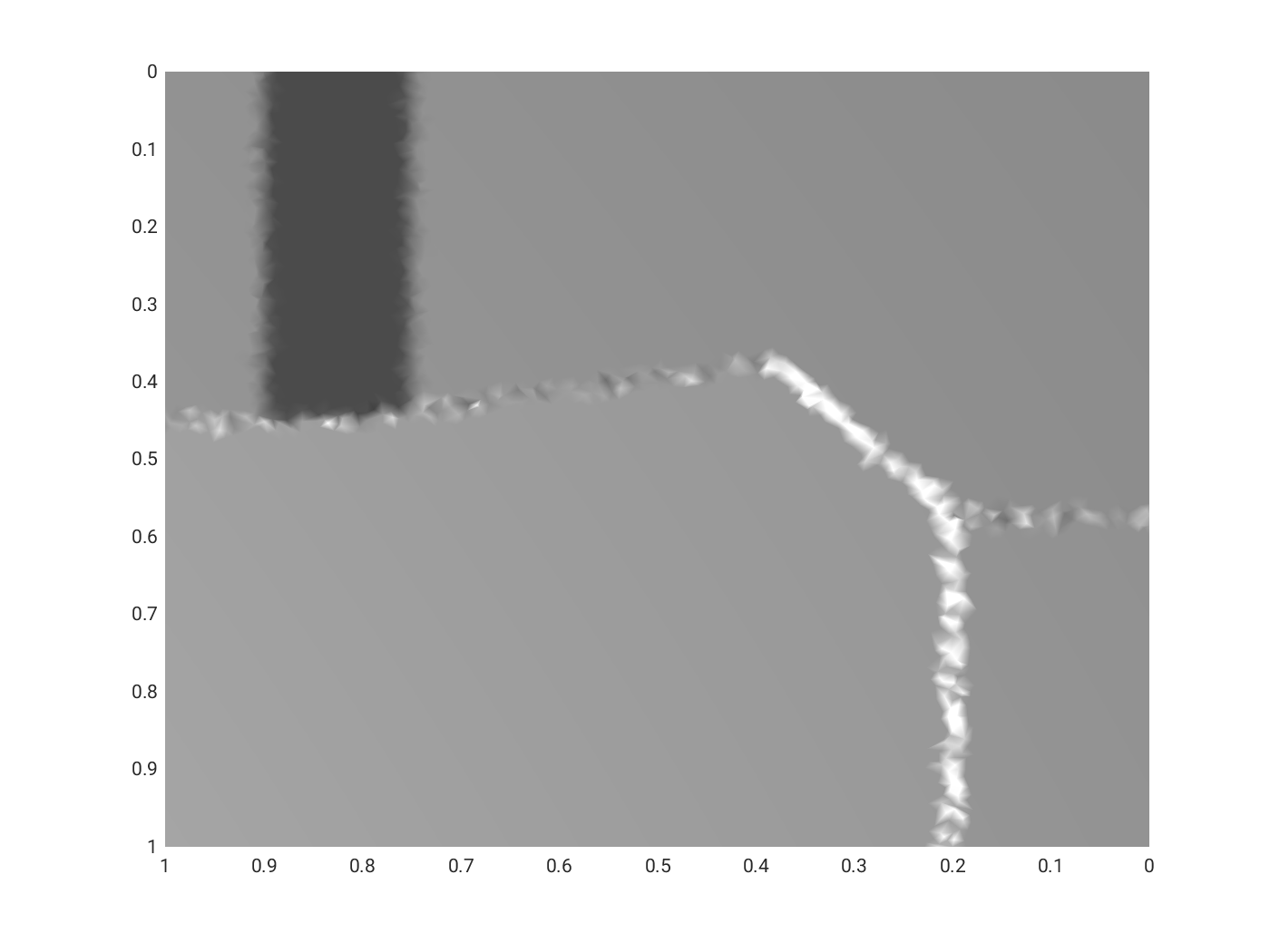}}
\subfigure[No viscosity (Algorithm~\ref{nonu}).]%
{\includegraphics[scale=.30]{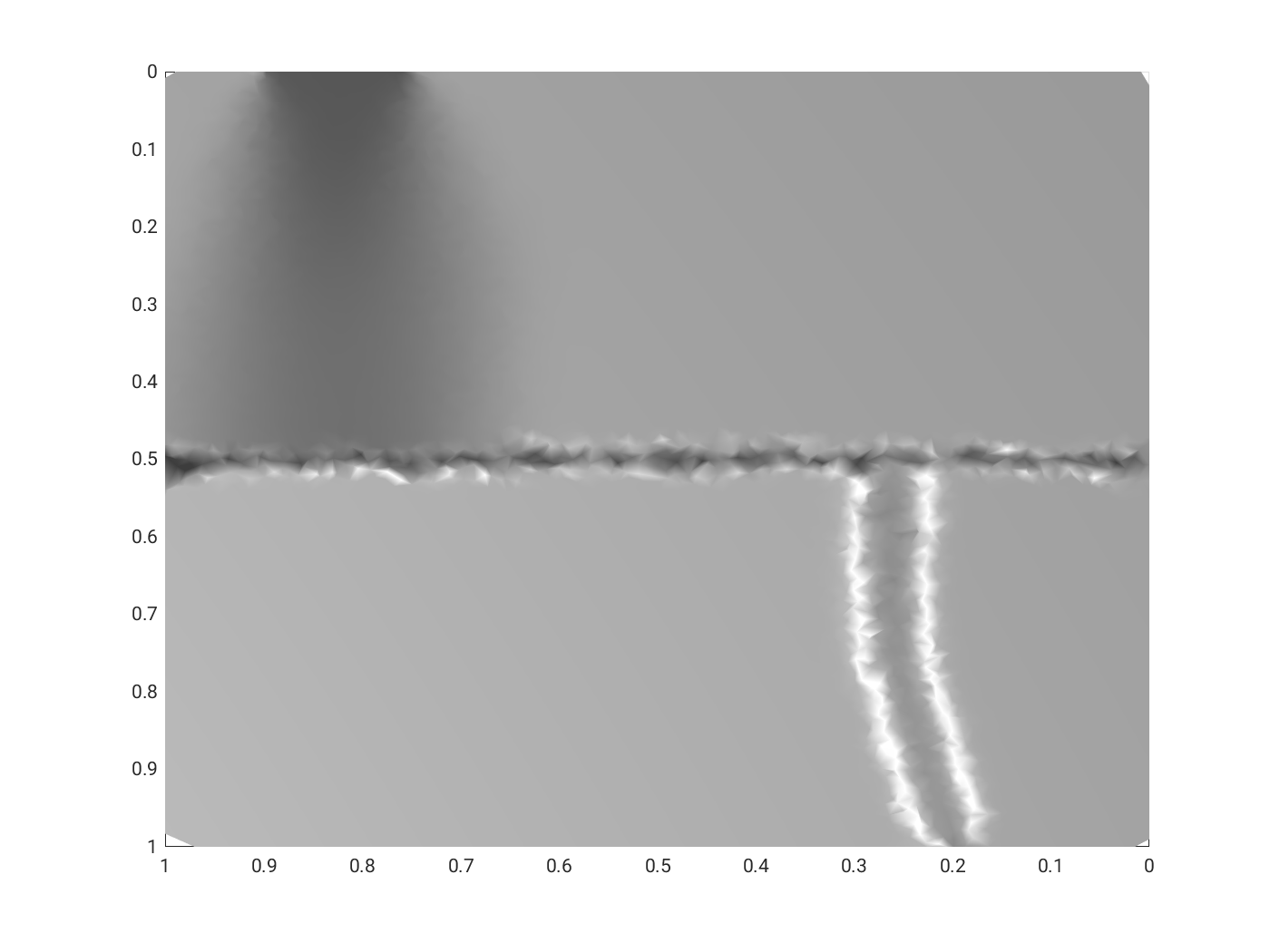}}
\subfigure[Constant viscosity (Algorithm~\ref{constnu}).]%
{\includegraphics[scale=.30]{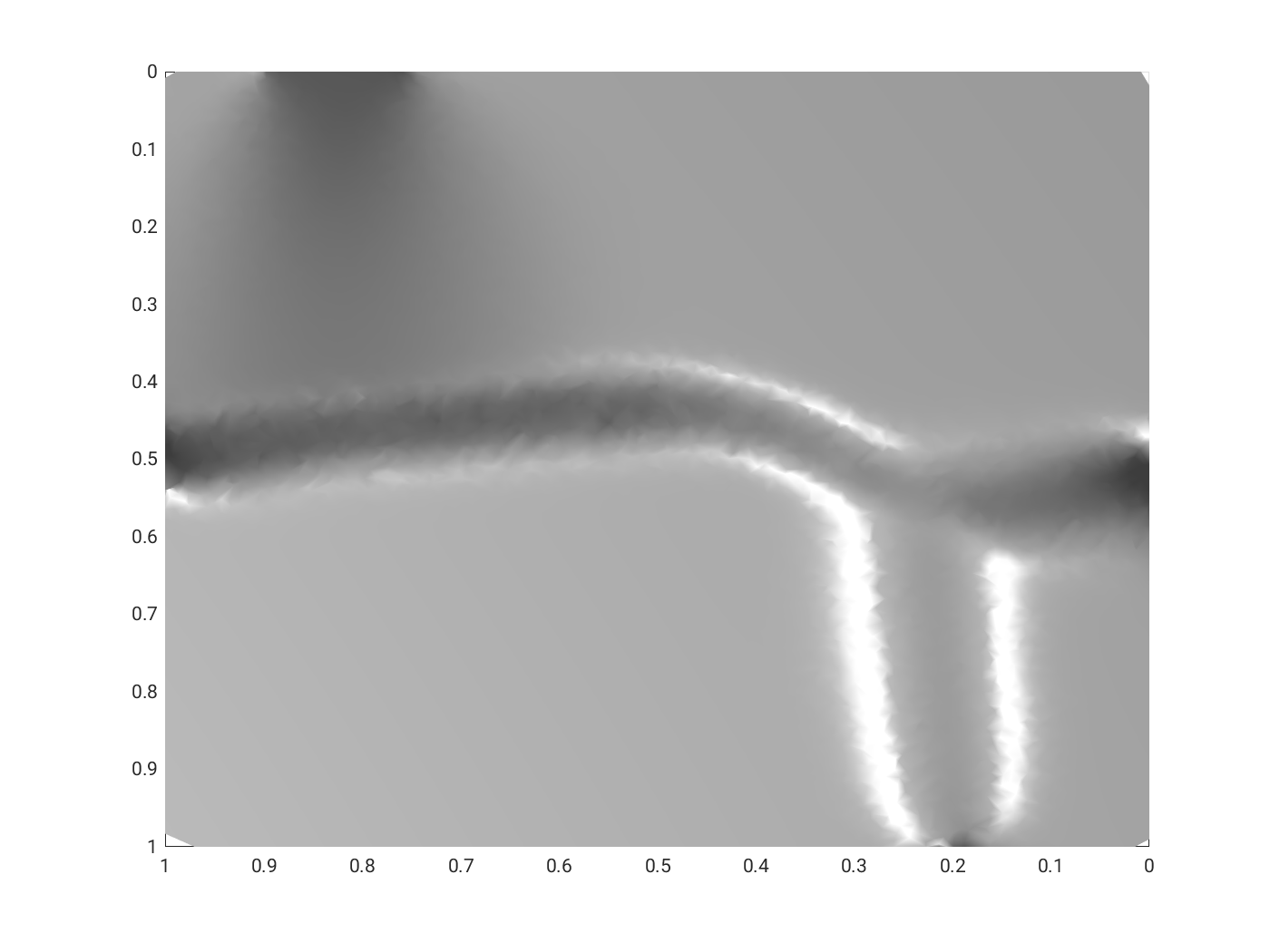}}
\subfigure[Adaptive viscosity (Algorithm~\ref{adaptnu}).]%
{\includegraphics[scale=.30]{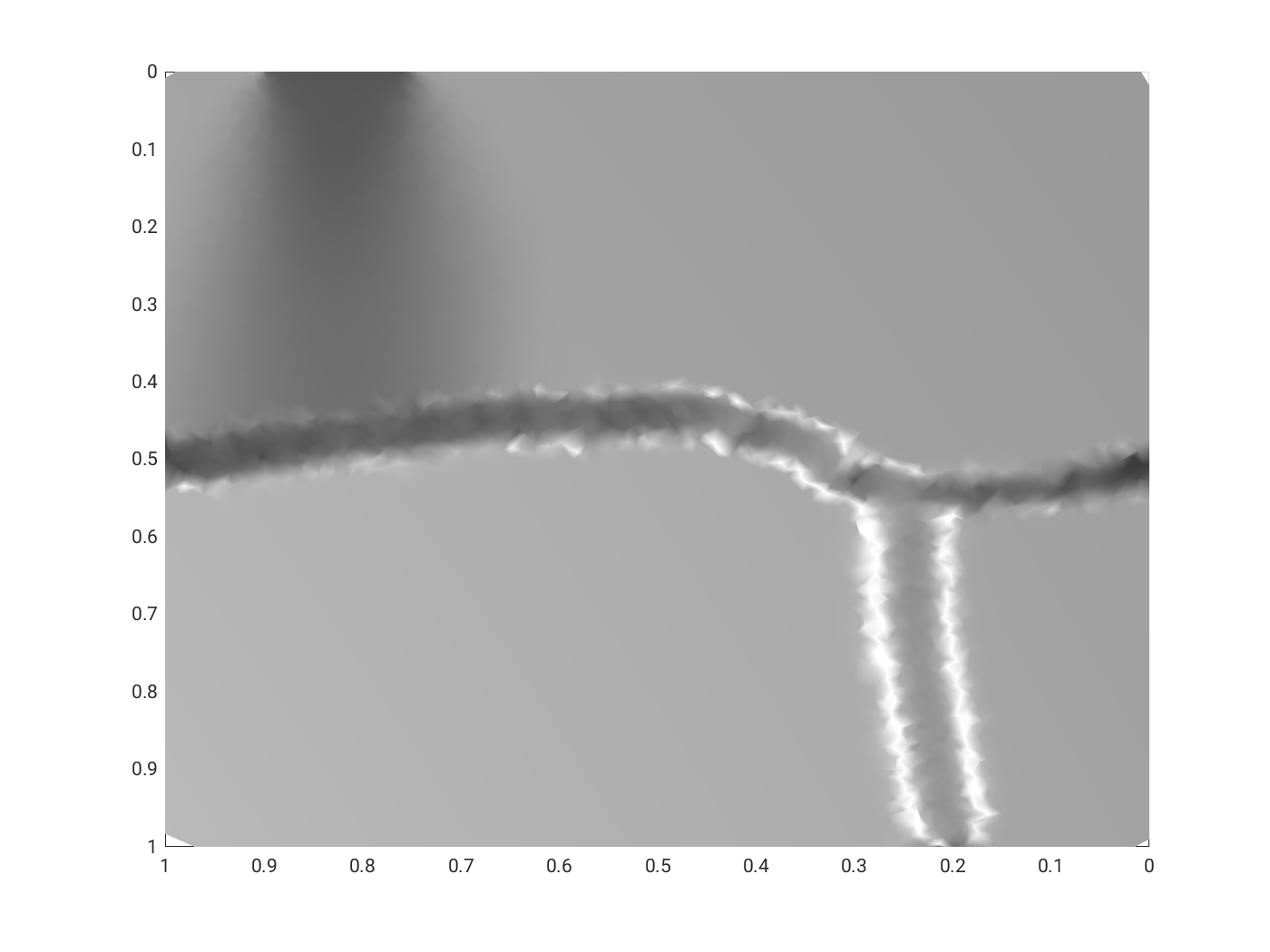}}
\end{center}
\caption{Example 1: Visual comparison of numerical solutions for $T=0.5$: view from above.}
\label{fig:ex2b_above}
\end{figure}

Numerical solutions at time $T=0.5$ obtained by Algorithms~\ref{nonu}--\ref{adaptnu} are presented in Figures~\ref{fig:ex2b}
and \ref{fig:ex2b_above}. \blue{Here, Figure~\ref{fig:ex2b} gives 3D visualizations of the discrete solutions
on the left and point cloud-based scatter plots over a transparent exact solution on the right, 
and Figure~\ref{fig:ex2b_above} provides the view from above, with the plot in \ref{fig:ex2b_above}(a) based on the exact solution on
the same Halton nodes as in Figure~\ref{fig:ex2b}. Note that the wrinkles near discontinuities in the left plots
of Figures~\ref{fig:ex2b} are visualization artefacts, whereas the actual solutions are discrete as shown on the right.} 
We observe that there are no
instabilities in either solution, which is justified by the maximum principle \eqref{locmax} guaranteed for any positive 
scheme. However, the positive scheme without viscosity (Algorithm~\ref{nonu}) completely misses an essential feature of the
solution by straightening out one of discontinuity curves, which confirms that the stability and consistency of the scheme on their own do not guarantee a physically correct
approximation when shocks are present. The shapes are much closer to the exact solution for Algorithms~\ref{constnu} and
\ref{adaptnu}, indicating that the artificial viscosity gives a crucial contribution to the accuracy of the method. 
Moreover, the comparison between the solutions obtained by the positive scheme with constant and adaptive viscosity, see (c)
and (d) in both Figures~\ref{fig:ex2b} and \ref{fig:ex2b_above}, reveals that the adaptive method significantly 
improves the solution in the crucial areas near discontinuities that look much sharper in the plots. In addition we provide  errors $E_1$ and $E_2$ for the three
solutions in Table~\ref{tabex2}. While they are significantly higher for Algorithm~\ref{nonu}, there is little difference 
between Algorithms~\ref{constnu} and \ref{adaptnu}. However, this quantitative comparison is somewhat misleading as it is
slightly in favor of Algorithm~\ref{constnu}, in contrast to the clear visual evidence giving a significant advantage to
Algorithm~\ref{adaptnu}.

\begin{table}[!h]
\centerline{
\begin{tabular}{c|cc} 
\midrule
{method}& $E_1$-error & $E_2$-error\\
\midrule
no viscosity (Algorithm~\ref{nonu}) & 1.18e-01 & 3.89e-01\\
constant viscosity (Algorithm~\ref{constnu}) & 6.43e-02 & 2.38e-01\\
adaptive viscosity (Algorithm~\ref{adaptnu}) & 7.04e-02 & 2.78e-01\\
\midrule
\end{tabular}
}
\caption{Example 1: Errors of solutions for $T=0.5$ obtained by Algorithms~\ref{nonu}--\ref{adaptnu}.}
\label{tabex2}
\end{table}

Note that our solutions are less accurate than those presented in \cite{guermond11}
for comparable node spacing, which is expectable since our methods are of first order, in contrast to the higher order methods employed in that paper.

\subsection{Example 2: Inviscid Burgers' equation with smooth initial condition}\label{ex2}

\begin{figure}[!t]
\begin{center}
\subfigure[Reference solution.]{\includegraphics[scale=.32]{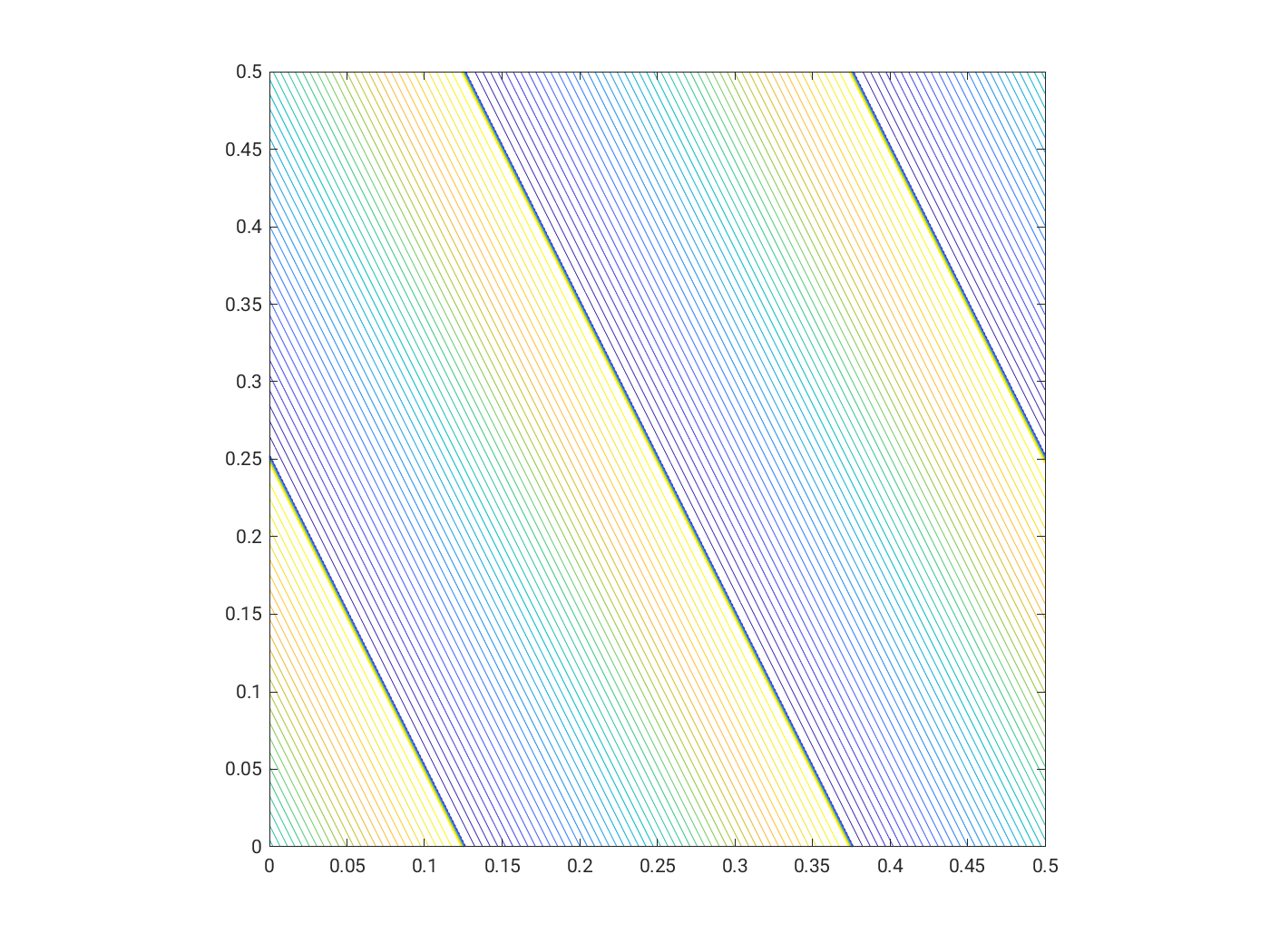}}
\subfigure[Monotone scheme, \cite{hestbook}.]{\includegraphics[scale=.32]{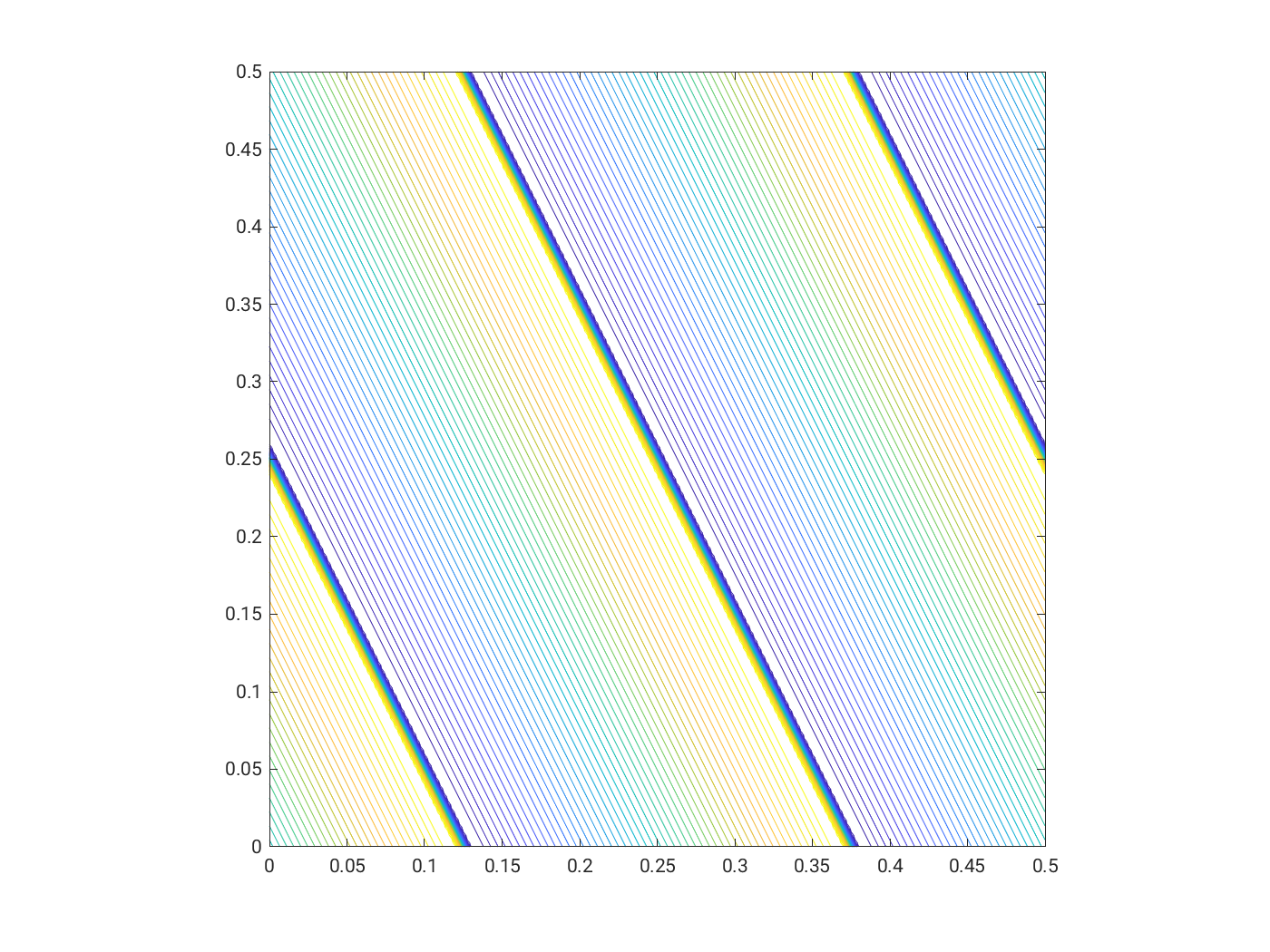}}
\subfigure[Meshless scheme (Algorithm~\ref{adaptnu}) on Cartesian grid.]{\includegraphics[scale=.32]{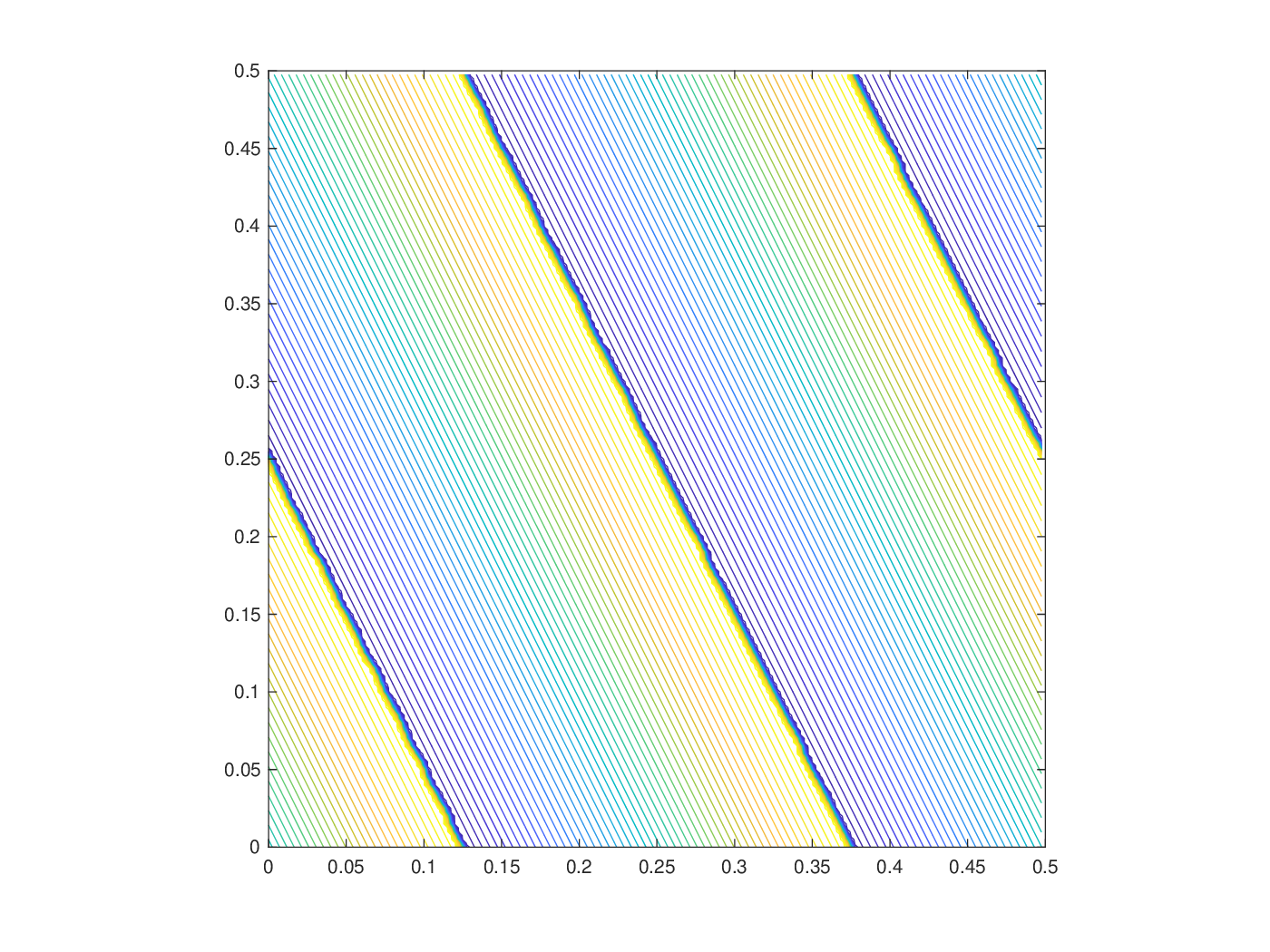}}
\subfigure[Meshless scheme (Algorithm~\ref{adaptnu}) on Halton points.]%
{\includegraphics[scale=.32]{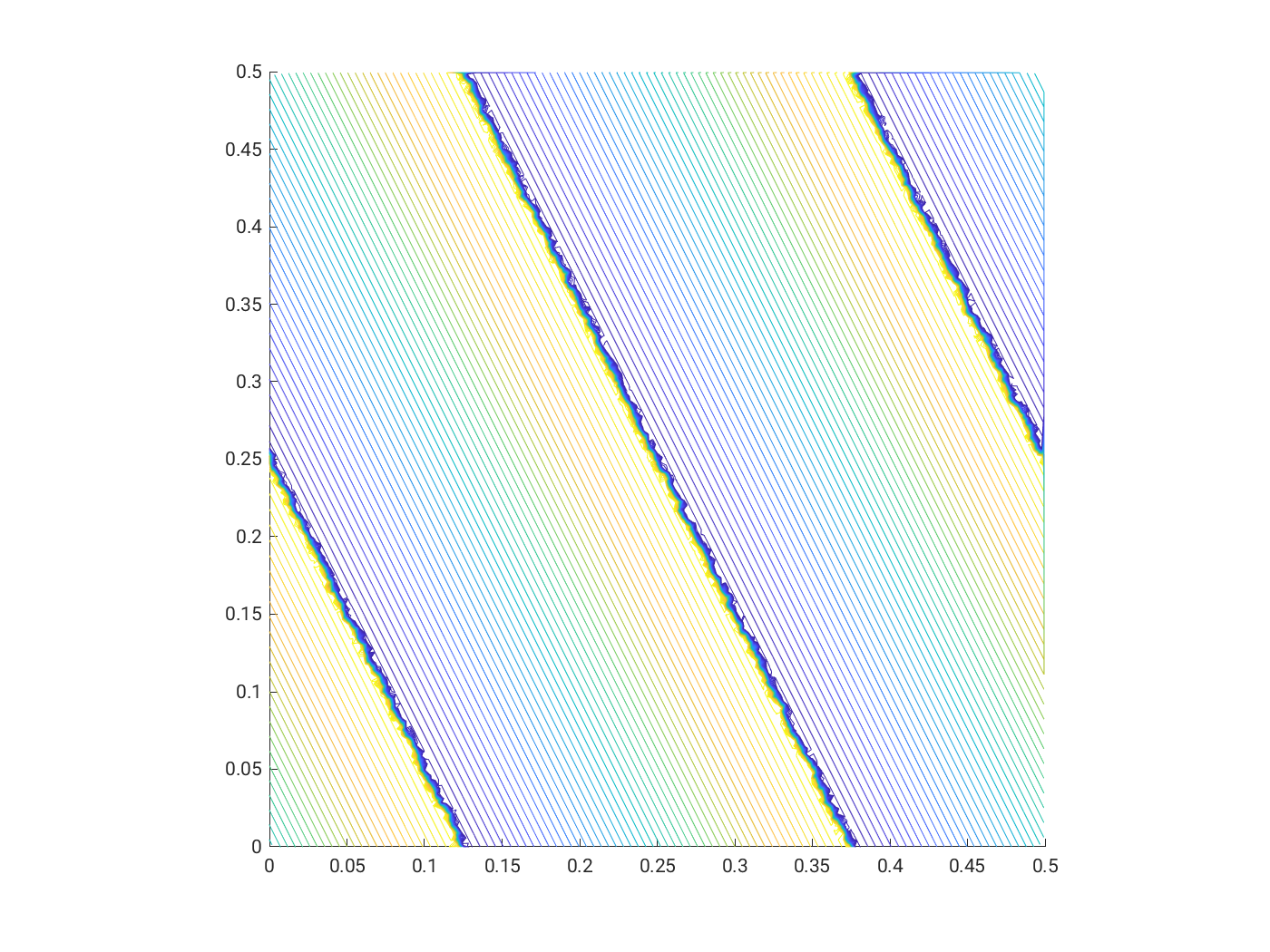}}
\end{center}
\caption{Example 2: Comparison of the contours of the numerical solutions with the reference solution ($T=0.1$).}
\label{fig:ex1cont} 
\end{figure}

In the second example we again consider inviscid Burgers' equation with the flux $\F(u) = \frac{1}{2}u^2\bfv$, 
$\bfv = (1\,,\,1)$, this time with periodic boundary conditions on the spatial domain $\Omega = [0,\,0.5]^2$,
and with the smooth initial condition
$$
u_0(\x) =  \sin(8\pi (x_1+x_2/2)).
$$
In spite of the initial condition being smooth, a shock appears as time advances. Note that this problem is an adaptation 
of the example given in Section 7.1.1 of \cite{hestbook}, where the flux is chosen slightly differently as $\F(u) = u^2\bfv$.
While no explicit exact solution is available, we consider as a reference a high-resolution solution obtained by a fifth order 
WENO scheme on the $400\times400$ Cartesian grid, see  Chapter 11 of the book \cite{hestbook} and
its supplementary code \texttt{BurgersWENO2D.m}. For comparison with our method we choose a first order monotone scheme 
(Script 7.2 in \cite{hestbook} and code \texttt{BurgersM2D.m}) on the $200\times200$ Cartesian grid 
with spacing $h=0.0025$.
 Contour lines of the reference solution and the solution by the first order monotone method at time $T=0.1$ are shown in
Figure~\ref{fig:ex1cont}(a,b).

We apply our positive meshless scheme with adaptive viscosity (Algorithm~\ref{adaptnu}) on the same Cartesian grid,
on a set of 40102 nodes in $[0,\, 0.5)^2$ obtained by generating $0.5^2/h^2=40000$ Halton points in $[0,\,0.5]^2$,
removing all points with distance less than $0.25h$ from the boundary, and then projecting to the lower and left sides of the boundary
the points with distance less than $h$ to them, \blue{and on a set of 40091 random nodes in $[0,\, 0.5)^2$ obtained in a
similar manner, with the components of the nodes from MATLAB's pseudo-random number generator}. 
In all three cases the set of nodes in $[0,\, 0.5)^2$ is extended periodically to a
larger domain, and the sets of influence $X_i$ are computed in the extended domain. The parameters of the scheme are 
chosen in the same way as before, in particular  $\Delta t=0.2h=5\cdot 10^{-4}$ and $\mu = 0.5h=1.25 \cdot 10^{-3}$.

\begin{figure}[!htbp]
\subfigure[Meshless scheme on Cartesian grid.]{\includegraphics[scale=.35]{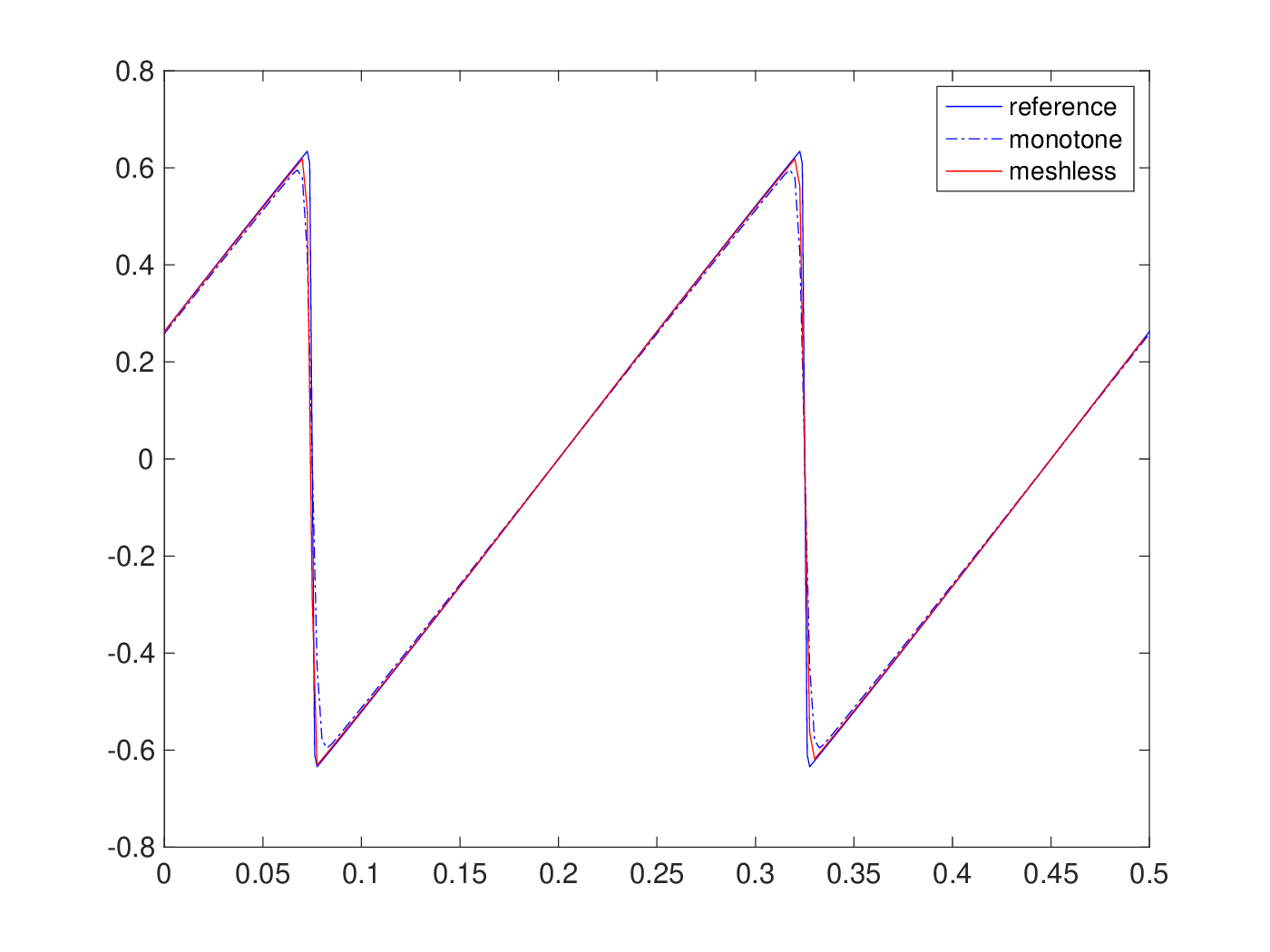}}
\subfigure[Meshless scheme on Halton points.]{\includegraphics[scale=.35]{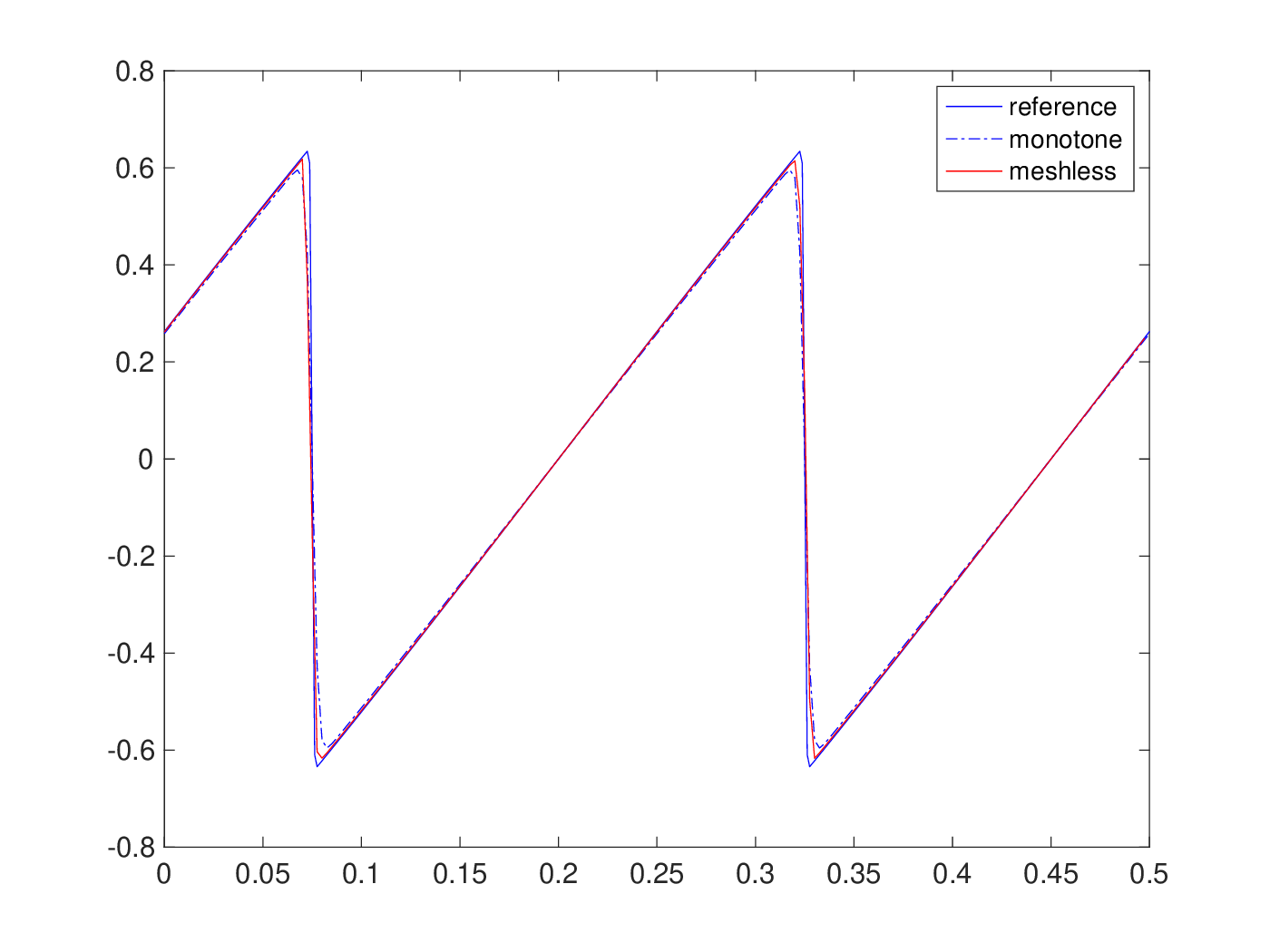}}
\subfigure[Grid: zoom no.~1.]{\includegraphics[scale=.17]{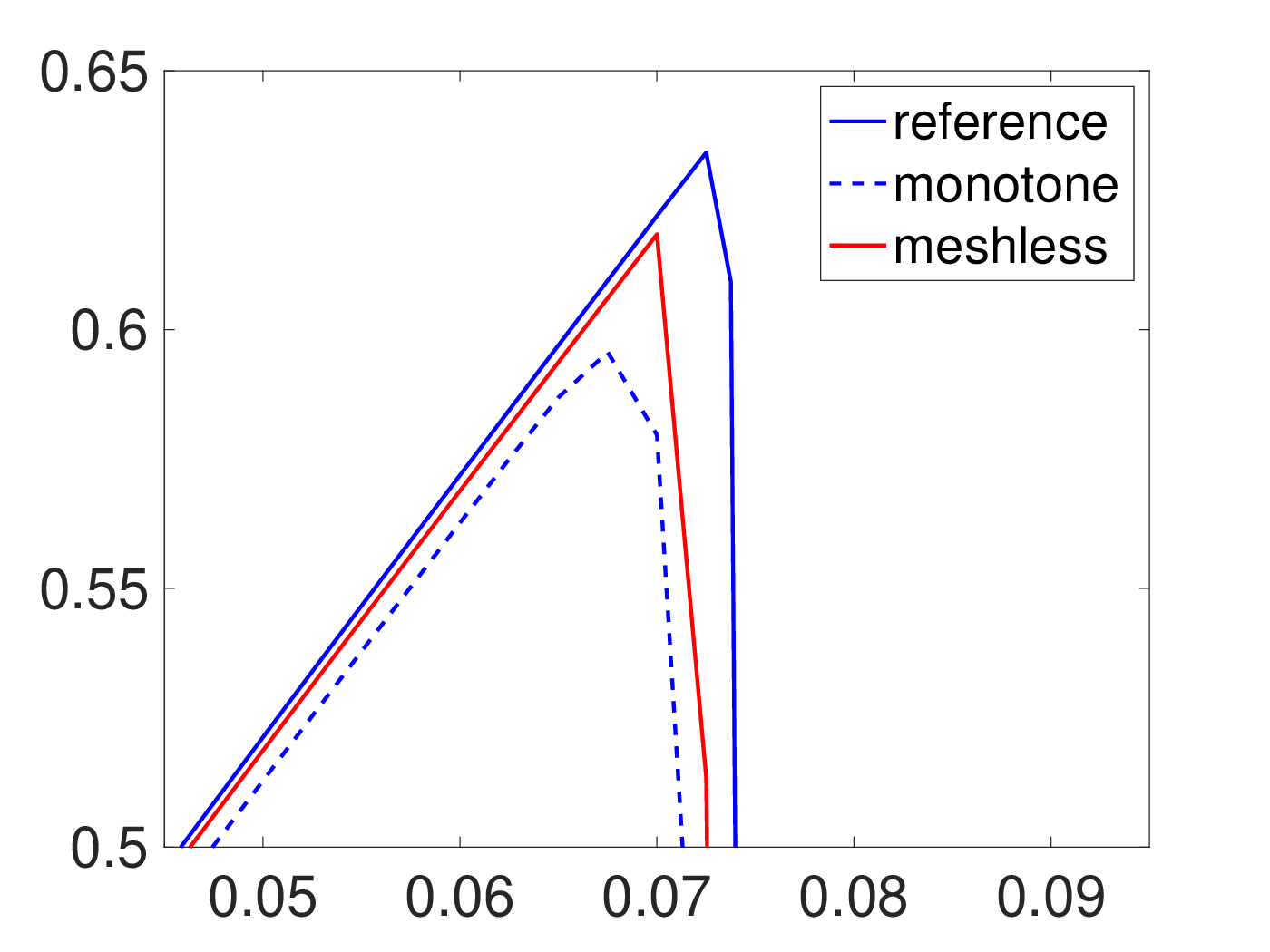}}
\subfigure[Grid: zoom no.~2.]{\includegraphics[scale=.17]{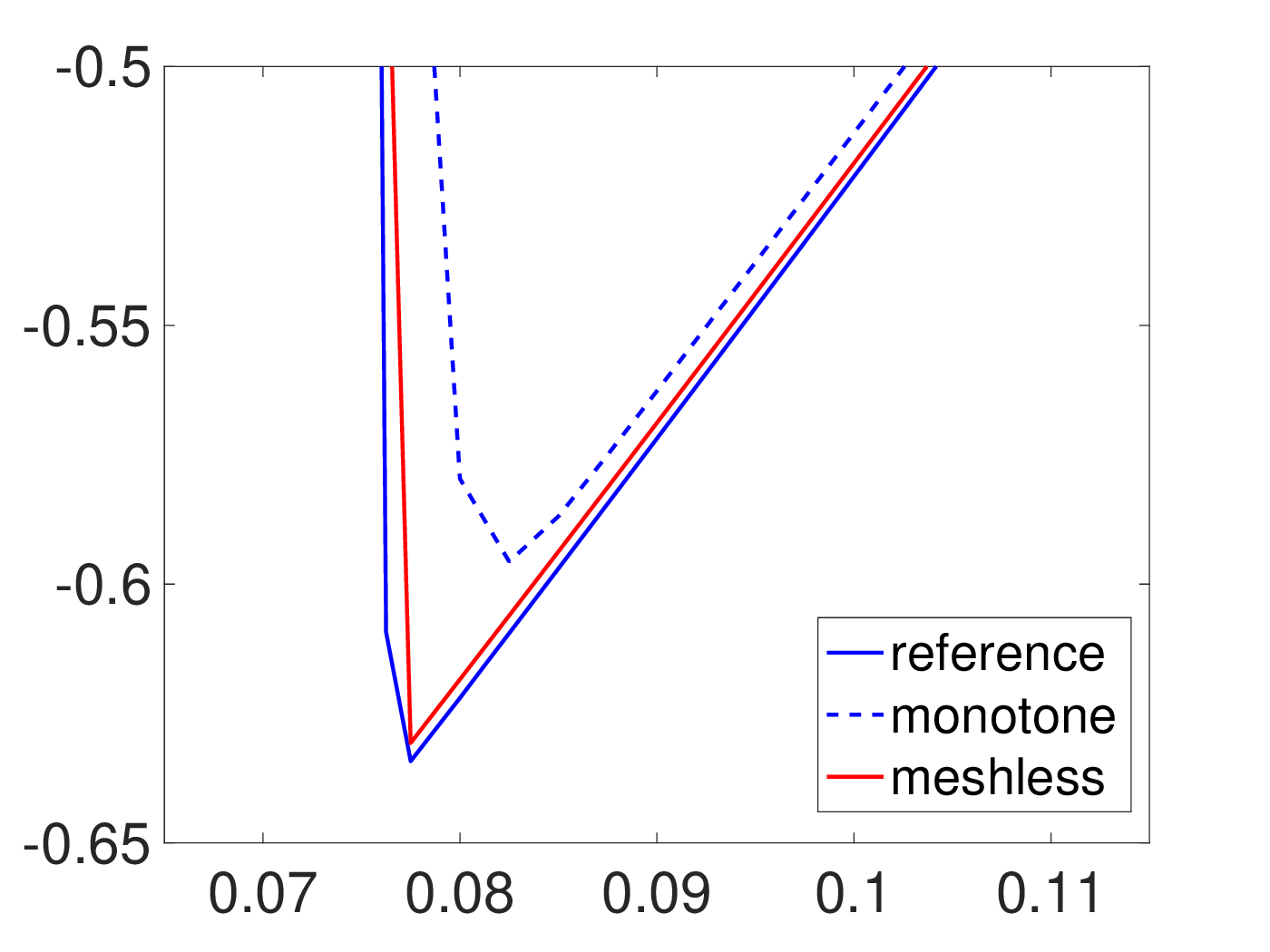}}
\subfigure[Grid: zoom no.~3.]{\includegraphics[scale=.17]{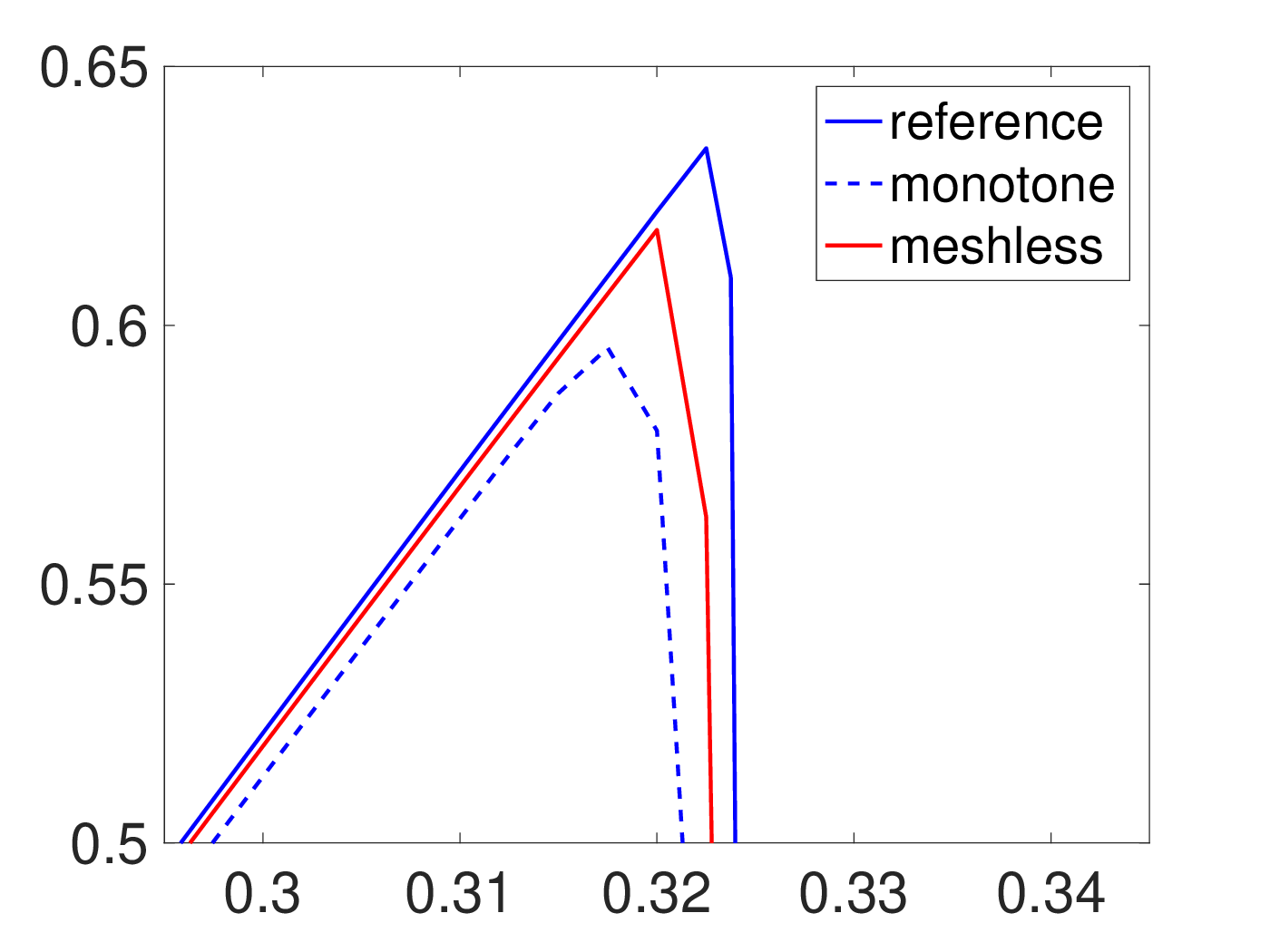}}
\subfigure[Grid: zoom no.~4.]{\includegraphics[scale=.17]{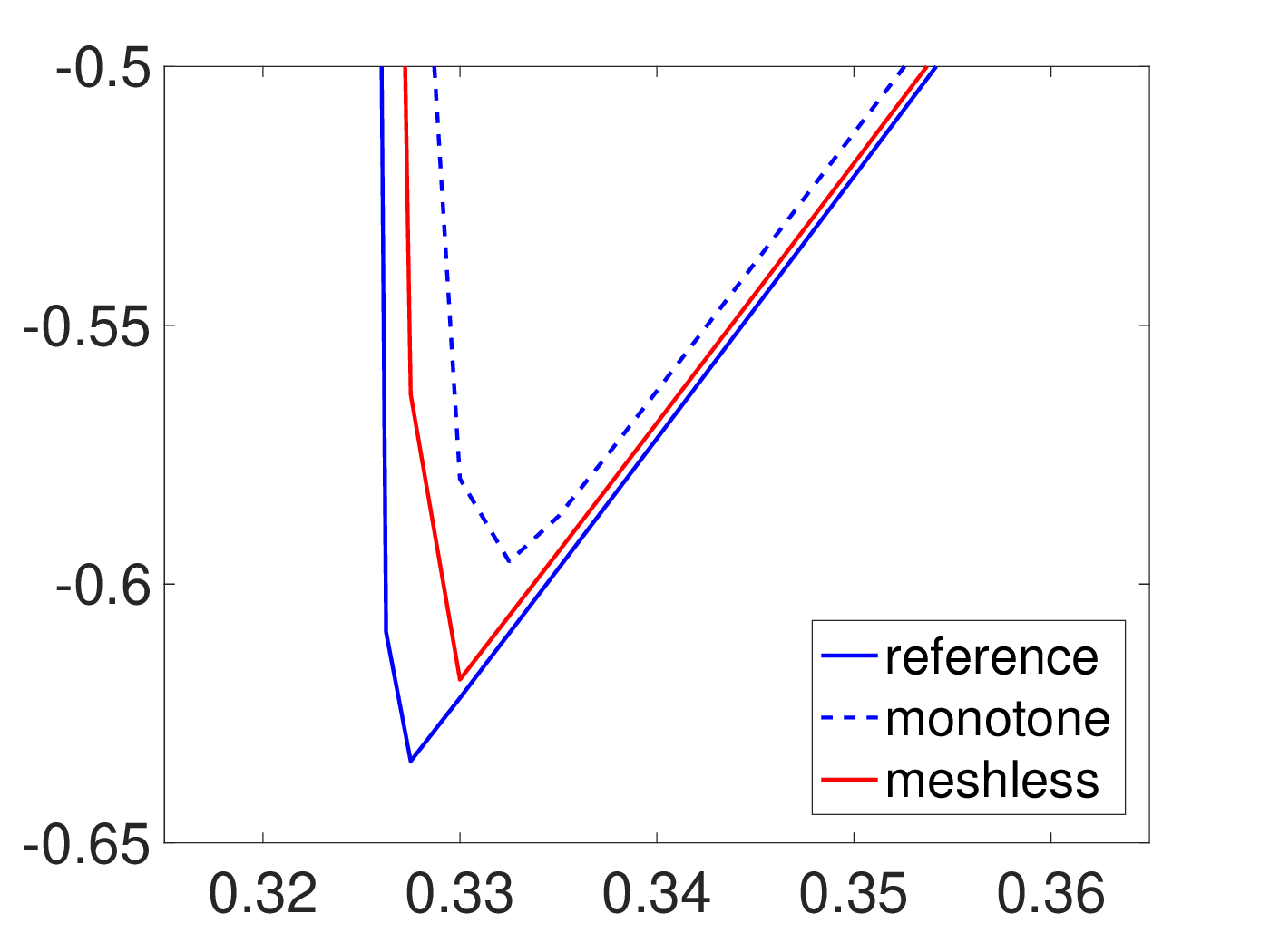}}
\subfigure[Halton: zoom no.~1.]{\includegraphics[scale=.17]{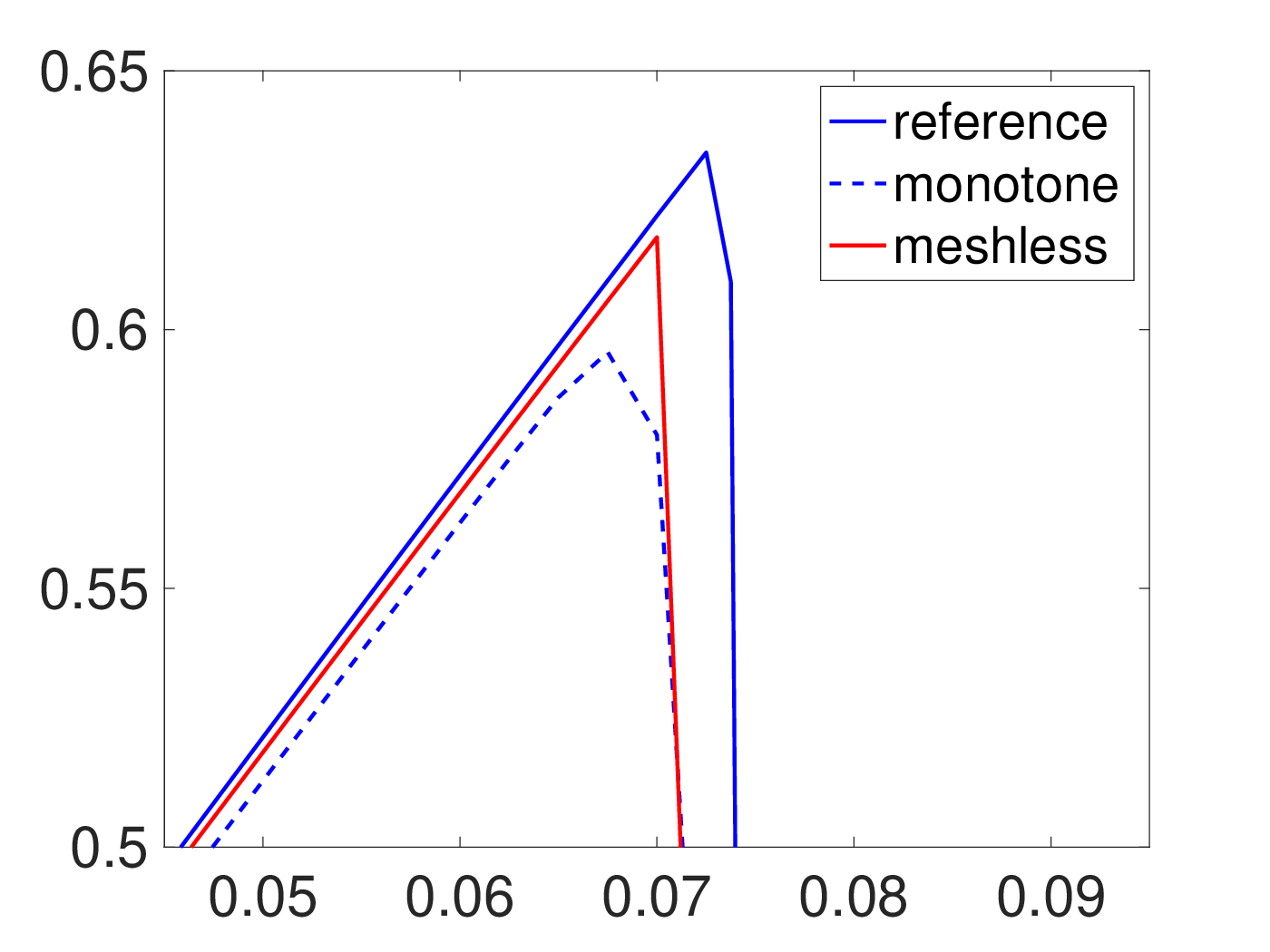}}
\subfigure[Halton: zoom no.~2.]{\includegraphics[scale=.17]{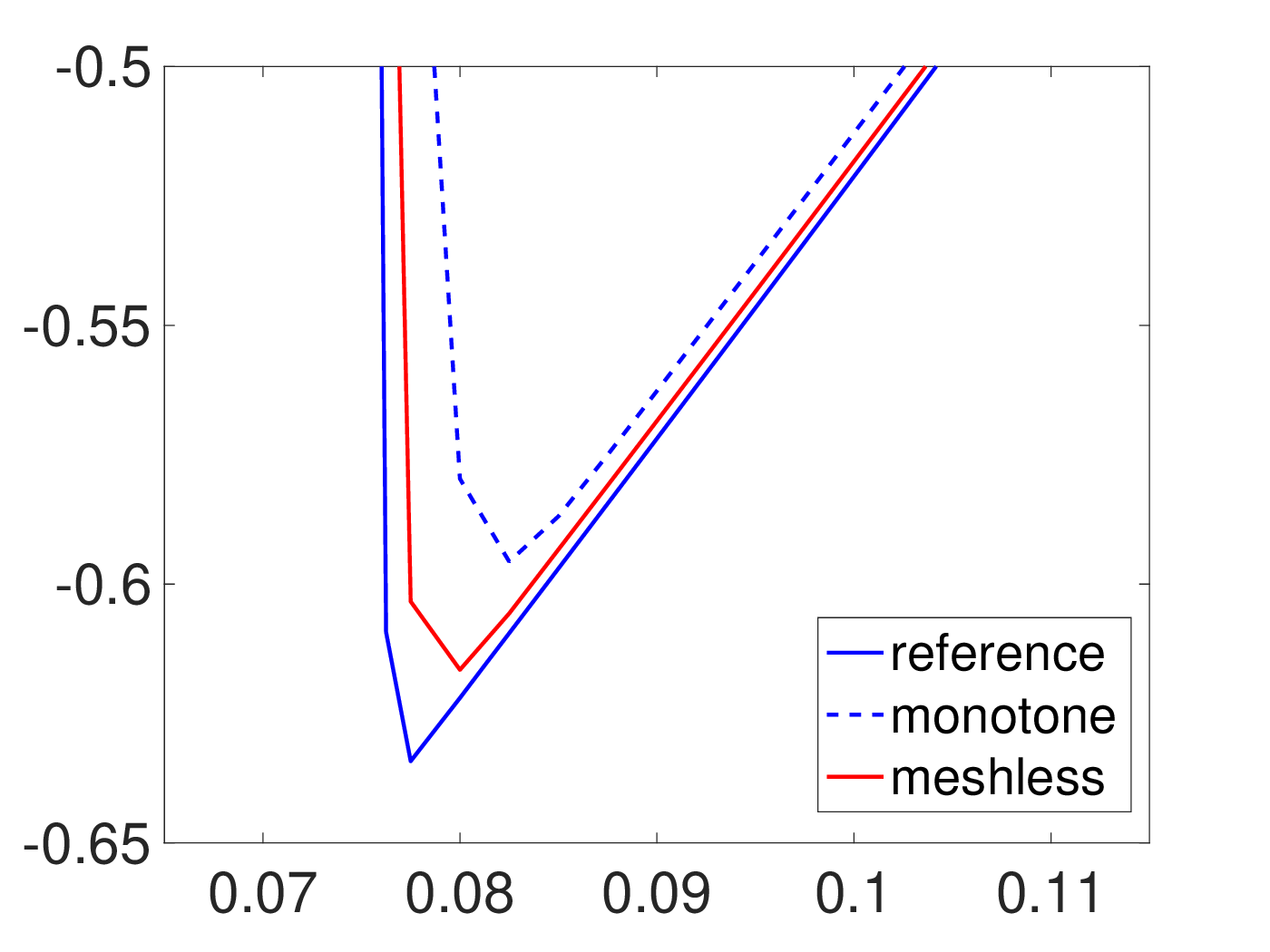}}
\subfigure[Halton: zoom no.~3.]{\includegraphics[scale=.17]{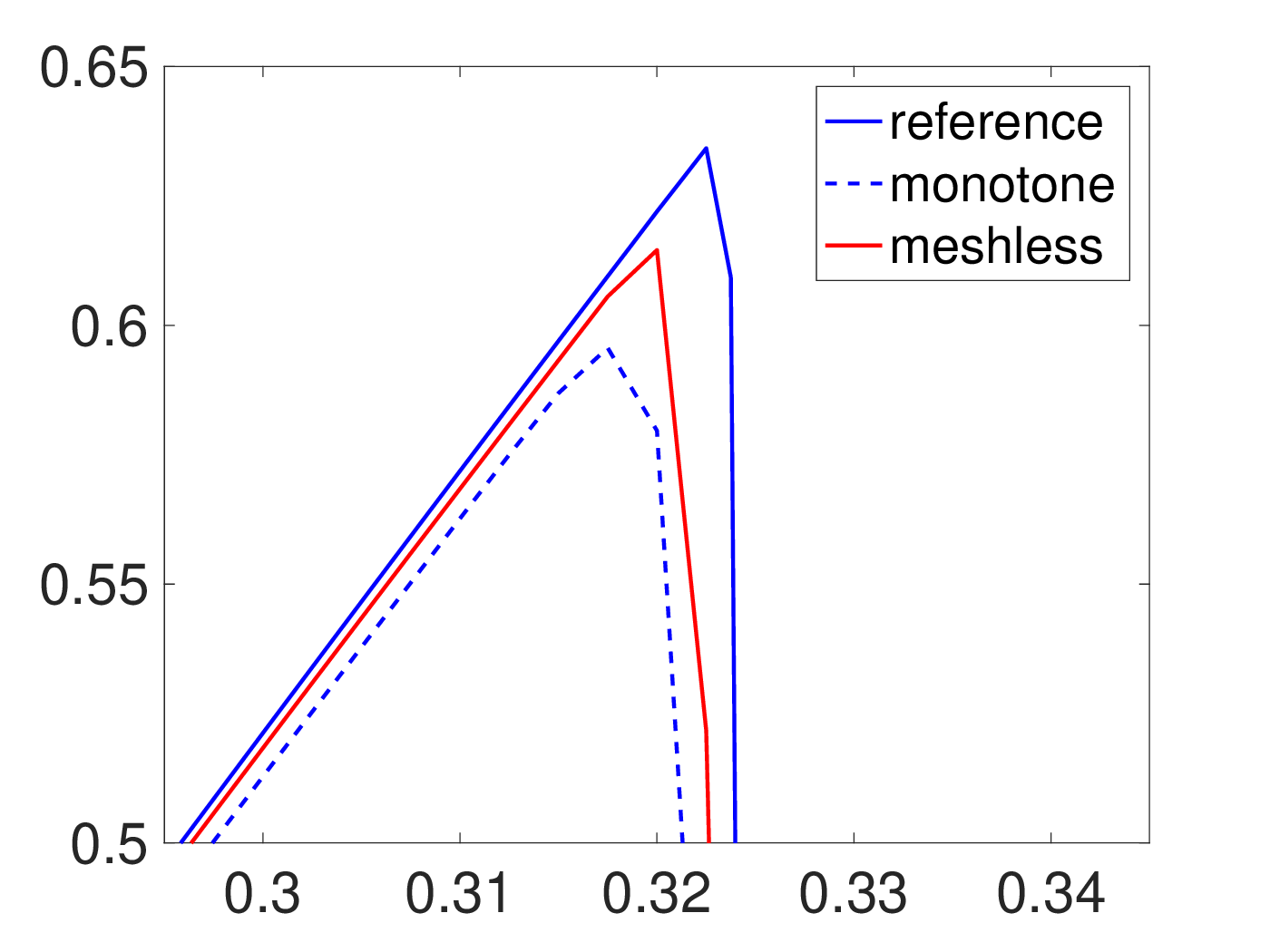}}
\subfigure[Halton: zoom no.~4.]{\includegraphics[scale=.17]{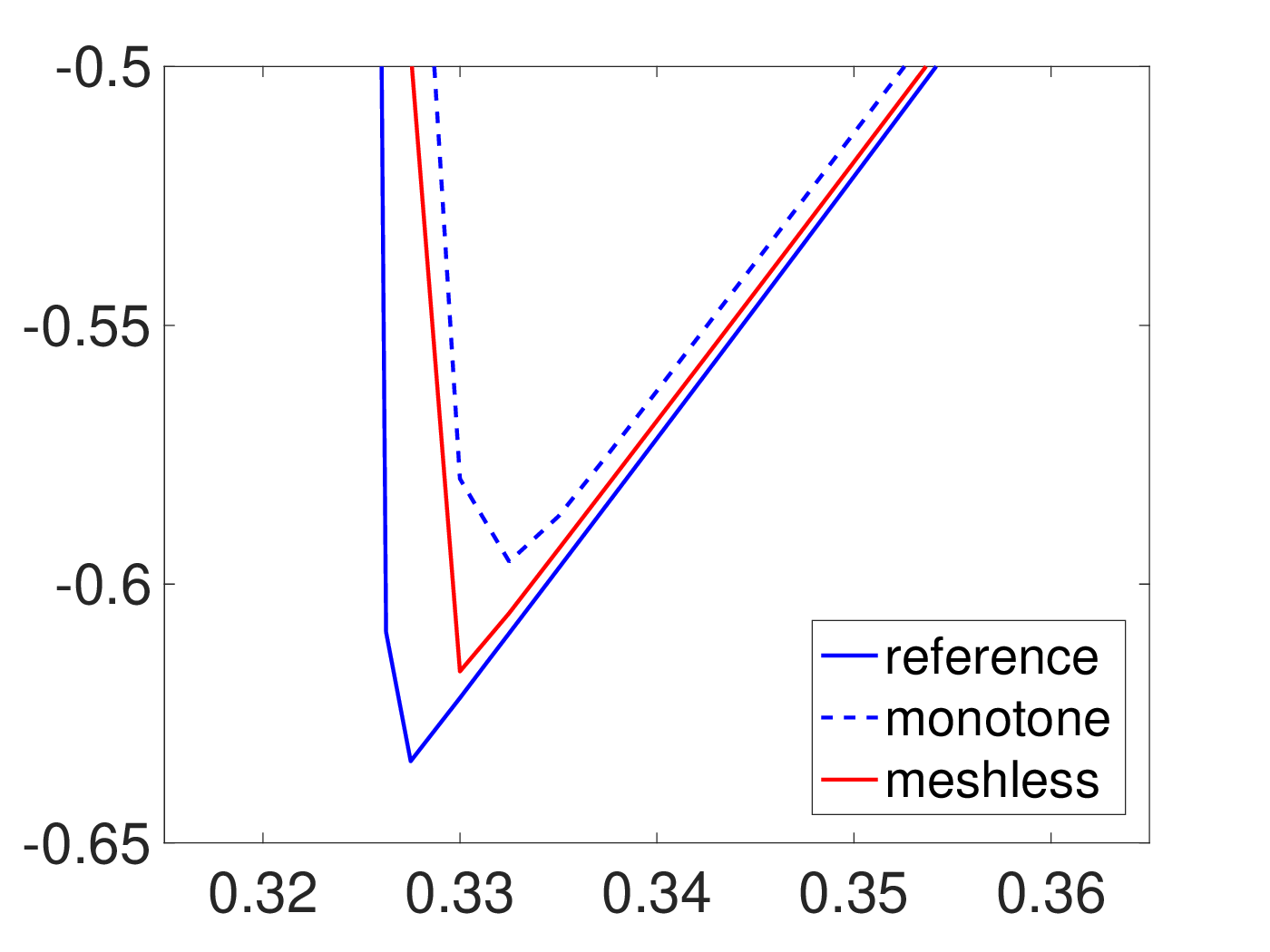}}
\subfigure[Random: zoom no.~1.]{\includegraphics[scale=.17]{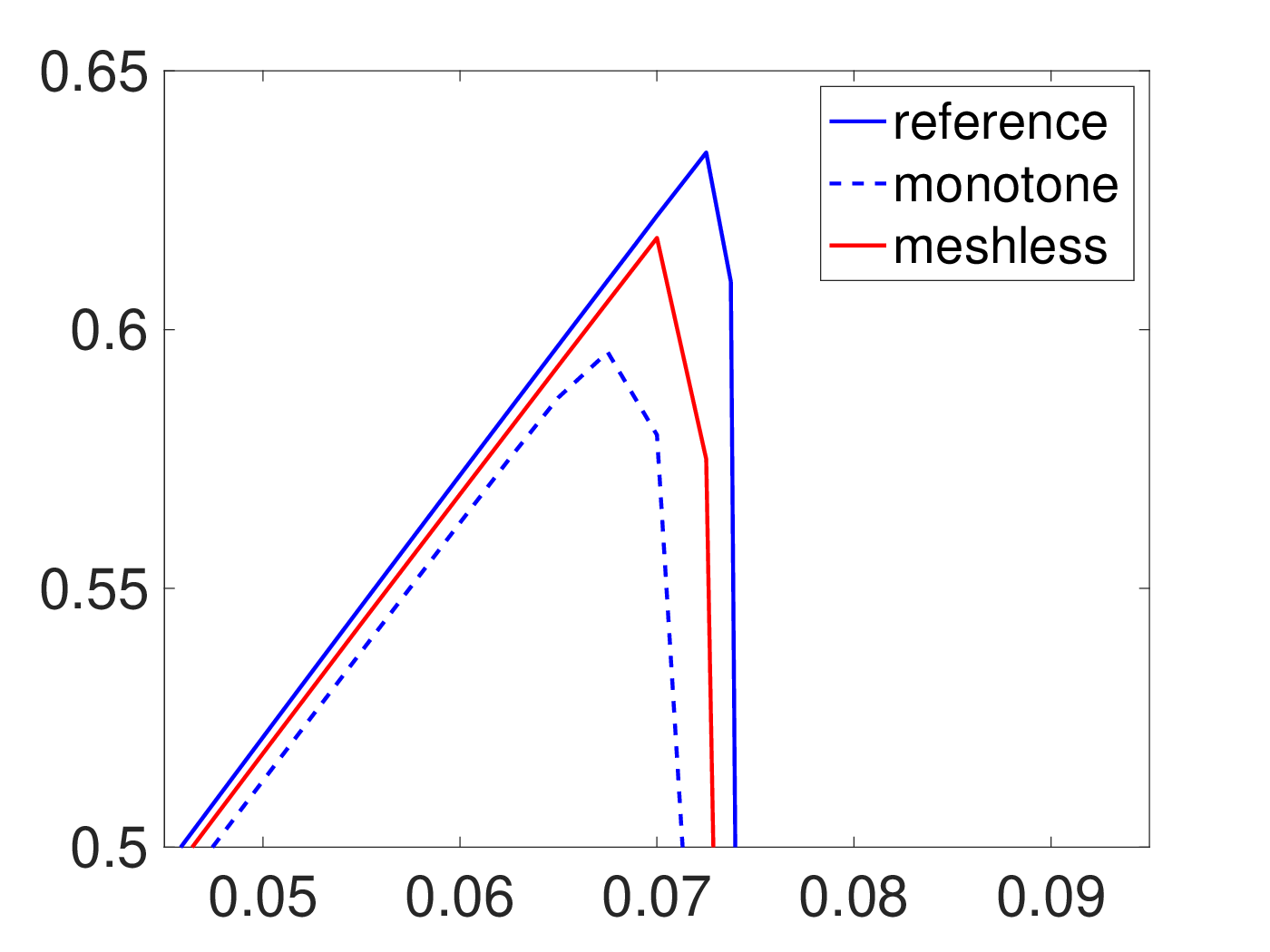}}
\subfigure[Random: zoom no.~2.]{\includegraphics[scale=.17]{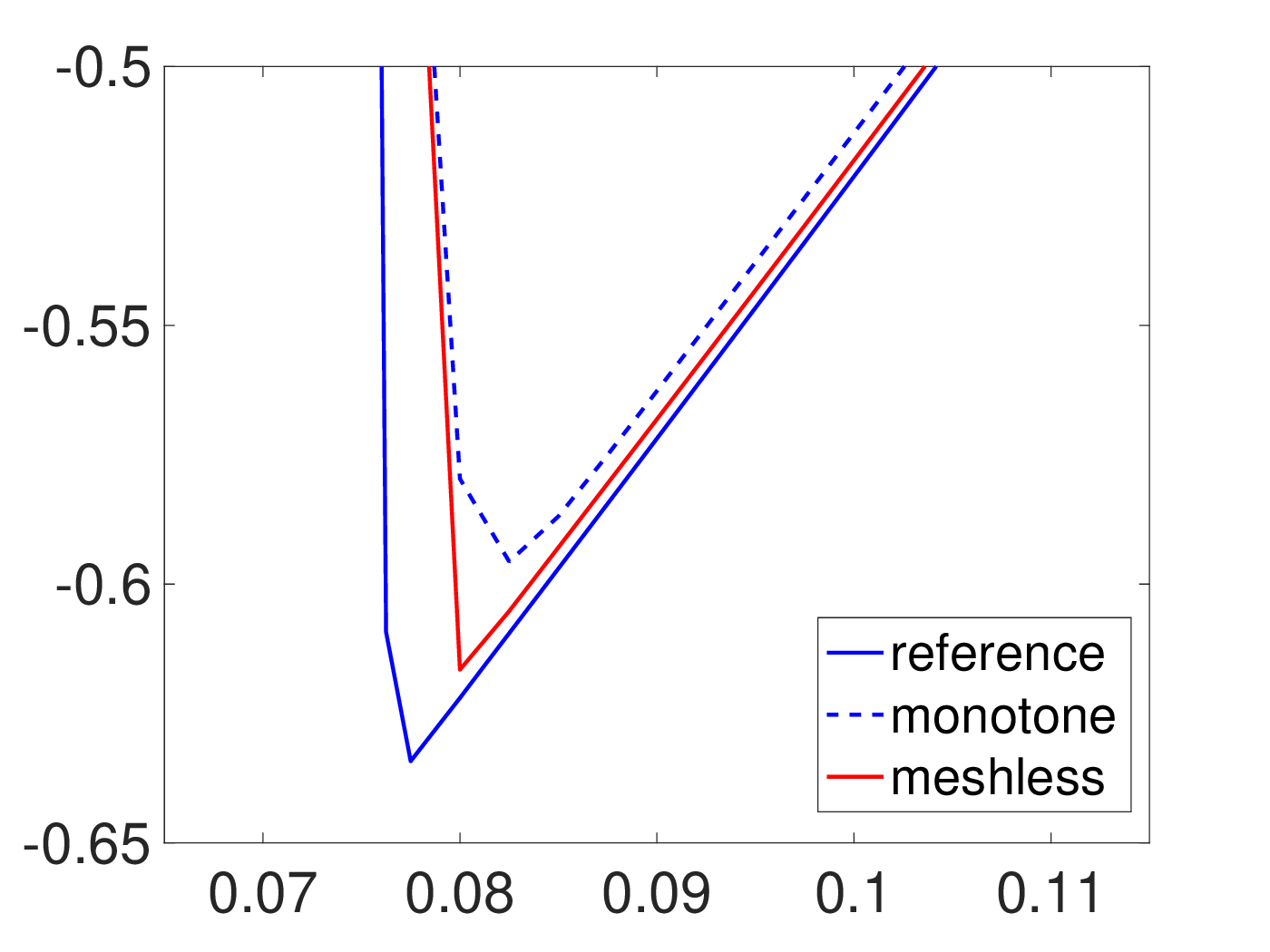}}
\subfigure[Random: zoom no.~3.]{\includegraphics[scale=.17]{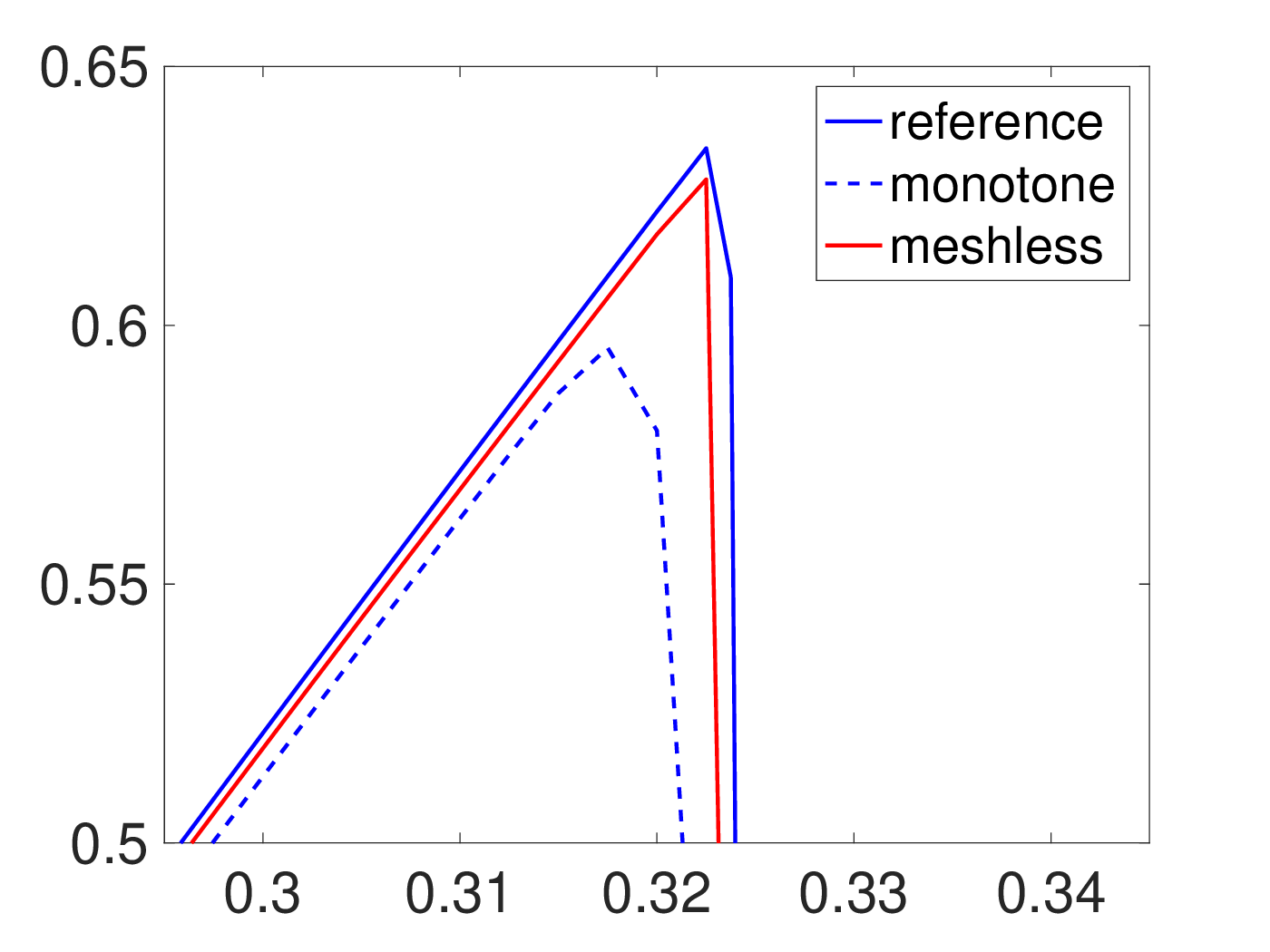}}
\subfigure[Random: zoom no.~4.]{\includegraphics[scale=.17]{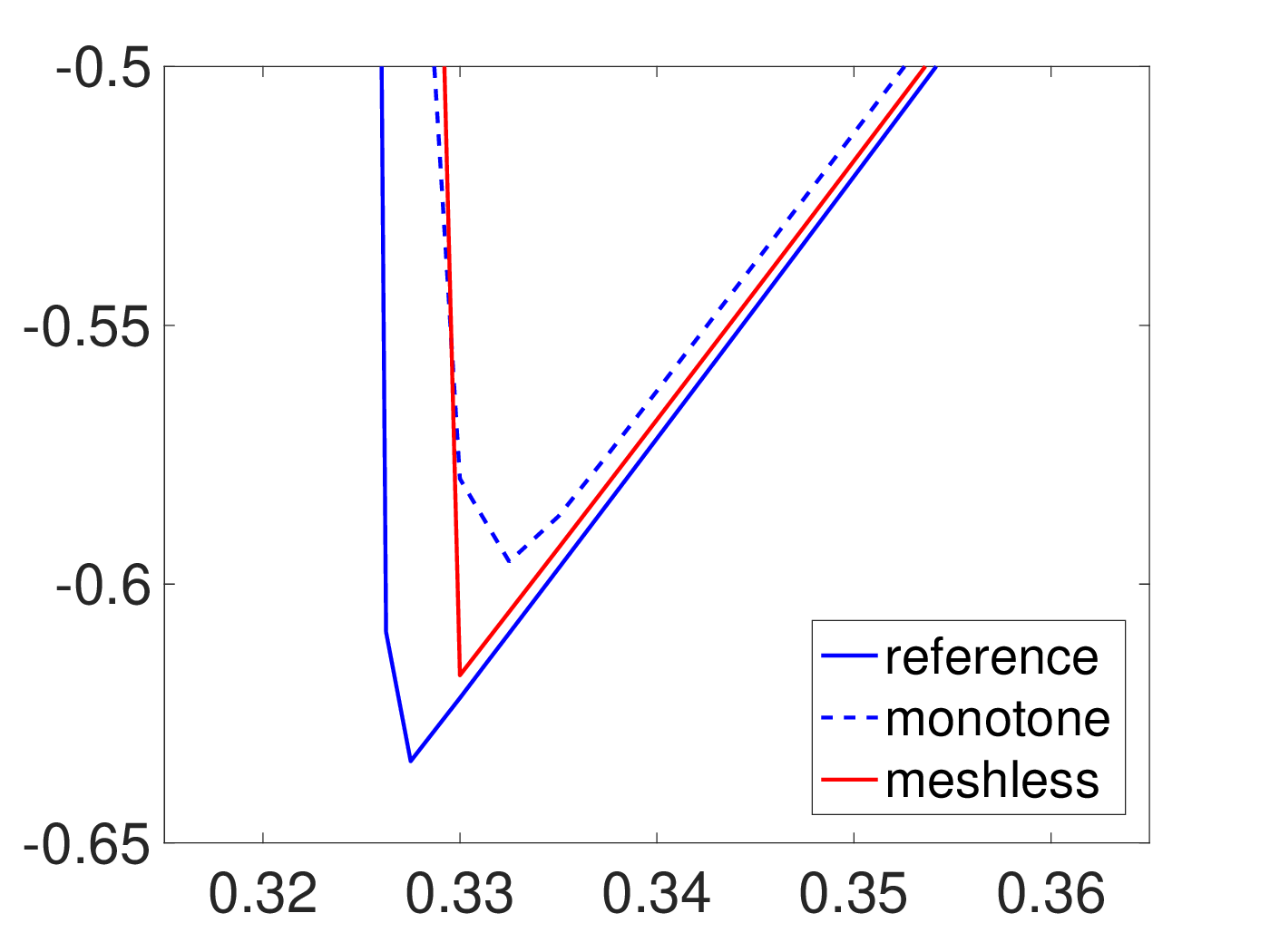}}
\caption{Example 2: Comparison of the cross sections at $x_2= 0.1$ of 
the reference solution and  the monotone scheme with  the meshless scheme of  Algorithm~\ref{adaptnu} on (a) the Cartesian grid  and (b) 
Halton points ($T=0.1$). Subfigures (c)--(n) show zoom-in plots of four critical portions of the cross sections 
near the shocks, comparing three meshless solutions (on the Cartesian grid, Halton and random points) 
to the reference solution and the solution by the monotone scheme.}
\label{fig:ex1cross} 
\end{figure}

The contour lines of the solutions at time $T=0.1$ obtained by our method on the grid and Halton nodes are presented in Figure~\ref{fig:ex1cont}(c,d). In
addition, we  compare the cross sections of all five solutions for $x_2= 0.1$ in Figure~\ref{fig:ex1cross}, and the errors
with respect to the reference solution in Table~\ref{table:tabex1}. Note that in order to compute the cross sections and
contour lines  in the case of Halton or random nodes we interpolate the solution on the $200\times200$ Cartesian grid by the piecewise
linear polynomial on the Delaunay triangulation of the nodes. The results of the meshless scheme on the
grid, Halton and random nodes
are very similar, and significantly better than those of the monotone scheme. This is clear from the errors in
the table and zoom-in plots in Figure~\ref{fig:ex1cross}(c-n). In particular, these results confirm remarkable robustness 
of our method with respect to the use of scattered points. Note that  the contours in Figure~\ref{fig:ex1cont}(d)
appear wobbling at the shocks, due to the irregular distribution of the Halton nodes in close proximity of the shock
lines.

\begin{figure}[!t]
\subfigure[Spacing 0.01.]{\includegraphics[scale=.17]{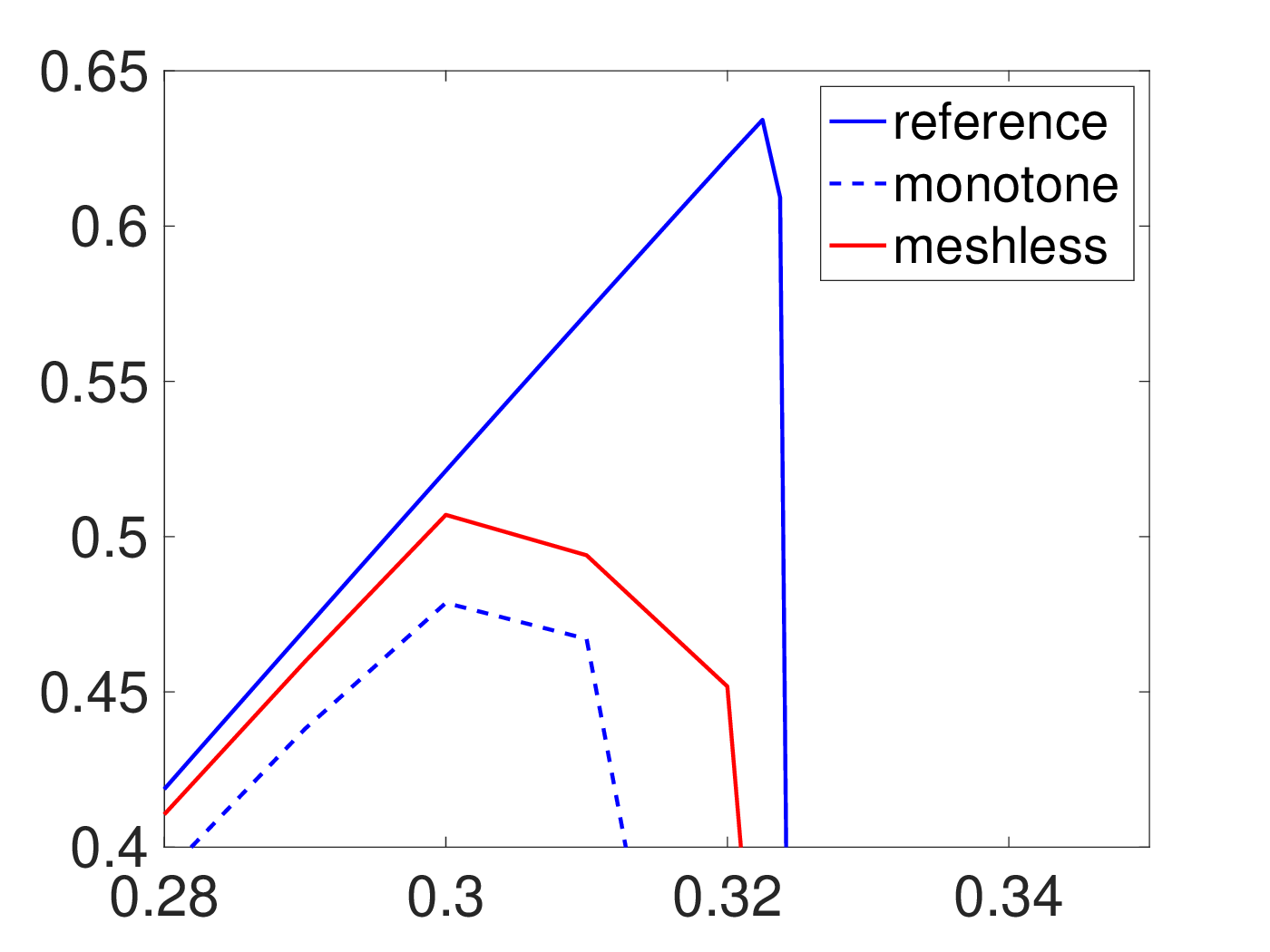}} 
\subfigure[Spacing 0.005.]{\includegraphics[scale=.17]{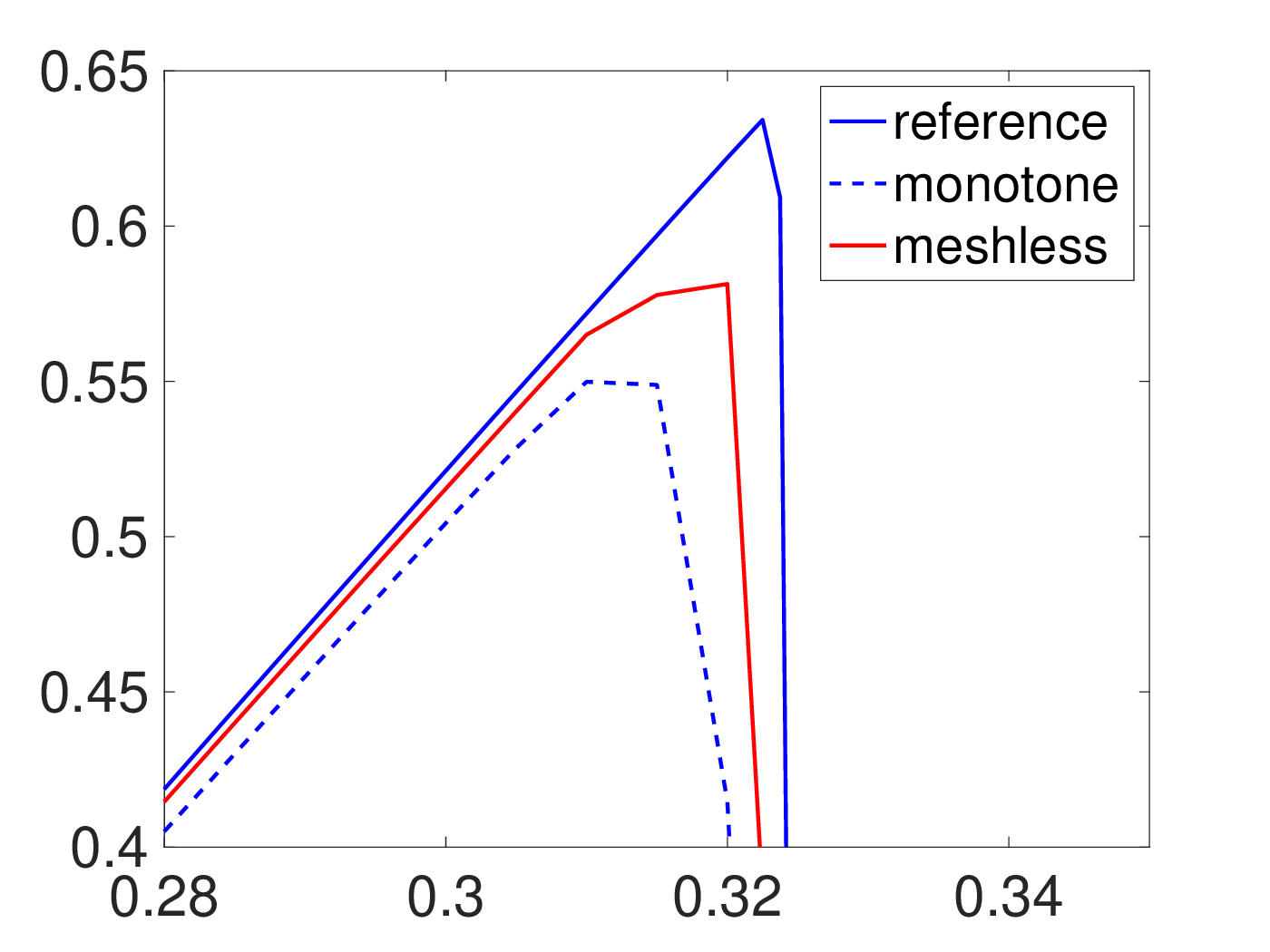}} 
\subfigure[Spacing 0.0025.]{\includegraphics[scale=.17]{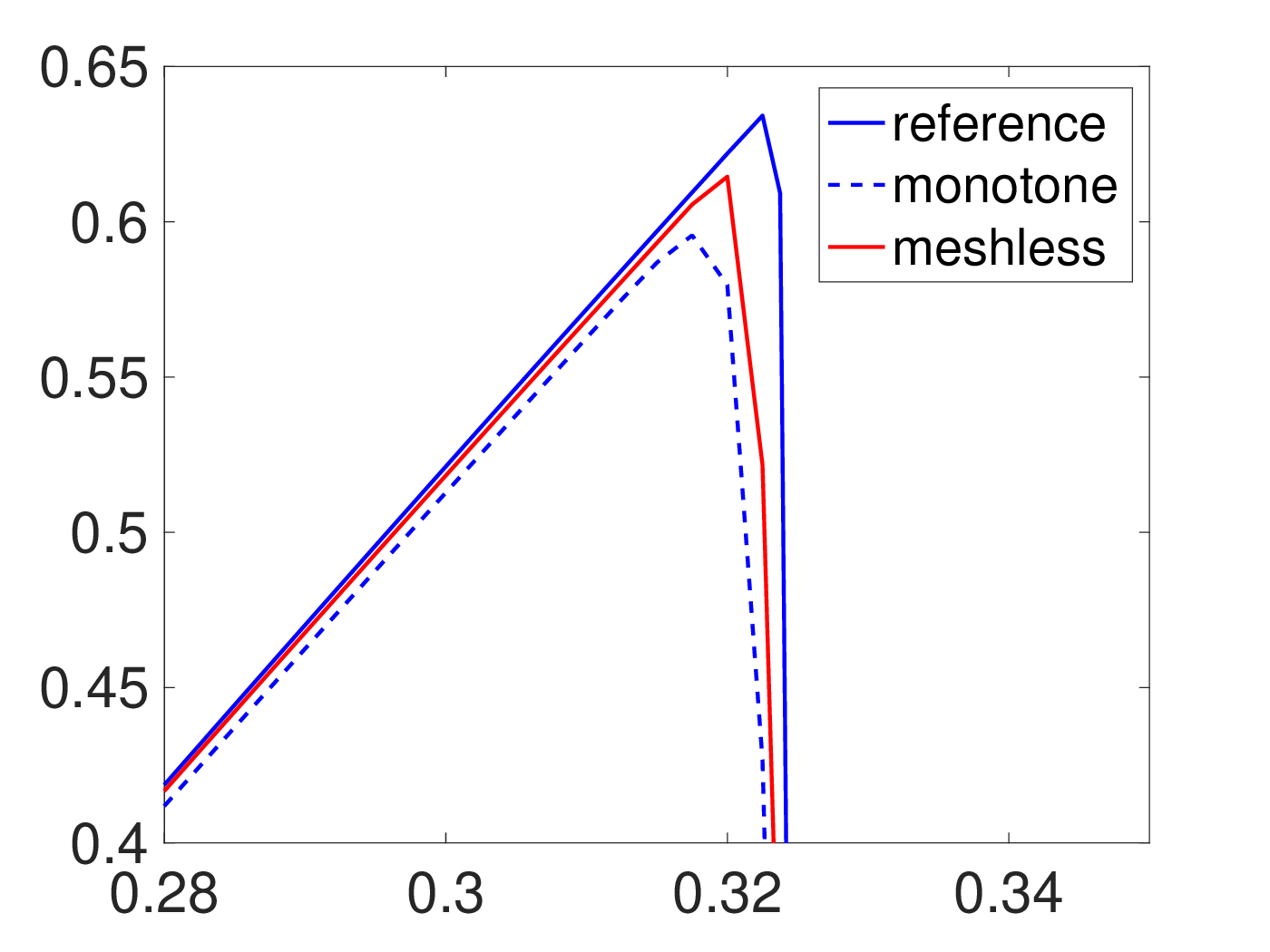}}
\subfigure[Spacing 0.00125.]{\includegraphics[scale=.17]{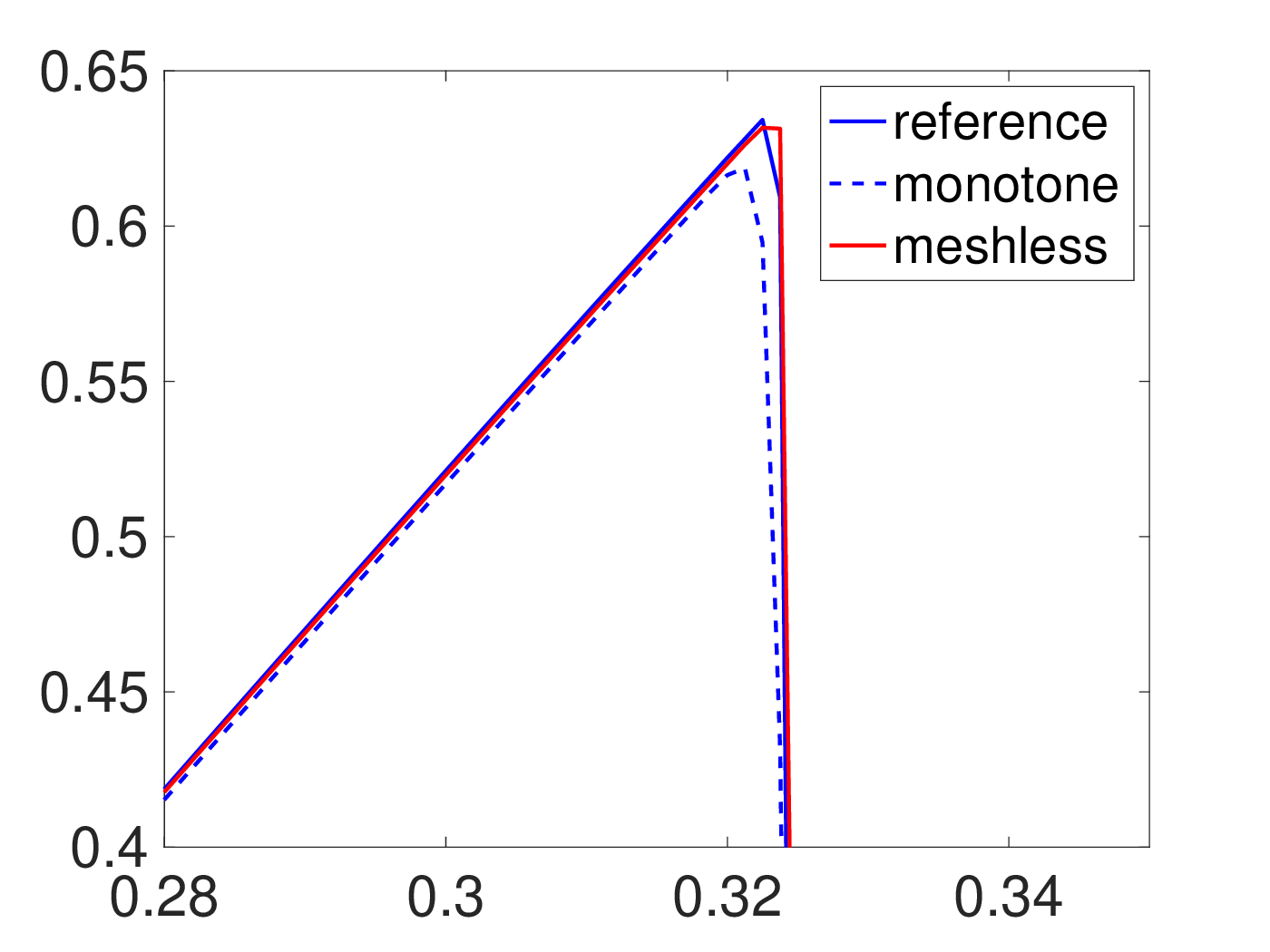}} 
\caption{Example 2: Comparison of the cross sections at $x_2= 0.1$ of 
the reference solution and  the monotone scheme with  the meshless scheme of  Algorithm~\ref{adaptnu} on 
 Halton points ($T=0.1$) with decreasing spacing in the vicinity
of the same corner as in the plots (e), (i) and (m) of Figure~\ref{fig:ex1cross}. } 
\label{fig:ex1cross_ref} 
\end{figure}

\blue{Finally, Figure~\ref{fig:ex1cross_ref} illustrates the performance of our method under refinement. In addition to the 
Halton points of the previous tests with spacing $h=0.0025$, we generate three more sets of Halton nodes with 
$h=0.01,\ 0.005$ and $0.00125,$ and compute solutions by Algorithm~\ref{adaptnu}. The four plots in Figure~\ref{fig:ex1cross_ref}
show the cross-sections of these solutions along with the reference solution and that by  the monotone scheme in the vicinity
of the same corner as in the plots (e), (i) and (m) of Figure~\ref{fig:ex1cross}. These results demonstrate a robust
convergence of the solution towards the sharp shocks.}

\begin{table}[!h]
\centerline{ 
\begin{tabular}{c|cc} 
\midrule
{method}& $E_1$-error & $E_2$-error\\
\midrule
monotone scheme (\cite{hestbook}, Script 7.2) & 1.21e-02&  4.38e-02\\
meshless scheme (Algorithm~\ref{adaptnu}), Cartesian grid& 6.08e-03& 2.79e-02\\
meshless scheme (Algorithm~\ref{adaptnu}), Halton points & 5.47e-03 & 2.71e-02\\
meshless scheme (Algorithm~\ref{adaptnu}), Random points & 6.45e-03 & 3.34e-02\\
\midrule
\end{tabular}
}
\caption{Example 2: Comparison of the errors (with respect to the reference solution) obtained by a standard monotone scheme
on the Cartesian grid, and our meshless method on the Cartesian grid, Halton and random nodes ($T=0.1$).}
\label{table:tabex1}
\end{table}

\subsection{Example 3: Rotating wave}\label{ex3}
\begin{figure}[!htbp]
\subfigure[Reference solution]{\includegraphics[scale=.55]{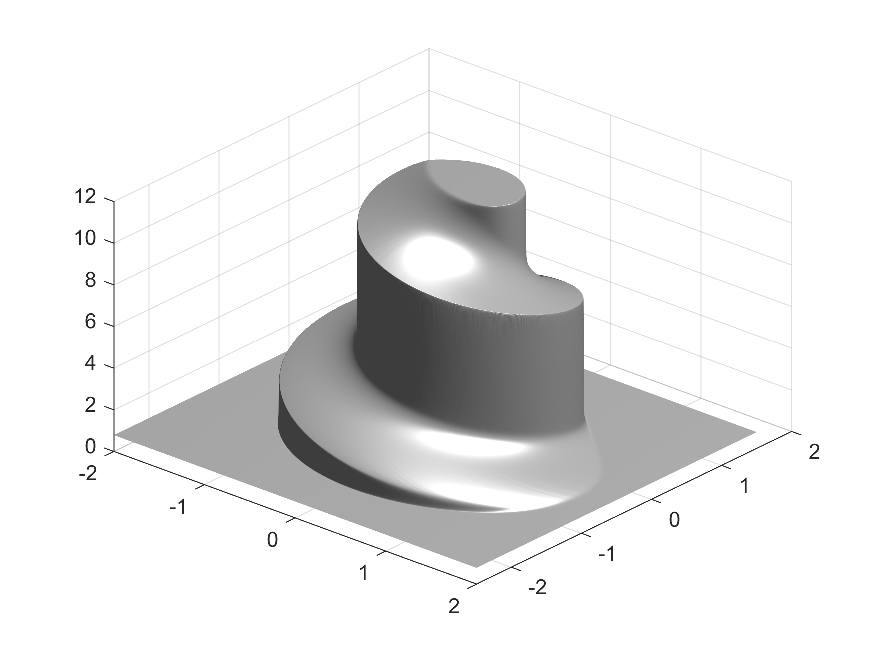}
\includegraphics[scale=.55]{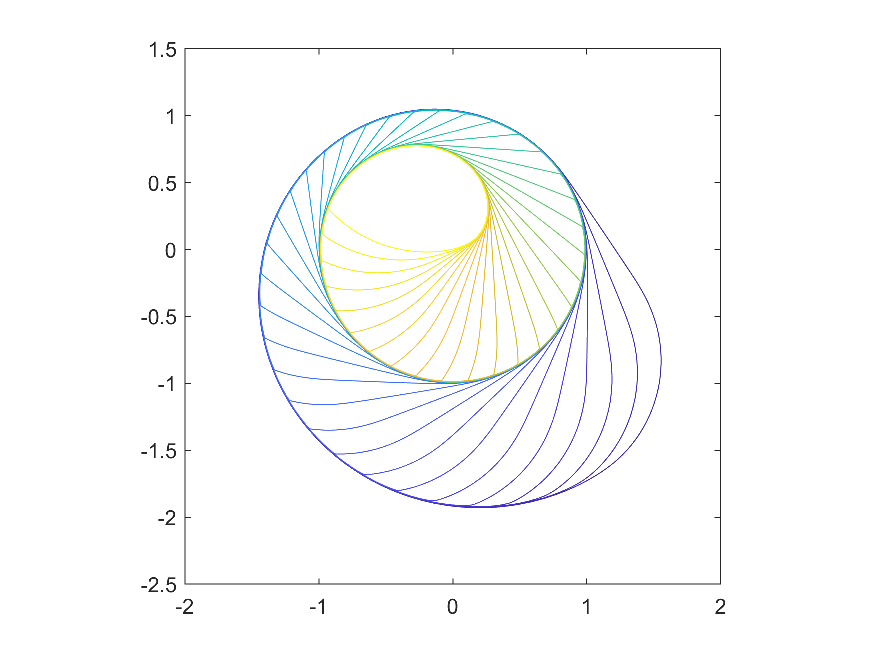}}
\subfigure[Monotone scheme (\cite{hestbook}, Script 7.2) on Cartesian grid]{\includegraphics[scale=.55]{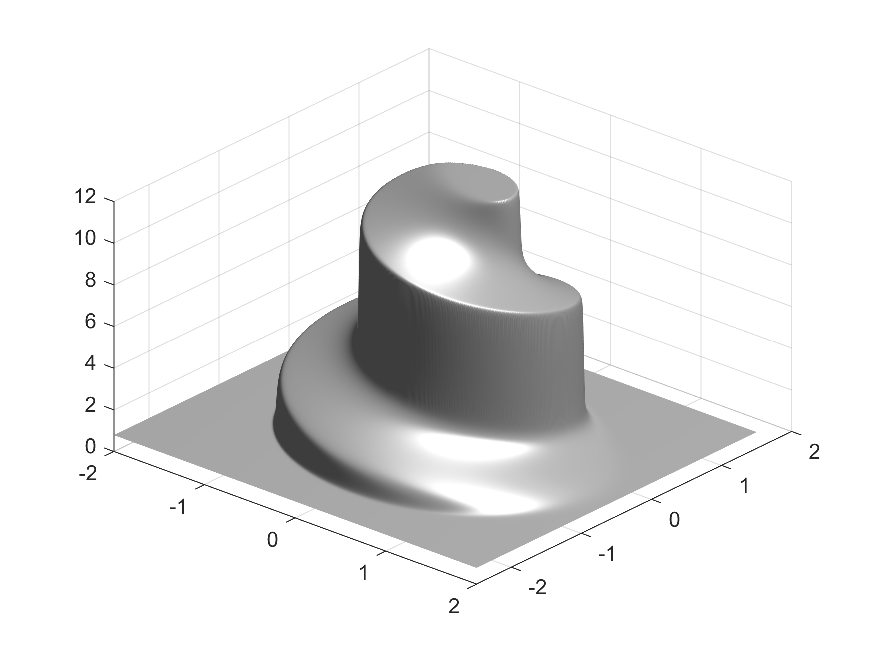}
\includegraphics[scale=.55]{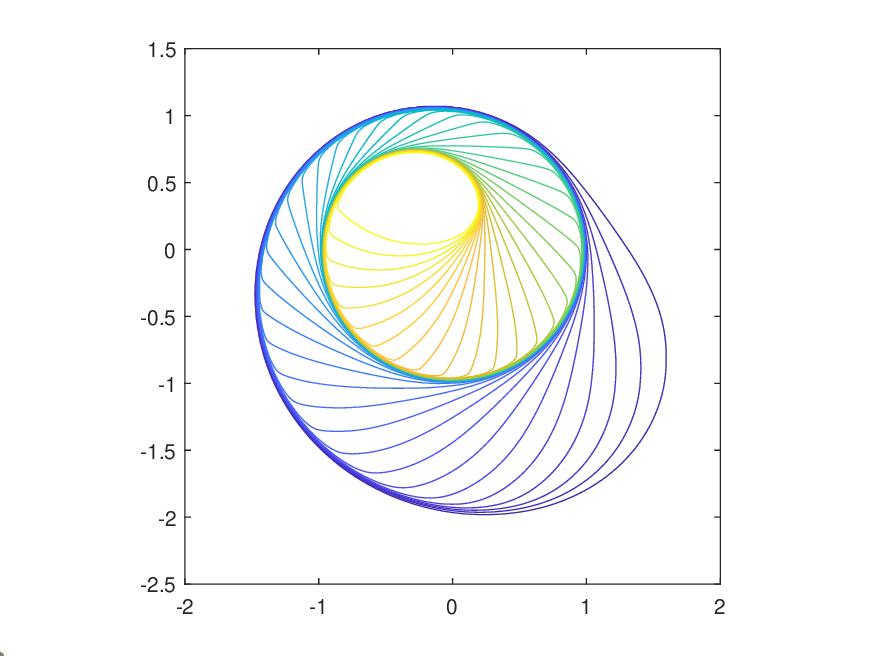}}
\subfigure[Meshless scheme (Algorithm~\ref{adaptnu}) on Halton points]{\includegraphics[scale=.55]{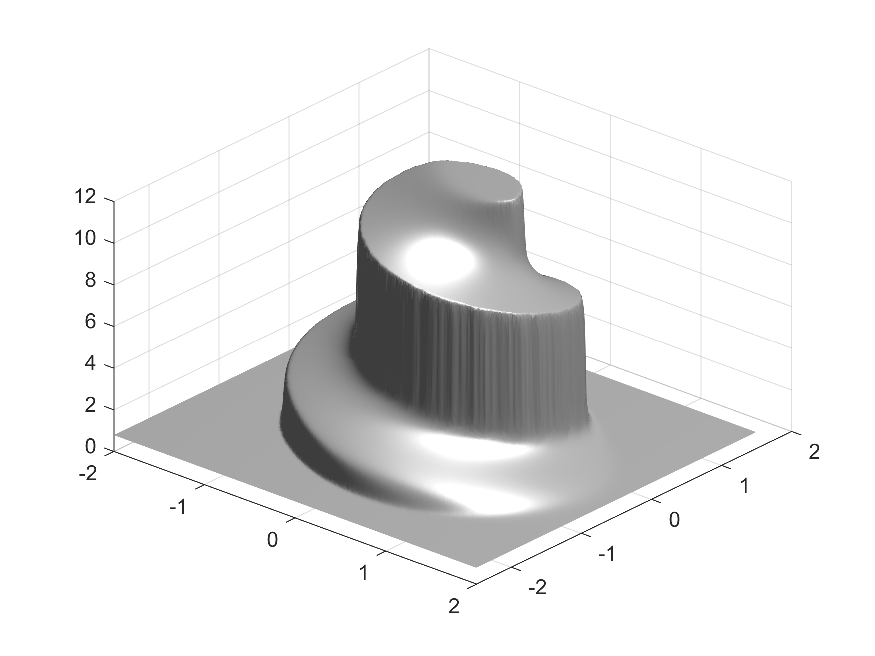}
\includegraphics[scale=.55]{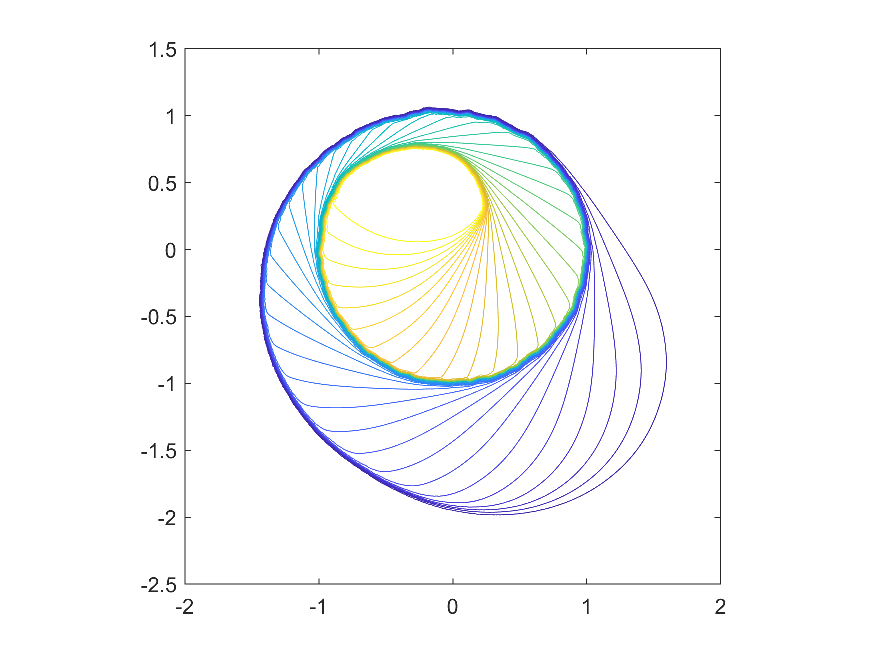}}
\caption{Example 3: (a) High-resolution reference solution  at $T=1$ and its contour lines. 
Comparison between the solutions at $T=1$ by (b) a monotone scheme on the $400\times400$ Cartesian grid
 and (c) our positive scheme with adaptive viscosity on $1.602\cdot 10^5$ irregular points.}
\label{fig:RW_comp}  
\end{figure}

As last example, we consider a more challenging equation with a non-convex flux ${\bf F}(u)=(\sin(u), \cos(u))$ on the spatial
domain $[-2, 2]\times [-2.5, 1.5]$, with periodic boundary conditions and initial condition
$$ 
u_0(\x) = \left\{ \begin{array}{ll} 3.5 \pi, & \mbox{if } \Vert \x \Vert_2 < 1, \cr 
0.25 \pi, & \mbox{otherwise}.   \end{array} \right. 
$$
The solution is clearly discontinuous from the beginning, but the shape of the discontinuity curve greatly changes as time
advances, making it an excellent test for our adaptive viscosity method. The exact solution is not available, but its
behavior is well known in the literature, see \cite{KPP07,guermond11,guermond17,hestbook}. This test problem has been found
challenging for many standard numerical schemes.
As a reference, here we use a high-resolution numerical solution obtained by employing the second order central scheme 
described in
\cite{hestbook}, Chapter 10, and the supplementary code \texttt{BurgersC2D.m}, on the $1600\times1600$ Cartesian grid, see
 Figure \ref{fig:RW_comp}(a).

\begin{figure}[!t]
\subfigure[Fault points (in red).]{\includegraphics[scale=.665, trim={1.8cm 0cm 0.8cm 0cm}]{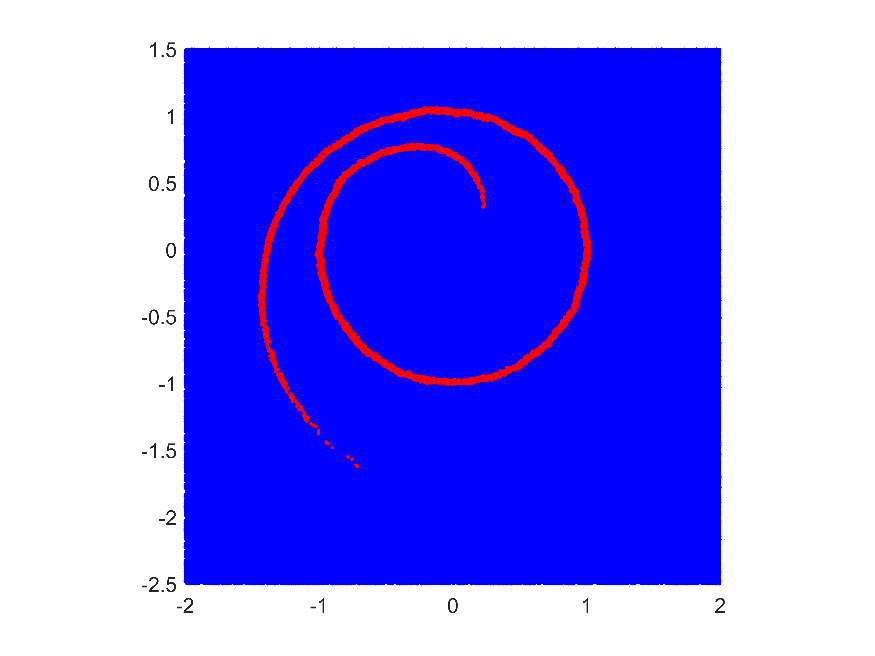}}
\subfigure[Values of $\mu_i$ obtained by formula \eqref{nui}.]{\includegraphics[scale=.665, trim={1cm 0cm 2cm 0cm}]{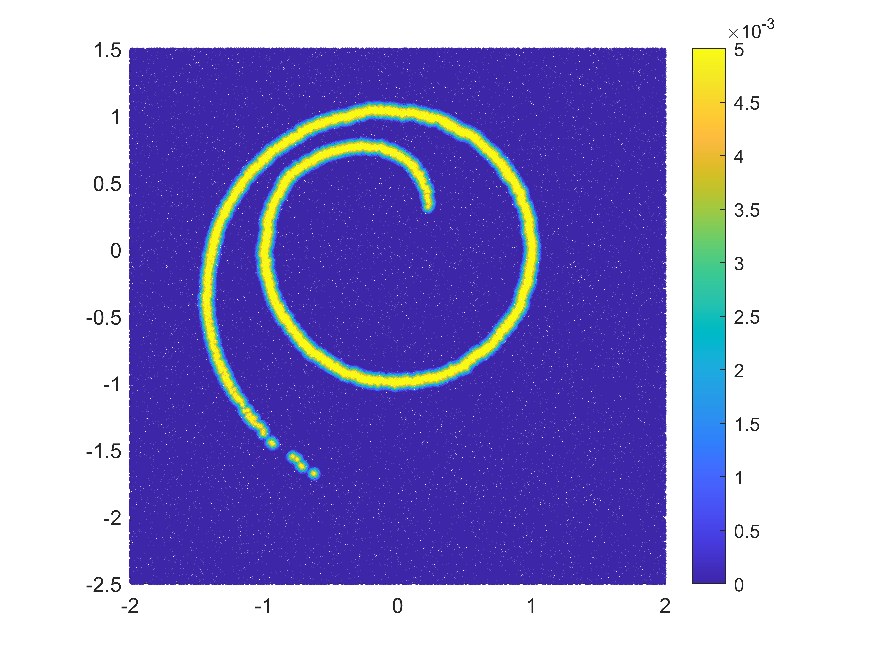}}
\caption{Example 3: Detected fault points and viscosity values at $T=1$.}
\label{fig:RW_fault}     
\end{figure}

We apply our positive scheme with adaptive viscosity (Algorithm~\ref{adaptnu}) on irregular nodes
obtained in the same way as in Example~2 by modifying $400^2$ Halton points in $[-2, 2]\times [-2.5, 1.5]$ near 
the boundaries. This corresponds to the spacing $h=0.01$, and hence we choose
$\Delta t=0.2h=2\cdot 10^{-3}$, $\mu = 0.5h=5\cdot 10^{-3}$.
The results at $T=1$ are again compared with the solution obtained by employing the first order monotone scheme of 
\cite{hestbook}, Script 7.2, on the $400\times400$ Cartesian grid.
 Figure \ref{fig:RW_comp}(b,c) compares the solutions at the final time 
and the corresponding contour lines. The results reproduce the correct shape, with our solution
having edges of discontinuity a bit less neat, which is expected having used Halton points. 
Despite this, Table \ref{table:tabex4} shows that the error with respect to the reference solution obtained with our method is somewhat better. 
\blue{Note that the maximum of the numerical solution is preserved at the same value $3.5\pi$ 
up to machine precision at all time steps.} Finally,
Figure \ref{fig:RW_fault} shows the detected fault points \blue{and the distribution of the viscosity 
values at the final time, indicating that viscosity is applied near the shocks as intended.}

\begin{table}[!h]
\centerline{ 
\begin{tabular}{c|cc} 
\midrule
{method}& $E_1$-error & $E_2$-error\\
\midrule
monotone scheme (\cite{hestbook}, Script 7.2) & 7.38e-03&  3.76e-02\\
meshless scheme (Algorithm~\ref{adaptnu}), Halton points & 6.51e-03&  2.93e-02\\
\midrule
\end{tabular}
}
\caption{Example 3: Comparison of the errors (with respect to the reference solution) obtained by a standard monotone scheme
on the Cartesian grid, and our meshless method on  Halton points ($T=1$).}
\label{table:tabex4}
\end{table}

